\input amstex
%This will run under Ams-Tex (without style files

%The following (from Ams-Tex) need to be defined 
%or redefined for plain Tex

%\plainfootnote   \,   \frac   
%\text   \pmb (= Poor man's bold)   \Cal (calligraphic font)
% \comment  \endcomment   \cases \endcases
% \matrix \endmatrix

%\font\Ti=Times
%\def\shat{{\Ti š}}  \font\eesp=ElysiumBook at 11.3pt

%this defines eth!

%%%%%%font definitions

\font\rm=cmr10 \rm

\font\bf=cmb10
\font\Rm=cmr9 at 11pt
\rm
\font\it=cmsl9 at 10pt
\font\sc=cmr7 %[track+70] 
 at 7pt
\def\Sc #1{{\sc \uppercase{#1}}}
\font\Rrm=cmr17 at 16pt
   \font\Rm=cmr12 at 11.5pt
%%%%%%%%%%
%\font\Frak=eufm10

  %%% (P,x^n)
\long\def\Pf{\par\noindent {\it Proof.} }
\def\({\left(}
\def\){\right)}
\def\st{such that }
\def\qed{\hfill$\bullet$\vskip 4pt}

\def\brcs#1{\left\{ #1\right\}}

\def\iso{\cong}
\def\wrt{with respect to }
\def\:{\,:}

\def\rk{\text{rank\,}}
\def\det{\text{det\,}}

\def\ker{\text{ker\,}}

\def\I{\text{I\,}}

\def\R{\text{\bf R}}
\def\N{\text{\bf N}}
\def\Z{\text{\bf Z}}
\def\Q{\text{\bf Q}}

\def\Arrow #1;#2.{#1\:#2 \to }
%#1: #2 -> 

\def\Set#1#2{\brcs{#1 \left|\vphantom{#1 #2} \right.#2}}
%set of #1 such that #2

\def\Oh#1{{\pmb O}\(#1\)}
%%%bigOh notation 

\def\oh#1{{\pmb o}\(#1\)}
%%%little oh notation

%product left and right
\def\Rrr#1,#2{{\Cal J}_{#1,#2}}%difference
\def\slfrac#1#2{{\raise -.07 ex\hbox{$^{#1}$}}\!/\raise .35 ex \hbox{${}_{#2}$}}
\def\ssf #1/#2{\slfrac {#1}{#2}}

\def\pd #1,#2.{\frac {\partial #1}{\partial #2}}

%%%%%%%statement of lemmas, propositions,
%%theorems, examples, ...  \Lem 
   \long\def\Lem
#1.#2\par{\vskip4pt{\baselineskip=13pt\font\it=cmsl12 at
11.5pt\Rm
   \noindent {\rm \uppercase{#1}} #2\vskip3pt

   }} 

\long\def\Proclaim #1.#2 \endproclaim{\vskip4pt{\baselineskip=13pt\font\it=cmsl12 at
11.5pt\Rm
   \noindent {\rm \uppercase{#1}} #2\vskip3pt

   }} 

\long\def\remark #1\endremark{\vskip 2pt \noindent {\it Remark\/} #1\par}

\long\def\Sectionhead #1.#2:\par #3{\vskip 4pt \noindent {\bf #1 #2}vskip 2pt\noindent\nospace #3}

\long\def\Title #1\par {\noindent{\Rrm #1}\vskip 9pt}
%main title

 \long\def\SubT #1.{\noindent {\it #1\/} } 
%subtitle
 
 \long\def\SecT
#1\par{\vskip 3pt \noindent {\bf #1}\vglue1pt
   \noindent}%section title

\long\def\subtitle #1.{\vskip 2pt \noindent {\it #1}}

\long\def\Rmk#1\par{\vskip 1pt \noindent {\it
Remark.} #1\vskip2pt}%remark

\long\def\Abstract #1\par{{\leftskip= 3 true cm \rightskip = 3 true cm \font\it=cmsl10 \font\rm=cmr10 \baselineskip = 10pt
\parindent=.35 true cm\rm\noindent 
{\it Abstract} #1\vskip 8pt

}}

\long\def\Author #1 \par{\noindent{\it #1}}
%%%%%%%%%macros for numbering
%%of propositions, lemmas,examples, theorems

\scrollmode\NoBlackBoxes
\magnification=1100
%%%%%%font definitions

\font\rm=cmr10 \rm

\font\bf=cmb10
\font\Rm=cmr9 at 11pt
\rm
\font\it=cmsl9 at 10pt
\font\sc=cmr7 at 7pt %[track+70] 
 
\def\Sc #1{{\sc \uppercase{#1}}}
\font\Rrm=cmr17 at 16pt
   \font\Rm=cmr12 at 11.5pt

\def\Sc #1{{\sc \uppercase{#1}}}

  %%% (P,x^n)
\long\def\Pf{\par\noindent {\it Proof.} }
\def\({\left(}
\def\){\right)}
\def\st{such that }
\def\qed{\hfill$\bullet$\vskip 4pt}

\def\brcs#1{\left\{ #1\right\}}
\def\Set#1#2{\brcs{#1 \left|\vphantom{#1 #2} \right.#2}}

\def\wrt{with respect to }
%\long\def\Lem #1 #2\par{\noindent {\Sc lemma #1.} {\Rm
%#2}\vskip 2pt}
\long\def\Lemma #1. #2\par{\noindent {\Sc  {#1.}} {\Rm
#2}\vskip 2pt}
\def\Arrow #1;#2.{#1\:#2 \to }

\def\Oh#1{{\pmb O}\(#1\)}
%%%bigOh notation 
\def\oh#1{{\pmb o}\(#1\)}
%%%little oh notation
\def\R{\text{\bf R}}
\def\N{\text{\bf N}}
\def\Z{\text{\bf Z}}
\def\Q{\text{\bf Q}}
 
\def\slfrac#1#2{{\raise -.07 ex\hbox{$^{#1}$}}\!/\raise .35 ex \hbox{${}_{#2}$}}
\def\ssf #1/#2{\slfrac {#1}{#2}}

    %notation: R_P, (P,x), pure traces etc
    %\C complexes \R reals \quotes#1
    %P will be of rad of cvg =1 with no negative coef 
%subordinate to: subeqv; f subeqv P means 
% all coefficients of P are nonnegative, and the 
%absolute values of coef of f are bdd above by etc.
    %equivalent to: eqv
    %explanation of 0 \leq f

\font\Rrm=cmr17 at 16pt
   \font\Rm=cmr12 at 11.5pt
   \long\def\Lem
#1.#2\par{\write1{#1,
p\,\folio\par}\vskip4pt{\baselineskip=13pt\font\it=cmsl12 \Rm
   \noindent {\rm \uppercase{#1}} #2\vskip3pt

   }} 

   \long\def\Title #1\par {\noindent{\Rrm #1}\vskip 9pt}
   \long\def\SubT #1.{\noindent {\it #1\/} }
   \long\def\SecT #1\par{\vskip 4pt \noindent {\bf #1}\vglue1pt
   \noindent}%section title

%%%%%%%%%macros for numbering

\def\oneone{1.1}

\def\twoone{2.1}

\def\throne{3.1}

\def\fouone{4.1}
\def\foutwo{4.2}
\def\fouthr{4.3}
\def\foufou{4.4}

\def\fivone{5.1}
\def\fivtwo{5.2}
\def\fivthr{5.3}
\def\fivfou{5.4}

\def\sixone{6.1}

\def\sevone{7.1}
\def\sevtwo{7.2}

\def\ninone{9.1}

\input diagrams

\def\flo #1{\lfloor #1 \rfloor}

\def\diag{\text{{\rm diag}}\,}
\def\rank{\text{{\rm rank}}\,}
\def\Q{\text{\bf Q}}

\def\ase{algebraic shift equivalence{}}
\def\Endc{\text{End}_c}
\def\GL{\text{GL}\,}
 
\def\ecrs{ECRS}

\let\hat=\widehat

\def\Aff{\text{Aff\,}}
\NoBlackBoxes

\def\One{1}
\def\oneone{1.1}

\def\twoone{2.1}

\def\throne{3.1}

\def\fouone{3.1}
\def\foutwo{3.2}
\def\fouthr{3.3}
\def\foufou{3.4}

\def\fivone{8.2}
\def\fivtwo{10.1}
\def\fivthr{8.3}
\def\fivfiv{5.3}
\def\fivsix{8.1}
\def\fivsev{10.4}
\def\fiveig{10.3}

\def\sixone{4.1}
\def\sixtwo{6.1}
\def\sixthr{6.2}
\def\sixfou{6.3}
\def\sixfiv{6.4}

\def\sevone{5.1}
\def\sevtwo{7.1}

\def\eigthr{8.3}
\def\eigfou{8.4}
\def\eigfiv{8.5}

\def\ninone{9.1}
\def\nintwo{9.2}
\def\ninthr{9.3}

\def\ninele{9.11}

\def\tenone{10.1}
\def\tentwo{10.2}

\def\eleone{11.1}
\def\eletwo{11.2}
\def\elethr{11.3}
\def\elefou{11.5}
\def\elefiv{11.4}
\def\elesix{11.6}
\def\elesev{11.7}
\def\eleeig{11.8}
\def\elenin{11.9}
\def\eleten{11.10}
\def\eleele{11.12}
\def\eletwe{11.13}
\def\elethi{11.11}

\def\I{\text{I\,}}

\def\flo #1{\lfloor #1 \rfloor}

\let\iso=\simeq
\def\ceil#1{\lceil #1 \rceil}

\def\tripnorm #1xxx{\left\|\hglue-.2ex\left|#1\right|\hglue-.2ex\right\|}

\def\flo #1{\lfloor #1 \rfloor}

\def\diag{\text{diag}\,}

\def\ecs{ECS}\def\ers{ERS}
\def\One{{\pmb 1}}

\Title Equal column sum and equal row sum dimension group realizations

\Abstract  Motivated by connections between minimal actions, especially Tšplitz, on Cantor sets, and dimension groups, we find realizations of classes of dimension groups as limits of primitive matrices all of which have equal column sums, or equal row sums.  

\Author David Handelman%
\plainfootnote{$^{1}$}{\rm Supported in part by an NSERC Discovery Grant.}

\noindent All groups are abelian, {\it free\/} means free as an abelian group, all partially
ordered groups are unperforated and torsion-free. Equivalence classes representing elements of the direct limit, $\lim \Arrow A_n; F_n.F_{n+1}$ are expressed as $[v,n]$, where $v$ belongs to the $n$th free abelian group.

Suppose $U$ is a noncyclic subgroup of the rationals, and let $\Arrow \tau;G.U \subset \R$ be an onto group homomorphism from a  torsion-free
group $G$ to $U$. We may impose an ordered group structure on $G$
 simply by declaring $g \in G^+\setminus \brcs{0}$ iff $\tau(g) > 0$. This makes $G $ into a simple dimension group with unique trace, and the trace is rational-valued; all such simple dimension groups are constructed in this manner. That is, $G$ is an extension (in the category of abelian groups) of a torsion-free abelian group $\ker \tau$ by the subgroup $U$ of the rationals. 

As $G$ is a countable dimension group, by [EHS], it has a representation as {\it ordered groups,} $G \iso \lim \Arrow A_n; \Z^{f(n)}. \Z^{f(n+1)}$ where $\Arrow f;\N.\N$ is a function, we take the usual coordinatewise ordering on each $Z^{f(n)}$ and impose the usual direct limit ordering. The $A_n$ have only nonnegative entries  (and, by telescoping, can be made strictly positive when $G$ is simple).
However, [EHS] does not give specific representations, that is, the matrices $A_n$ cannot be  constructed from the argument, except by extremely complicated machinations. Here we consider the case that $G$ be  of rank $k+1$ (so $\ker \tau$ is rank $k$), and provide explicit realizations for $G$ with the ECS property ({\it equal column sums\/}: each of the nonnegative (or strictly positive) matrices $A_n $ has all of its column sums equal).

With an \ecs\ realization, there is a canonical choice of trace, namely (up to scalar multiple), the sequence of normalized multiples of constant rows; in this case, we say the trace has an \ecs\ realization. We  show that for general dimension groups with order unit, a trace admits an \ecs\ realization iff it is faithful, rational-valued, and good (in the sense of Akin, after translation to dimension groups as in [BeH]). In this case, there is no control on the matrix sizes, but we do not require simplicity.

When we take the transposes of the matrices used for \ecs\ realizations, and thus obtain \ers---equal row sum---realizations, the resulting dimension groups run over all possible simple dimension groups of finite rank with unique trace (not generally rational-valued) which could admit an \ers\ realization. These are very closely related to Tšplitz systems (pairs $(X,T)$ consisting of a Cantor set and a minimal self-homeomorphism, which is an almost everywhere one-to-one extension of an odometer), as explained to  me by Chris Skau, whose question about \ers\ realizations motivated this paper.

An \ers\ realization of a simple dimension group $G$ \wrt a (noncyclic rank one subgroup) $H$ \st $\tau(H) \neq
0$ and $G/H$ is torsion-free is an ordered group isomorphism $\Arrow \phi; G.\lim \Arrow A_n; \Z^{f(n)}. \Z^{f(n+1)}$ where $\Arrow f; \N.\N$ is some function, $A_n$ are nonnegative integer matrices each having equal row sums,  the direct limit ordering is imposed, \st $\phi(H) = \cup_n [\pmb 1_{f(n)},n]\Z$, with $\pmb 1_{f(n)}$ the column consisting of ones. If $f(n) = s$ for all $n$, the realization is of size $s$. 

If $G$ and $H$ are as in the prevous paragraph, and there exists an \ers\ realization of $G$ \wrt $H$ that is also \ecs, then we refer to this as an \ecrs\ realization of $G$ \wrt $H$. 

We establish following results on \ecs, \ers, and \ecrs\ realizations.

\Lem . Let $G$ be a dimension group.
\item{(i)} If $G$ is simple, of rank $k+1$ with unique trace $\tau$, and $\tau(G)$ is a subgroup of the rationals, then there exists an \ecs\ realization of $G$ of size $k+2$. (Theorem \sixone)
\item{(ii)} Let $\tau$ be a trace on $G$ with $\tau(G) \subseteq \Q$. Then there exists an \ecs\ realization of $G$ representing $\tau$ if and only if $\tau$ is good (in the sense of Akin, as translated to the dimension group setting [BeH]) and faithful (that is, $\ker \tau \cap G^+ = \brcs{0}$. (Theorem \sixtwo(b))
\item{(iii)} Suppose $G$ is simple and has  unique trace $\tau$, and $H$ is a noncyclic rank one subgroup of $G$ \st $\tau(H) \neq 0$ and $G/H$ is torsion-free. 
\itemitem{(a)} If $\rk G = k+1$, then there exists a size $k+2$ \ers\ realization of $G$ \wrt $H$. (Theorem \sevtwo(a))
\itemitem{(b)} There exists an \ers\ realization of $G$ \wrt $H$. (Theorem \sevtwo(b))
\item{(iv)} Suppose $G$ is as in (iii), and in addition, $\tau(G) \subseteq \Q$.
\itemitem{(a)} If $\tau(G)$ has no primes of infinite multiplicity (that is, $\tau(G)$ is not $p$-divisible for any prime $p$), then $G$ admits an \ecrs\ realization \wrt $H$ if and only if $\lambda:=|\tau(G)/\tau(H)| \geq \rk G$; when  $\lambda < \infty$, there is an \ecrs\ realization of size $\lambda$ but none of other sizes. (Theorem \elenin)
\itemitem{(b)} If $\tau(G) $ has a prime of infinite multiplicity (that is, $\tau(G)$ is $p$-divisible for some prime $p$), then $G$ admits an \ecrs\ realization \wrt $H$; this can be constructed to be bounded if $|\tau(G)/\tau(H)| < \infty$. (Theorem \eletwe)   

Part (ii) above applies to all dimension groups (with an order unit), but the other parts require simplicity and unique trace. 

Much of the time, we work in the category of abelian groups with group homomorphism to the
reals:  a torsion-free abelian group $G$ together with a group
homomorphism $\Arrow t; G.\R$ \st $t(G)$ is dense in $\R$; we denote this
$(G,t)$. Isomorphism in this category is the obvious one, $(G,t) \iso
(G',t')$ if there exists a group isomorphism $\Arrow \phi; G.G'$ \st
$t'\phi$ is a nonzero scalar multiple of $\tau$. Automatically, this
induces an isomorphism $\ker t \iso \ker t'$. 
 
Suppose that $G$ and $G'$ are noncyclic simple dimension groups with unique trace,
$\tau$ and $\tau'$ respectively. Then $G \iso G'$ as ordered groups if
and only if $(G,\tau) \iso (G',\tau')$ as abelian groups with real-valued
group homomorphism. One way is trivial. Conversely, suppose $\Arrow \phi;
G.G'$ is a group isomorphism
\st $\tau' \phi = \lambda \tau$ for some nonzero real $\lambda$. By
replacing $\phi$ by $-\phi$ if necessary, we may assume $\lambda > 0$.
Then $\phi$ is an isomorphism of ordered groups.
 
To see this, we note that $g \in G^+ \setminus\brcs{0} $ iff $\tau(g) >
0$; this is equivalent to $\tau' (\phi(g)) > 0$, which is equivalent to
$\phi(g) \in ( G')^{+} \setminus \brcs{0}$. As $\phi$ is a group
isomorphism, this says both $\phi$ and $\phi^{-1}$ are order preserving,
hence $\phi$ is an order isomorphism.
 
Hence to decide if two simple dimension groups with unique trace are order
isomorphic, it is sufficient to find a group isomorphism between them that
scales the trace(s). This makes life simple, at least when the dimension
group has unique trace.
 
The dimension groups we will be considering for \ecs\ realizations have an additional
property, that the range of their trace is (up to nonzero scalar multiple)
a subgroup of the rationals. So we consider them as groups with
real-valued group homomorphism, $(G,t)$ \st $t(G) $ is rank one and dense
in $\R$ (so up to scalar multiple, $t(G) = U \subseteq \Q$).
 
Although we will often be talking about extensions of abelian groups, $0 \to C \to G \to
U \to 0$, it is too restrictive to deal with the classification as
extensions (that is, within $\text{Ext}^1 (C, U)$); instead, we are
dealing with the coarser classification, isomorphism for maps  $G \to U$,
where we are allowed to multiply by $\pm 1$ (and if $U$ is $p$-divisible,
by powers of $p$). There are still generically uncountably many isomorphism classes of
these, since $\text{Aut} (C)$ and $\text{Aut} (U)$ are usually countable.

As usual, if a group or ordered group is given as a direct limit, $\lim \Arrow M_n;F_n.F_{n+1} $ (typically, $F_n$ will be free abelian groups, and if the ordered direct limit is required, the entries of $M_n$ will be nonnegative), then elements of the direct limit can be written as equivalence classes, $[a,n]$ where $a \in F_n$, and the equivalence relation is generated by $[a,n] = [M_n a,n+1]$. 

\SecT 1 Via subsemigroups

For this section, $G$ need only  be a partially ordered  group with positive cone $G^+$. Let
$P$ denote the set of nonnegative integers. If $\brcs{a_i} \subseteq G^+$,
we denote by $\sum a_i P$, the set of sums $\Set{\sum a_i n(i)}{n(i) \in
P}$, the semigroup (or subsemigroup) generated by $\brcs{a_i}$.
 
Let $\brcs{S_n}_{n\in \N}$ be a collection of subsemigroups of $G^+$ with
$S_1 \subseteq S_2 \subseteq S_3 \subseteq \dots$ \st $G^+ = \cup S_n$. Suppose
$S_n$ is generated by $\brcs{a_i^{(n)}}_{i=1}^{f(n)}$. Since $S_n
\subseteq S_{n+1}$, we can find an $f(n+1)\times f(n)$ matrix $A_n$
(called a {\it transition\/} matrix)  with entries from $P$ \st for all
$i$,
$$
a_i^{(n)} = \sum_{j=1}^{f(n+1)} (A_n)_{ji}a_j^{(n+1)}; \tag*
$$
there is usually a great deal of choice available for the matrix entries,
since there is no assumption of any sort of unique decomposition. Note the subscript $ji$, not $ij$.
 
Let $F_n = \Z^{f(n)}$, the free abelian group on $f(n) $ generators
(denoted $e_i^{(n)}$, but when  superscripted ${}^{(n)}$ is understood,
it is removed), equipped with the usual coordinatewise ordering. Now form
the direct limit dimension group from the $A_n$s, $H = \lim \Arrow A_n;
F_n .F_{n+1} $.  Define $\Arrow \psi_n; F_n. G$ via $\psi_n (e_i) =
a_i^{(n)}$. This is a well-defined positive group homomorphism from $F_n $
to $G$. The condition in (*) is precisely what we need in order that
$\psi_{n+1}\circ A_n = \psi_n$. Hence the family $\brcs{\psi_n}$ induces a
positive group homomorphism $\Arrow \Psi; H.G$ (explicitly, $\Psi [v, n] =
\psi_n (v)$ where $v \in F_n$).
 
Since $G^+ = \cup S_n$, $\Psi(H^+) = G^+$; if we assume $G$ is directed
(as we may as well), that is, $G = G^+ - G^+$, then $\Psi$ is onto. If
$\Psi$ is one to one, then it is an isomorphism of ordered groups (in
particular, $G$ is a dimension group), and we have a realization for it
using the matrices $A_n$. If  $\rk H \leq \rk G < \infty$,
then $\Psi$ is automatically an isomorphism (since an onto homomorphism
from a torsion-free abelian group of finite rank to a torsion-free group
of the same or more rank is automatically one to one), then $\Psi$ is an
isomorphism.
 
The construction of $\Psi$ depends on the choice(s) of the generators for
the semigroups $S_n$, and then on the matrices $A_n$; different choices
for the matrices (even fixing the generators of all the $S_n$) can result
in different $\Psi$ functions, some of which may be one to one while
others need not be.
 
We summarize this in one gigantic statement.
 
\Lem Lemma \oneone. Suppose that $G$ is a partially ordered abelian group with an
increasing set of subsemigroups, $S_1 \subseteq S_2 \subset \dots$ \st
$G^+ = \cup S_n$, and suppose that $A_n$ is the transition matrix
associated to a choice of generators for  $S_n \subset S_{n+1}$, with each
$S_n $ generated by $\brcs{a_i^{(n)}}$. Form the dimension group $H = \lim
\Arrow A_n ; F_n. F_{n+1}$. {\par}
\noindent \item{(a)} There is a unique positive  group
homomorphism $\Arrow \Psi;H.G$ \st $[e_i^{(n)},n] \mapsto a_i^{(n)}$; moreover, 
$\Psi(H^+) = G^+$.
\item{(b)} If $\Psi$ is one to one, then it is an isomorphism of
ordered abelian groups, and thus $G$ is a dimension group. 
\item{(c)} If $G$ is
torsion-free \st $\rk H \leq \rk G $ and $\rk H < \infty$, then $\Psi$ is
an ordered group isomorphism.

\SecT 2 Realizing $G$ as \ecs\ (free kernel that splits)

Over this and the next few sections, we deal with the simple dimension group $G$ of rank $k+1$ having unique trace $\tau$, and in addition, $\tau(G) = U$ is a rank one (necessarily noncyclic) subgroup of the reals. For expository reasons, we proceed in three steps. 

This section deals with a rather special 
case, that $\ker \tau$ be free of rank $k$ and the extension splits. In the next section , we drop the splitting property (but maintain freeness of the kernel); finally, we deal with the general case, wherein $\ker \tau$ is an arbitrary rank $k$ torsion-free abelian group, and the extension by $U$ is arbitrary. We could   go straight to the general case, but this would have resulted in a very complicated argument. Instead, as we proceed through the cases, we find the extra complications can be dealt with in a relatively smooth manner. 

Here we deal with the easiest case, $G = U \oplus \Z^k$ where $G^+\setminus\brcs{0} = \Set{(u,w)}{u> 0}$.
Although we know that $G = U \oplus \Z^k$ is a dimension group, and $G^+\setminus \brcs{0} =
\Set{(u,w)}{u \in \Q^{++} \text{ and } w \in \Z^k}$, and thus is a limit
of free abelian groups with their coordinatewise limit by [EHS], the
latter does not give an explicit form. Here we obtain from a natural (but
not the most natural) subsemigroups of $G^+$, an explicit realization
with all the matrices being size $k+2$ and column stochastic (all column
sums equal for each matrix; this is abbreviated \ecs).  The following is the result of this section.

\Lem Proposition \twoone. Let $G = U \oplus \Z^k$ where $U$ is a noncyclic subgroup of the rationals, and $G$ is the simple dimension group obtained from the map $G \to U$. Then $G$ can be realized as a direct limit (in the category of ordered abelian groups) $\lim \Arrow A_n;\Z^{k+2}. \Z^{k+2}$ where $A_n$ are primitive and \ecs. 
 
First, we find a suitable representation of $\Z^k$ as a union of $k+1$
subsemigroups. For $1 \leq i \leq k$, let $\epsilon_i$ denote
the standard basis vector of $\Z^k$, and set $\epsilon_{k+1} = -\sum
\epsilon_i$. Obviously $\sum_{i=1}^{k+1} \epsilon_i = \pmb 0$ and it is
easy to verify that $\sum_{j=1}^{k+1} \epsilon_i P = \Z^k$.
 
Now let the supernatural number of $U$ be given. We may block (telescope) all the
primes and their powers that appear, so that we have a sequence of
positive integers $p_1, p_2, \dots,$ with $p_n> (k+1)^2$ for all $n$  and
$U \iso \lim \Arrow \times p_i ;\Z.\Z$. Let $q_n = \prod_{i=1}^n p_i$. Now
define the elements, for $i = 0, 1,2,\dots, k+1$,
$$
a_i^{(n)} = \cases
\( \frac1{q_n}, \pmb 0 \) & \text{if $i=0$}\\
\( \frac1{q_n}, \epsilon_i \)& \text{if $1 \leq i \leq k+1$} \\
\endcases
$$
Set $S_n = \sum_{i=0}^{k+1} a_i^{(n)} P$, so that $f(n)$ is constant
in $n$ with value  $k+2$. Now we can write (in lots of different ways) $a_i^{(n)}$ as
a nonnegative linear combination of the $a_i^{(n+1)}$, for example,
$a_0^{(n)} = p_{n+1} a_0^{(n+1)}$ and $a_i^{(n)} = \sum_{j \neq 0,2}
a_j^{(n+1)} + 2a_i^{(n+1)} + (p_{n+1}-k-1) a_0^{(n+1)}$ (this exploits the
facts that $\sum_{i=1}^{k+1} \epsilon_i = \pmb 0$ and $p_{n+1} > k+1$; in
fact, we assumed $p_{n+1} > (k+1)^2$ which we will need later). This
yields that $S_n \subseteq S_{n+1}$; the matrices resulting from these
representations are not suitable for our purposes, as the resulting map
$\Psi$ is not one to one.
 
It is elementary that $G^+ = \cup S_n$; an arbitrary element of $G^+$
is of the form $ x = (b/q_n, v)$ where $b$ is a positive integer and $v
\in \Z^k$. Let $d$ be the maximum absolute value of the coordinates of $v$
(as an element of $\Z^k$, i.e., the usual coordinates), and find $l$ so
that $p_{n+1} \cdot p_{n+2}\cdot \dots \cdot p_{n+l} > (k+1)d(b+1)$. We
can find positive integers $r(1), r(2), \dots, r(k+1)$ with
$\sum_{i=1}^{k+1} r(i)\epsilon_i = v $ \st $\sum r(i) < d(k+1)$. Then
$$
\(\frac{b}{q_n}, v\) = \sum_{i=1}^{k+1} r(i) \(\frac 1{q_{n+l}},
\epsilon_i\) + \(b(p_{n+1} \cdot p_{n+2}\cdot \dots \cdot p_{n+l} )-
\sum_{i=1}^k {r(i)} \) \(\frac 1{q_{n+l}}, 0\)
$$
expresses $x$ as an element of $S_{n+l}$.
 
Now we make a very particular choice of the transition matrices, $A_n$;
not only do they have to satisfy (*), but they have to be rank $k+1$ (or
less, but strictly less is not possible, except for finitely many $n$).
Since the matrices are  all square of size  $k+2$, the rank condition
turns out to be not so onerous, especially since imposing the obvious
constraint on the trace will force the rank condition to apply.
 
Temporarily drop the subscript $n$ on some of the variables, so we will
obtain a matrix $A$ whose large eigenvalue is $p$ (corresponding to
$p_n$); we insist that $p > (k+1)^2$. We write,
$$\eqalign{
a_0^{(n)} & = (p -k -1) a_0^{(n+1)} + \sum_{i=1}^{k+1} a_i^{(n+1)}\cr
a_i^{(n)} & = (p-1)a_{0}^{(n+1)} + a_{i}^{(n+1)} \qquad \text{for $1 \leq
i \leq k$}\cr
a_{k+1}^{(n)} & = (p-k^2-k-1)a_0^{(n+1)} + \sum_{i=1}^k k a_{i}^{(n+1)} +
(k+1)a_{k+1}^{(n+1)}.\cr
}$$
These relations are trivially verified by using $\sum_{i=1}^{k+1}
\epsilon_i = \pmb 0$. The corresponding matrix (we act from the left, so
each equation gives rise to a column), is rather simple to describe (but
really tedious to \TeX). There is a $k \times k$ identity matrix occupying
most of the space.
$$
A = \left[  \matrix
p-k-1 & p-1 & p-1 & p-1 & \dots & p-1 & p-1& p- k^2 -k -1 \\
1 &&&&&&& k\\
1 &&&&&&& k\\
\vdots &&& \I_k &&&& \vdots\\
1 &&&&&&& k\\
1 &&&&&&& k\\
1 & 0 & 0 & 0 &\dots  & 0 & 0  & k+1\\
\endmatrix \right].
$$
A consequence of the equations is that the column sums are all $p$. If we
sum all but the first column, the result is $k+1$ times the first column
(as follows from
$p- k^2 - k - 1 + k(p-1) = (k+1)(p-k-1)$). Hence $\rk A < k+2$, that is,
$\rk A \leq k+1$.
 
Now restore the subscript $n$; we obtain square matrices $A_n$ of size
$k+2$ with nonnegative entries, \st $\rk A_n \leq k+1$, each with large
eigenvalue $p_n$ (this latter is not needed now). Then the rank of $H =
\lim \Arrow A_n;\Z^{k+2}.\Z^{k+2}$ is at most $\liminf \rk A_n \leq k+1$,
and so the positive map $\Arrow \Psi;H . G$  obtained from this  sequence
of relations is an isomorphism of ordered groups.
This yields an ECS realization of $G$. \qed

\SecT 3 Arbitrary extensions by free abelian groups 

Now we try to find explicit realizations of dimension groups $G$ that are
given as arbitary extensions of $\Z^k$ by $U$ a subgroup of the rationals.
Explicitly, we have a short exact sequence of abelian groups $0 \to \Z^k
\to G \to U \to 0$; regarding $U$ as a subgroup of the reals, the quotient
map $\Arrow \tau ;G.U$ yields the ordering: for nonzero $g$ in $G$, $g \in
G^+$ iff $\tau(g)> 0$. This describes all dimension groups with unique
trace, whose values lie in the rationals, and whose kernel is free abelian
of finite rank. The previous constructions of course dealt with the case
wherein $\tau$ is split.

\Lem Proposition \fouone. Let $G$ be a noncyclic simple dimension group with unique trace $t$ \st 
$t(G) := U$ is a subgroup of the rationals, and \st $\ker t$ is free of rank $k$. Then $G$ admits an \ecs\ realization by primitive matrices of size $k+2$.

Given the data $(G,\tau, U, \ker \tau)$ we can construct semigroups $S_n
\subseteq S_{n+1} \subset \dots$ of $G^+$ with the property that $G^+ =
\cup S_n$.  To begin with, write $U = \lim \Arrow \times p_n ; \Z.\Z$;
form $q_n = \prod_{i=1}^n p_i$, and make an initial selection, one for each
$n$, of $g_n \in \tau^{-1}(1/q_n)$. Then $g_n - p_{n+1}g_{n+1} \in \ker
\tau$, so we can write $g_n = p_{n+1}g_{n+1}  + v^n$ for unique $v^n \in
\ker \tau = \Z^k$ (obviously, $v^n$ depends on the selection of the
sequence $(g_n)$. The sequence $(v_n, p_{n+1})$ determines the isomorphism
class of $G$, but by itself, this is not that useful.

The obvious candidate for the subsemigroup $S_n$ is the semigroup
generated by 
$$\brcs{g_n, g_n + \epsilon_1, \dots, g_n + {\epsilon_k}; g_n +
\epsilon_{k+1}},
$$ 
where $\epsilon_i$ have their usual meaning: standard
basis elements if $i < k+1$ and $\epsilon_{k+1} = -\sum_{j=1}^k
\epsilon_j$; it is convenient to define $\epsilon_0 = 0$, so we can write
$S_n = \sum_{i=0}^{k+1} (g_n + \epsilon_i)P$. Unfortunately, there is no
guarantee that $S_n \subseteq S_{n+1}$ (in other words, that the matrix entries be nonnegative), largely because at this stage, we
have no control over $v^n$.
 
The idea is to make a better choice of $g_n$, and then a telescoping
(amounting to telescoping the $p_n$), and if we are careful, we will obtain $S_n \subseteq
S_{n+1}$ for all $n$, and the corresponding transition matrices can be
written down; in fact, we will write down the transition matrices, verify
the entries are all nonnegative, from which it follows that the $S_n$ are
increasing.

Let us see what we need to obtain this; we will write down the relations
between the generators, and hope for the best. Fix $n$ and order  $(g_n +
\epsilon_i)$ according to the subscript of the $\epsilon_i$, with $ 0 \leq
i \leq k+1$. The relations are given by $g_n + \epsilon_i = \sum_j A_{ji}
(g_{n+1} + \epsilon_j)$, where $A_{ji}$ are integers, hopefully
nonnegative, and this forces various equations.
 
Since $\tau(g_n+ \epsilon_i)  $  are all equal to $1/q_n =
p_{n+1}/q_{n+1}$ and $\tau (g_{n+1} + \epsilon_i) = 1/q_{n+1}$, we deduce
that for all $j$,  $\sum_{i} A_{ji}  = p_{n+1}$, that is, the row sums of
$A^T$ are all $p_{n+1}$, so that the column sums of $A$ are all $p_{n+1}$.
(This is a useful way of    calibrating the matrix---I am
always confused as to whether it should be $A$ or $A^T$, and keeping in
mind that the column sums must be equal determines which it is.)
 
Now fix $i$; we have the equation, $g_n + \epsilon_i = p_{n+1} g_{n+1} +
\sum_j A_{ji} \epsilon_j$. Using $g_n = p_{n+1}g_{n+1}  + v^n$, we have,
for all $i=0,1,\dots, k+1$ (suppressing the subscript $n$ on $A_n$, as
otherwise it gets too crowded),
$$
v^n + \epsilon_i =  \sum_{j=0}^{k+1} A_{ji} \epsilon_j.
$$
When $i = 0$ (so $\epsilon_i = 0$), we obtain
$$\eqalign{
A_{j,0} - A_{k+1,0} &= (v^n)_j \qquad \text{if $j > 0$} \cr
A_{00} &= p_{n+1} - \sum_{i> 0} A_{i,0} = p_n - (k+1)A_{k+1,0} -
\sum_{i=1}^k (v^n)_i.
}$$
Already we see a problem; the coefficients are suppose to be nonnegative,
and so we require $p_n \geq (k+1)A_{k+1,0} +  \sum_{i=1}^k (v^n)_i$ with
$A_{k+1,0} \geq 0$ (we have no controlÑyetÑon the sum of the coefficients
of $v^n$). Anyway, we continue; for $1 \leq i \leq k$,
$$\eqalign{
A_{i,i}  & = A_{k+1,i} + (v^n)_i + 1 \cr
A_{j,i} & = A_{k+1,i} + (v^n)_j \qquad\text{if $ j \neq i$ and $1 \leq j
\leq k$}\cr
A_{0,i} & = p_{n+1} = \sum_{j=1}^{k+1} A_{j,i} = p_{n+1} - (k+1)A_{k+1,i} 
-\sum_{i=1}^k (v^n)_i - 1.\cr
}$$
Finally, with $i = k+1$,
$$\eqalign{
A_{j,k+1}  &= A_{k+1,k+1} + (v^n)_j -1 \qquad \text{if $k+1 > j > 0$} \cr
A_{0,k+1} &= p_{n+1} - \sum_{i> 0} A_{i,k+1} - k = p_n - (k+1)A_{k+1,k+1}
- \sum_{i=1}^k (v^n)_i -k.\cr
}$$
Now set $a_i = A_{k+1,i}$ (obviously this depends on $n$, but for now we
suppress the sub/superscript); then all the entries are linear in the
choice of $a_i$. If the entries do happen to be nonnegative, then the
resulting matrix $A_n = (A_{ij})$ (order of the subscripts reversed) will implement the
embedding $S_n \subseteq S_{n+1}$. The resulting matrix is similar to the
preceding ones, in that the interior $k \times k$ matrix is $v\cdot \pmb 1^T
 + \pmb 1\cdot (a_1,\dots, a_k) + \I_k$ where $\pmb 1$ is the column of
size $k$ consisting of ones,  we regard $v$ as a column, and $\cdot$
represents the usual product of matrices. Notice that $v\cdot \pmb 1^T$ is
$k \times k$ but $v^T \pmb 1$ is just the sum of the coefficients of
$v$, $\sum_{i=1}^k (v^n)_i$. We sometimes suppress the sub/superscripts $n$ or $n+1$
in $v^n$ and $p_{n+1}$, and the implicit superscripts in $a^{(n)}_i$.

$$
A_n = \(\matrix p_{n+1} - (v^n)^T\pmb1 - (k+1)a_0 & (p_{n+1} - (v^n)^T\pmb1)\pmb1^T - (k+1)(a_1,a_2,\dots,a_k) & *\\
v^n + a_0 \pmb 1 &  v^n\pmb 1^T + \pmb 1 (a_1, \dots, a_k) +\I_k & v^n  + (a_{k+1}-1)\pmb 1\\
a_0 & a_1, a_2, \dots , a_k & a_{k+1} \\
\endmatrix\), \tag1
$$
where the $(0,k+1)$ entry (the upper right; left blank, because of overflow) is $p_{n+1} - v^T\pmb1 - (k+1)a_0 +k$. The column sums are all $p_{n+1}$, as follows from the choice of generators of the subsemigroups. 
\comment
$$
A_n=\left[\matrix
p- (k+1)a_0 - v^T\cdot \pmb 1^T & p- (k+1)a_1 - v^T\cdot \pmb 1^T  - 1 &
 \dots & p- (k+1)a_k - v^T\cdot \pmb
1^T  - 1 & p- (k+1)a_{k+1} - v^T\cdot \pmb 1^T  - k \\
v_1 + a_0 & v_1 + a_1 + 1 & \dots  & v_1 +
a_{k} & v_1 +a_{k +1} -1\\
v_2 + a_0 & v_2 + a_2  +1& \dots & v_2 +
a_{k} & v_2 +a_{k +1} -1\\
\vdots&\vdots&\ddots&\vdots&\vdots& \\
v_k + a_0 & v_k + a_1   & \dots & v_k + a_{k}
+1 & v_k +a_{k +1} -1\\
a_0 & a_1  &\dots &  a_k & a_{k+1} \\
\endmatrix\right]
$$
\endcomment
Without yet worrying about positivity or rank, we can calculate the
eigenvalues and their geometric multiplicities, by   explicitly
computing the left eigenvectors.
 
First, $\pmb 1_{k+2}^T$ is the left eigenvector for $p$. Next, define
${}^{\perp}v:= \Set{w\in \Z^{1\times k}}{wv = 0}$ (we use $\Z^m$ to mean
$\Z^{m\times 1}$, that is, the default is columns). For each $u \in
{}^{\perp}v$, the row of size $k+2$, $(0, u, -u^T \pmb 1)$, is a left
eigenvector for the eigenvalue $1$. If $v \neq \pmb 0$ (as we are assuming
implicitly anyway), then ${}^{\perp}v$ is rank $k-1$, and thus even the
geometric multiplicity of $1$ as an eigenvalue is at least $k-1$.
 
This leaves two eigenvalues. We may find $u_0 \in \Q^{1 \times k}$ \st
$u_0 v = 1 - p + v^T\pmb 1_k$; then $(1, u_0, k+1 - u_0\pmb 1_k)$ is
another left  eigenvector for the eigenvalue $1$ (we may multiply by an
integer and so obtain an integer eigenvector if desired), and since its
first coordinate is not zero, it is not in the $\R$-span of the previous
eigenvectors for $1$; hence the multiplicity of $1$ is at least $k$.
 
There is one remaining eigenvalue, in addition to $p, 1^k$, and it is
easily determined from the trace; the trace of the matrix is $p +k +
\sum_{i=1}^{k+1} a_i - (k+1)a_0$, hence the last remaining eigenvalue is
$\sum_{i=1}^{k+1} a_i - (k+1)a_0$. Since we want the rank of the matrix to
be $k+1$, we are free to choose any selection of integers $a_i$ \st
$\sum_{i=1}^{k+1} a_i = (k+1)a_0$ (that is, $a_0$ is the average of all
the others). When this is imposed, we see quickly that the corresponding
relation holds for the columns, that is, the sum of all but the first
column is $k+1$ times the first. In other words, if we set $z = (k+1, -1,
-1, \dots, -1)^T \in \Z^{k+2}$, then $Az = \pmb 0$. Moreover, $z$ is
independent of the choice of $n$ (that is, $A_n z = \pmb 0$ for all $n$).
 
A particular consequence is that $W:= {}^{\perp}z = \Set {w \in \Z^{1
\times (k+2)}}{wz = 0}$ is a common $A_n$-invariant subgroup (on the left,
of course, meaning $W A_n \subseteq W$ for all $n$); moreover, the eigenvalues of $A_n$ restricted to this subgroup
are exactly $p,1^k$ (the zero eigenvector has conveniently been
eliminated, since $z$ spans, as a real vector space, the right
zero-eigenspace of all the $A_n$).

Now we modify the sequence $(g_n)$ and corresponding $(v^n)$ to permit a selection of integers
$a_{i,n}$ (and with $\sum_{i=1}^{k+1} a_i = (k+1)a_0$) so that the matrix $A_n$ has only nonnegative entries.

\comment
First, we make the easy observation that we can arrange that $v^{2n}$ has
all its coordinates in $\brcs{-(p_{n+1}-1, \dots, -1,0}$ and $v^{2n+1}$
has all its coordinates in $\brcs{0,1,\dots, p_{n+1}- 1}$.
To see this, we proceed by induction, one step at a time. With $n$ fixed,
replace $g_n$ by $g'_n = g_n + w_n$ for some $w_n \in \Z^k$, leaving
the other $g_i$s unchanged. This affects just two of the recurrence
relations,
$$\eqalign{
g'_n &= p_{n+1}g_n + v^n  + w_n\cr
g_{n-1}& = p_n (g'_n ) + v^{n-1} - p_n w_n \cr
}$$
Thus $(v^{n-1})' = v^{n-1} - p_n w_n$, and the new $(v^n)' = v_n + w_n$. Now
we can iterate this with $n$ replaced by $n+1$, replacing $g_{n+1}$ by
$g_{n+1} + w_{n+1}$ so that $v^n + w_n - p_{n+2} w_{n+1}$ has all of its
coordinates in the desired interval. This sequence of one at a time
changes yields a new sequence of ${g_n}'$s (each requiring $n-1$
iterations, and subsequently are unchanged) with the desired properties
for the $v^n{}'$.
 
Relabel the new sequences $(g_n), (v_n)$. Now we telescope in pairs,
starting with the first. This means with $p_n' = p_{2n-1}p_{2n}$, we
define via $g_n' = p'_{n+1} g'_{n+1} + v'{}^n$; to see what the coefficients
are, we note $g_{2n} = p_{2n+1}g_{2n+1} + v^{2n} = p_{2n+1}p_{2n+2}
g_{2n+2} + p_{2n+1} v^{2n+1} + v^{2n}$. With $g'_n = g_{2n}$, we have
$v'{}^n = p_{2n+1} v^{2n+1} + v^{2n}$. Now we recall that $v^{2n+1}$ and
$v^{2n}$ have opposite signs uniformly in each coordinate (where we regard
zero as having opposite sign to itself. It follows that $ \|
v'{}^n\|_{\infty}  \leq p_{2n+2}(p_{2n+1}-1) $. In particular, $p'_{n} -
\| v'{}^n\| \geq p_{2n+2}$.
 
Start with an initial telescoping of the integers giving $U$ so that the
corresponding $p_n \to \infty$. Now apply the two at a time construction
to ensure that the signs of all the coefficients of the coboundaries
alternate in $n$, and perform the telescoping as in the previous
paragraph. Relabel, so that we are in the situation in which $p_{n+1} - \|
v^n\|_{\infty} \to \infty$.
 
%%%%%%%special case $k=1$
Here $k =1$, and so $v^n$ is replaced by the integer $v(n)$ and
$\epsilon_i$ is replaced by $\epsilon$. Suppose $g_n = p_{n+1} g_{n+1}
+v(n)\epsilon$.
 
\Lem Lemma. Suppose that $G$ is an extension of $\Z$ by $U$, and suppose
that there is representation of $U$ by $\brcs{p_n}$ \st  for all $n$,
$v(n) \leq -1$ and $p_{n+1} \geq |v(n)| + 5$. Then $G$ is represented by
the sequence of matrices $A_n$, and by the sequence consisting of their
transposes.
$$
A_n = \left[\matrix
p_{n+1} + v(n) - 2 & p_{n+1}  + v(n)+ 1 & p + v(n)-5\\
1 & 0 & 2 \\
1- v(n) & -1 -v(n) & 3 - v(n) \\
\endmatrix\right]
$$

\Pf The column sums of $A_n$ are all $p_{n+1}$, and the second row less
the third is $v(n), 1+v(n), -1 + v(n)$, so that the relations among the
$g_n, g_n + \epsilon, g_n - \epsilon$ are satisfied. We thus obtain, with
$S_n $ generated by $\brcs{g_n, g_n + \epsilon, g_n - \epsilon}$, that
$S_n \subseteq S_{n+1}$. It is trivial that $G^+ = \cup S_n$, so we obtain
an onto positive map from $H:=\lim \Arrow A_n ;\Z^3. \Z^3$ to $G$. The sum
of the second and third columns is twice the first, so $\rk A_n \leq 2$
(and obviously equality holds). Since $\rk G = 2$, the induced map $\Arrow
\Psi; H.G$ is an order isomorphism. This representation of $G$ has equal
column sums, \ecs.
\endcomment

\comment  
We can deal with $k=1$ and almost all of the case $k=2$.
 
First, when $k=1$; we may assume (after telescoping, if necessary) that
$p_{n+1} \geq 7$ for all $n$. Now each $v^n$ is a single integer, which we
may assume lies in the interval $-p_{n+1}/2 \leq v^n < p_{n+1}/2$. We have
to choose $a_0, a_1, a_2\geq 0$ with $a_2 + a_1 = 2a_0$ \st the resulting
matrix is primitive (in particular, nonnegative).
 
If $v^n\leq 0$, set $a^n_1 =a^n_2 = a^n_0 = 1 - v_n$; the resulting matrix
is then
$$
A_n = \( \matrix  p_{n+1} + v^n - 2 & p_{n+1} + v^n -3 & p_{n+1} +  v^n -1\\
1 & 2 & 0\\
-v^n + 1 & -v^n + 1 & -v^n +1 \\
\endmatrix\).
$$
We note that $p_{n+1} + v^n \geq p_{n+1}/2 \geq \slfrac 72 > 3$, so all
entries are nonnegative (and all but one is positive), and this case is
done.
 
If $ 0 < v^n < p_{n+1}/2$, set $a^n_0 = a^n_1 = a^n_2 = 1$; the resulting
matrix is
$$
A_n =  \( \matrix  p_{n+1} - v^n - 2 & p_{n+1} - v^n  -3 & p_{n+1} -  v^n
-1\\
v_n+ 1 & v_n + 2 & v_n +0\\
1 & 1 & 1 \\
\endmatrix\).
$$
For this to be nonnegative, it is sufficient that $p_{n+1} \geq v_n +3$,
but since $v_n < p_{n+2} /2$ and the latter is at least $\slfrac 72$, we
are done here too.
 
\Lem Corollary. Let $G$ be a simple dimension group with unique trace, 
\st the trace rational-valued and its kernel is $\Z$. Then $G$ can be
realized as a sequence of size three \ecs\ matrices. In addition, if $G
\to \tau(G)$ is nearly split, then $G$ admits an \ers\ representation of
size three.
 
Now we deal with the case of $k = 2$. By the spread result, it suffices to
find a telecoping so that $s(v) < p - k-2$. This is easy to obtain simply
by arranging that signs in each coordinate alternate in $n$, and telescope
in pairs.
 
However, when $k \geq 4$, the vectors approximately given by $(-p/2,0,
0,0)$ have spread $p/2$, and there does not seem to be anything that can
be done in the general case. Explicitly, suppose $p_n$ is a sparse set of
primes (say with $p_n/p_{n+1} < 1/(n+3)$),  $k = 4$, and $v^n = (-(p_n
-1)/2, 0,0,0)^T$. It is plausible that the corresponding simple dimension
group with unique trace not only cannot be realized by a sequence of size
five square nonnegative matrices of the form we have described, but cannot
be realized by a sequence of primitive matrices of any size. [In contrast,
it is known that all dimension groups with unique trace which is
additionally rational-valued---and many more---can be realized as a
sequence of matrices with equal column sum, but without bound on the size
of the rectangular matrices.]
 
The reason for using a sparse sequence of primes is to make sure that
telescoping won't improve the spread (much);  the normalized spread, that
is $s(v_n)/p_{n+1} \sim  \slfrac 12$ hardly changes at all under
telescoping as $n \to\infty$.
 
We excluded the case $k =3$ from this, because it seems to be on the
borderline: $(-p/3,0,p/3)$ has spread $2p/3$ (and we cannot reduce this by
adding multiples of $p$ to the coordinates independently of each other);
we must have $a_1, a_2, a_3 \geq p/3 +1$, which forces $a_0 \geq p/3 + 1$
(since $a_0$ is the average of the others); then $(k+1)a_i > p +1$ while
$v_1 + v_2 + v_3 = 0$; however, in this case, the excess is tiny compared
to $p$.
 
If $v \sim (-p/2, 0, p/4)$ we can replace it by $(p/2,0,p/4)$ and now the
sum of the entries is less than $p$ and the entries are nonnegative, so
this can be realized (if $p > 4k$). Can probably realized all cases $k =3$
(telescope so never actually hit $p/3$ or whatever, ...

More generally, we should define the spread $\mod p$ (where $p$ can be any
positive integer), $s_p(v)$: find the representative of the $\mod p$
equivalence classes of $v_i - v_j$ that is largest among those less than
or equal $2p/3$. Now define $\sigma_p(v) = s_p(v)/p$, and then we consider
the sequence $\brcs{\sigma_{p_{n+1}} (v_n)}$. Sufficient for the entire
sequence to be represented is that $\sigma_{p_{n+1}}(v_n) <1/(k-1)$ for
all but finitely many $n$.
 \endcomment
\def\GL#1{\text{GL($#1$)}}
 Let $G$ be given by the sequence $(p_{n+1}, v^n)$. Let $E \in \GL{k,\Z}$
and $W \in \Z^k$; then the   group extension given by the sequence 
$(p_{n+1}, Ev^n + (p_{n+1}-1)W)$ is equivalent. To see this, form the
square matrices of size $k+1$,  $A_n = \(\smallmatrix  p_{n+1} & \pmb 0\\
v^n &  \I_k\\\endsmallmatrix\)$ and $F = \(\smallmatrix  1 & \pmb 0\\ W
&  E\\ \endsmallmatrix\)$. Then $F^{-1} = \(\smallmatrix  1 & \pmb 0\\
-E^{-1}W &  E^{-1}\\ \endsmallmatrix\)$, and
$$
D_n := F A_n F^{-1} = \(\matrix  p_{n+1} & \pmb 0\\ Ev^n + (p_{n+1}-1)W &
 \I \\ \endmatrix
\).
$$
 
Now $G \iso \lim \Arrow A_n; \Z^{k+1}. \Z^{k+1}$ as abelian groups, and
the map on each copy of $\Z^{k+1}$ given by $F$ induces a group
isomorphism from $G$ to $G' := \lim \Arrow D_{n}; \Z^{k+1}.\Z^{k+1}$. The
corresponding data for the sequence of $D_n$s is $(p_{n+1}, Ev^n +
(p_{n+1}-1)W)$ and these maps preserve the map to $U$. Of course the drawback with this equivalence relation is
that it applies to all the $v^n$s at once.
 
\Lem Lemma \foutwo. Suppose $V \in (\R^+)^k$. There exist $E \in \GL{k,\Z}$
and $W \in \Z^k$ \st  the coefficients of $EV -W$ are all nonnegative, and
sum to less than one.
 
\Pf Define $W^0 \in \Z^k$ via $(W^0)_i = \flo {V_i}$ (the floor function). Then $V_0 := V -
W^0$ has all its coefficients nonnegative and strictly less than one. If
either all or all but one of the entries of $V_0$ is zero, we are done.
 
Otherwise, let $s = \max\brcs{V_i}$ and $t = \max \brcs{V_i \setminus
\brcs{s}}$ (the notation is not very clear, but if there is a tie for
maximum, then $t = s$). We apply the division algorithm to $s$ and $t$;
there exists an integer $m > 1$ \st $ s = mt + s'$ where $0 \leq s' < t$;
this is implemented by an elementary transformation, hence by an element
of $\GL{k,\Z}$, and the new vector (replacing $s$ in one its positions by
the smaller $s'$) either has strictly smaller maximal entry, or the
multiplicity of its maximal entry has been reduced. If in the resulting
vector, there are still more than one entry, we can continue the process.
 
 It is easy to see that it either terminates in a single nonzero entry
(which occurs precisely when all the nonzero entries of $V$ are rational
multiples of each other), or we can make the maximal entry as small as
we like, say less than $1/k$. Either way, we have constructed $E$ as a
product of elementary transformations (hence in $\GL{k,\Z}$) \st
$E(V-W^0)$ has only nonnegative entries and whose entries sum to less
than $1$. Now set $W = EW^0 \in \Z^k$.
\qed
 
\Lem Lemma \fouthr. Suppose that the extension $G$ of $\Z^k$ by $U$ is implementable by
$(p_{n+1}, v^n)$ \st  $v^n/p_{n+1}$ converges (in $\R^k$) and $p_{n+1} \to
\infty$. Then the corresponding dimension group $G$ with unique trace
being the map to $U$ is realizable as a limit $\Arrow A_n ; \Z^{k+2}.
\Z^{k+2}$ with $A_n$ of the form (1)  primitive.
 
\Pf Set $V = \lim v^n/p_{n+1}$. By the preceding, there exists $E \in
\GL{k,\Z}$ and $W \in \Z^k$ \st $EV - W$ has only nonnegative entries
adding  to $\lambda < 1$. Now the extension corresponding to $(p_{n+1},
v^n)$ is equivalent to $(p_{n+1}, (v^n)':=Ev^n - (p_{n+1}-1)W)$, so it
suffices to show that $(p_{n+1}, (v^n)')$ can be realized by a sequence of
primitive matrices of the form (1).
 
We observe that $(v^n)'/p_{n+1} = Ev^n/p_{n+1} - W(p_{n+1}-1)/p_{n+1}$,
and this sequence converges to $EV - W$. Thus, given $\epsilon <
\min\brcs{(1-\lambda)/3(k+1), \lambda/3k}$,  for all sufficiently large
$n$, we have $- \epsilon \pmb 1 \leq (v^n)'/p_{n+1} $ and the sum of the
entries  is less than $\lambda + \epsilon$.  Thus 
$$
-\epsilon p_{n+1} \leq (v^n)'_i    \qquad \text{and} \qquad  \sum 
(v^n)'_i  < \lambda + \epsilon p_{n+1}.
$$
Set $a_0 = a_1 = \dots = a_{k+1} $ to be $1$ if $\min (v^n)'_i \geq 0$ and
equal to $1 -  \min (v^n)'_i$ otherwise.
 
If $\min (v^n)'_i \geq 0$, then we note that from $\sum(v^n)'_i  \leq
(\lambda + \epsilon)p_{n+1}$, we obtain an upper bound on the sum, 
$p_{n+1}\mu$, where $\mu =  1 - (1-\lambda)(1-1/3k)$ (what is important is
that the coefficient is bounded above away from one, uniformly in
sufficiently large $n$). Then $p_{n+1} - \sum(v^n)'_i - (k+1)-1 \geq
p_{n+1}(1- \mu) - (k+2)$. Since $p_{n+1} \to \infty$, for all further
sufficiently large $n$, this expression is positive; thus the matrix $A_n$
in  (1) has only nonnegative entries.
 
If $\min (v^n)'_i < 0$, then for any $j$,
$$\eqalign{
p_{n+1} - \sum(v^n)'_i - (k+1)a_j-1 & \geq p_{n+1} - p_{n+1} (\lambda +
\epsilon) - (k+1) - (k+1)p_{n+1} \epsilon - 1\cr
& \geq p_{n+1} (1 - \lambda - (k+2)\epsilon) - (k+1).\cr
}$$
Now $1 - \lambda - (k+1)\epsilon > 1 - \lambda - (1 - \lambda)/3  =
2(1-\lambda/3) > 0$. Hence by further increasing $n$, we have that for all
sufficiently large $n$, the matrix entries of $A_n$ are nonnegative. We
can always delete a finite number of the matrices at the outset. Because $(k+1)a_0 
= \sum_{i=1}^{k+1} a_i$, the rank of each $A_n$ is $k+1$, so that Lemma \oneone\ applies.
\qed
 
\Lem Lemma \foufou. Let $G$ be a group extension of $\Z^k$ by noncyclic $U \subseteq \Q$
with data $(p_{n+1}, v^n)$. Then there is an equivalent representation,
$(q_{n+1}, (v^n)')$, \st $q_{n+1}$ is increasing,  $q_{n+1} \to \infty$,
and $(v^{n})'/q_{n+1}$ converges.
 
\Pf First, we may make an initial telescoping, and thus may assume that
$p_{n+1}$ are increasing to infinity at the outset. Now we perform the following
substitution transform, to ensure that the resulting $v_n$ entries are all
between $0$ and $p_{n+1}-1$. Suppose we have done this up to $n =m-1$;
that is, $g_i = p_{i+1}g_{i+1} + v^i$ for $1 \leq i \leq m-1$. Now set
$g_{m+1}' = g_{m+1}+ u^m$, the $u^m$ to be determined. Then we have $g_m =
p_{m+1}g_{m+1}' + v^m - p_{m+1} u^m$ (and the subsequent relations, for
larger $m$, are also affected, but we come to them by the induction
argument). We can obviously choose $u^m \in \Z^k$ so that all the entries
of $v^m - p_{m+1}u^m$ lie in the set $\brcs{0,1,2,\dots, p_{m+1}-1}$. This
completes the induction, and allows us to assume that the newly relabelled
$v^{n}$ satisfy $\pmb 0 \leq v^n \leq (p_{n+1} - 1)\pmb 1$.
 
In particular, with the current notation, $\brcs{v^n/p_{n+1}}$ is a
bounded sequence in $[0,1]^k$. Hence there exists a subsequence indexed by
$n(i) \in \N$
\st $\brcs{v^{n(i)}/p_{n(i)+1}}$ converges, say to $V \in [0,1]^k$. The
integers $n(1) < n(2) < n(3) < \dots \to \infty$ suggest a telescoping; set $M_j =
\( \smallmatrix p_{j+1} & \pmb 0 \\ v^j & \I_k \\ \endsmallmatrix \)$,
discard the $M_j$ for $j < n(1)$, and define
$$
M^{n(i)} =        M_{n(i+1)-1}  \cdot \dots     \cdot M_{n(i)+1}\cdot
M_{n(i)};
$$
the upper left entry is $q_i = \prod_{j=1}^{n(i+1)- n(i)} p_{n(i)+j}$, and
of course the lower right $k \times k$ is the identity matrix. The column
of size $k$ to the left of the identity is obtained by an easy induction
argument. This yields (and probably
better is to prove this by induction):
$$\eqalign{
\frac{v^{(i)}}{q} &= \frac {v^{n(i)}}{p_{n(i)+1}} +\frac
{v^{n(i)+1}}{p_{n(i)+1}p_{n(i)+2}}  +\frac
{v^{n(i)+2}}{p_{n(i)+1}p_{n(i)+2}p_{n(i)+3}} + \dots  \cr
\left \| \frac{v^{(i)}}{q} - \frac {v^{n(i)}}{p_{n(i)+1}}\right
\|_{\infty} & \leq \sup \left\| \frac{v^j} {p_{j+1}} \right\| \( \frac
1{p_{n(i)+2}}   + \frac 1{p_{n(i)+2} p_{n(i)+3}}  + \dots  \) \cr 
& \leq \delta \( 
\frac 1{p_{n(i)+2}} \( 1 + \frac 1{p_{n(i)+2}} + \( \frac
1{p_{n(i)+2}} \)^2 +\( \frac 1{p_{n(i)+2}} \)^3  + \dots \)\)\cr
& =\frac {\delta}{  p_{n(i)+2}-1 }. \cr
}$$
The $\delta$ was obtained from $\brcs{v^n/p_{n+1}}$ being a bounded
sequence. We have used $p_{n+2} \geq p_{n+1}$ to convert the estimate into
a geometric series. Finally, since $p_{n+1} \to \infty$, the sequence 
$\brcs{\frac{v^{(i)}}{q_i}} $ converges.
\qed

The previous two  results yield Theorem \fouone.

\SecT 4 Arbitrary extensions

In this section, we deal with arbitrary   extensions of $U$ by arbitrary finite
rank torsion-free abelian groups, instead of $\Z^k$; that is, let $C$ be a
rank $k$ torsion-free abelian group, $U$ a noncyclic subgroup of the
rationals, and consider extensions $0 \to C \to G \to U$. We realize
the dimension group $G$ with strict ordering induced by $G \to U \subset
\R$ as a limit of \ecs\ matrices of  size  $k+2$.
The arguments are rather similar to those of the previous section, but involve a couple of extra features. 

Let $C = \ker t$; this has rank $k$, so we can write $C$ as a limit (as
abelian groups) $\Arrow B_n; \Z^k.\Z^k$ for some choice of $B_n$ with
$\det B_n \neq 0$. We can incorporate the identity as many times as we
wish, say $B_1, \I_k , \dots, \I_k, B_2, \I_k, \dots, \I_k, B_3, \dots$,
and this gives the same abelian group (the idea is that we will be telescoping the $k+1 $ size matrices, and we want to ensure that the absolute column sums of the $B_n$ are $\oh{p_{n+1}}$). Re-indexing, we can obtain $G$ (the
extension) as the abelian group direct limit arising from the square 
matrices of size $k+1$,  $\(\smallmatrix p_{n+1} & \pmb 0 \\  v^n & B_n \\
\endsmallmatrix \)$, and we can assume that $\| B_n \|_{\infty,\infty} = \oh{\sqrt{p_{n+1}}}$ (or anything reasonable). 
 
When we realize the corresponding semigroup coming from the relations, we
obtain rather similar matrices to those previously encountered. Let
$\epsilon^n_i$ be the standard basis elements for $\Z^k$ at the $n$th
level, and define $\epsilon^n_0  = \pmb 0$
and $\epsilon^n_{k+1} = - \sum_{i=1}^k \epsilon^n_i$. 

We can express $G$ (as an abelian group with real-valued homomorphism) as the limit, $G = \lim \Arrow M_n:= \(\smallmatrix p_{n+1} & \pmb 0 \\  v^n & B_n \\
\endsmallmatrix \); \Z^{k+1}.\Z^{k+1}$, where the common left eigenvector $x = (1,0,0, \dots,0)$ induces the map to $U$. Call the $n$th copy of $\Z^{k+1}$, $F_n$, so the elements of $G$ are the equivalence classes $[a,n] = [M_n a, n+1]$  where $a \in F_n$, and the map to $U$ is given by $t:[a,n] \mapsto xa/p_1\dots p_{n}$. In particular, the corresponding positive cone (of the dimension group, once we put the strict ordering arising from $t$) will be $\cup t^{-1}(1/p_1\dots p_n)P$.

Set $g_n = [x^T := (1,0, \dots,0)^T, n]$; so  $t(g_n) = 1/p_{1}\cdots p_n$. Then $g_n - p_{n+1}g_{n+1} = [M_n x^T - x^T, n+1]$, and this is simply the column whose top entry is $0$ and the rest of which is $v^n$, which we write as $\sum_{i=1}^{k} v^n_i \epsilon^{n+1}_i$. 

Take as our generators for the positive cone at the $n$th level, $g_n + \epsilon^n_i$ (now $i= 0, 1,\dots, k+1$ to incorporate the two extra elements required); we rewrite this as $g_n + \sum_{j=1}^k (B_n)_{ji}\epsilon^{n+1}_j$ (arising from the effect of $B_n$; note that the coefficients are transposed). From the equation $g_n = p_{n+1}g_{n+1} + \sum_{i=1}^{k} v^n_i \epsilon^{n+1}_i $, we want to find a $(k+2) \times (k+2)$  matrix $A^n:= (A^n_{ij})$ ($0 \leq i,j\leq k+1$; this is not supposed to represent the $n$th power of some matrix $A$, but is merely superscripting the indexÑthe previous notation was $A_n$; we frequently drop the super/subscript $n$ for simplicity) \st  for all $i$,
$$
g_n + \sum_{j=1}^k (B_n)_{ji}\epsilon^{n+1}_j = \sum_{j=0}^{k+1} (A^n)_{ji} (g_{n+1} + \epsilon_j^{n+1}).
$$
We are free to choose the entries of $A^n$ (subject of course to positivity constraints) so long as this set of equations holds. First, from $t(g_n + \epsilon^n_j) = 1/p_{1}\dots p_n$, the column sums of $A^n$ must all be $p_{n+1}$.  Hence we require $\sum_{j=0}^{k+1} A^n_{ji} = p_{n+1} $ for all $i$. Using the relation $g_n - p_{n+1}g_{n+1} = \sum_{i=1}^{k} v^n_i \epsilon^{n+1}_i $, the previous equations become
$$
\sum_{j=0}^{k+1} A^n_{ji}  \epsilon_j^{n+1}= \sum_{l=1}^k v^n_l \epsilon_l^{n+1} + \cases 0 & \text{if $i=0$}\\
\sum_{1}^k (B_n)_{ji}  \epsilon_j^{n+1} & \text{if $1 \leq i \leq k$}\\
-\sum_{l=1}^k \sum_{j=1}^k (B_n)_{jl} \epsilon_j^{n+1} & \text{if $i= k+1$.} 
 \endcases 
$$
Restrict to the case $j>0$ (each entry in the top row is determined by the remaining ones in its column, as the column sums are all $p_{n+1}$); taking coordinates, we obtain the following equations, successively obtained by setting $i =0$, $1 \leq i \leq k$, and $i= k+1$,
$$\eqalign{
A_{j0} &= v^n_j + A_{k+1,0} \cr
A_{ji} & = v^n_j + (B_n)_{ji} + A_{k+1,j} \qquad \text{for $1 \leq i \leq k$} \cr
A_{j,k+1} & = v^n_j - \sum_{l=1}^k (B_n)_{jl} + A_{k+1,k+1}.\cr
 }$$ 

Set $a_i = A_{k+1,0}$ (of course, we should really write this as $a_{(i,n)} = A^n_{k+1,0}$ to indicate dependence on $n$), so we obtain these as free parameters, which determine all the rest of the entries (and $A_{0,i} = p_{n+1} - \sum_{j=1}^{k+1} A_{ji}$). The matrix $A^n$ has the following form (recalling that $\pmb 1$ is the column of $1$s of size $k$). The following display is labelled (**), but because of the width, the label is not obvious. 
$$
A^n =\(\matrix p_{n+1} - v^T\pmb1 - (k+1)a_0 & (p_{n+1} - v^T\pmb1)\pmb1^T - (k+1)(a_1,a_2,\dots,a_k) & \\
v^n + a_0 \pmb 1 & B_n + v^n\pmb 1^T + \pmb 1 (a_1, \dots, a_k) & v^n - B_n \pmb 1 + a_{k+1}\pmb 1\\
a_0 & a_1, a_2, \dots , a_k & a_{k+1} \\
\endmatrix\),\tag{**}
$$
where the $(0,k+1)$ entry (the upper right; left blank, because of horizontal overflow) is $p_{n+1} - v^T\pmb1 - (k+1)a_0 +\pmb1^T B_n \pmb1$. The matrix $B_n + v^n\pmb 1^T + \pmb 1 (a_1, \dots, a_k) $ appearing in the middle is a $k \times k$ block. If we sum all of the columns except the leftmost, we obtain $(k+1)v^n + \(\sum_{i=1}^{k+1} a_i\)\pmb 1$; thus if we impose the condition $ \sum_{i=1}^{k+1} a_i = (k+1)a_0 $, the rank of the matrix $A^n$ is at most $k+1$. 

The only restriction to deal  with is positivity. We find a telescoping of $M_n$, so that in the resulting telescoping and transformation, the column sums of the corresponding $B_n$s plus the sum of the entries of $v$ plus $(k+1) a_i$ is less than the (new) $p_{n+1}$ obtained from the telescoping.

We begin, as in the preceding case, with the original relations, $g_n = p_{n+1} g_{n+1} - v^n$. Replace $g_{n+1}$ by $g_{n+1}' = g_{n+1} + u^n$ (where $u_n$ is to be determined), so that the new relation is $g_n = p_{n+1}g_{n+1}' + v^n - p_{n+1}u^n$. So we may choose $u^n$ so that $v'_n = v^n - p_{n+1}u^n$ has all its entries in $\brcs{0,1,2, \dots, p_{n+1}-1}$. Now the relation for $g'_{n+1}$ in terms of $g_{n+2}$ can  be  adjusted, and we continue by induction. Relabelling everything in sight (including the matrices $M_n$), we are now in the situation that $v^n \geq 0$ and $\| v^n \| < p_{n+1}$. 

Since $\brcs{v^n/p_{n+1}}$ is a bounded set of $\R^k$, it contains a convergent subsequence, say $v^{n(i)}/p_{n(i) +1} \to V \in [0,1]^k$. 
This yields an obvious telescoping; set $M^{(i)} = M_{n(i+1)} \cdot M_{n(i+1)-1}\dots M_{n(i)+1}$ and $q_{i+1} = \prod_{j= n(i)+1}^{n(i+1)} p_j$, and $B^{(i)} = B_{n(i+1)} \cdot B_{n(i+1)-1}\dots B_{n(i)+1}$; then $M^{(i)} = \(\smallmatrix q_{i+1} & \pmb 0 \\ v^{(i)} & B^{i}\endsmallmatrix \)$. 

The column $v^{(i)}$ has a relatively simple expression, 
$$
v^{(i)} = q_{i+1} \(\frac {v^{n(i+1)}}{p_{n(i+1)+1}}  + 
\frac {B_{n(i+1)-1}v^{n(i+1)}-1}{p_{n(i+1)+1}p_{n(i+1)}}+ 
\frac {B_{n(i+1)-1}B_{n(i+1)-2}v^{n(i+1)-1}-1}{p_{n(i+1)+1}p_{n(i+1)}p_{n(i+1)-1}} + \dots
\).
$$
Hence 
$$
\left\| \frac{v^{(i)}}{q_{i+1}} - \frac {v^{n(i+1)}}{p_{n(i+1) + 1}}\right\| \leq \sum \left\| \frac{v^{n(i+1) - j}}{p_{n(i+1)-j+1}}\right\| \cdot 
\left\| \frac{B_{n(i+1) - j+1}}{p_{n(i+1)-j+2}}  \right\| \cdot \dots \cdot 
\left\| \frac{B_{n(i+1) - 1}}{p_{n(i+1)}} \right\|
$$
(The norm on the matrices is the maximum absolute column sum, which is either $\infty-\infty$ or $1-1$.) Since we have made the norms of $B_n$ be $\oh{p_{n+1}}$, this goes to zero. Hence $v^{(i)}/q_{i+1} \to V$.
As before, we find $W \in \Z^{k}$ and $E \in \GL{k,\Z}$ \st the absolute sum of the entries of $ EV - W$ is less than one. Now conjugate the matrices simultaneously with $D = \(\smallmatrix 1 & \pmb 0 \\ W & E\endsmallmatrix \)$ as before, and the new matrices are of the form $\(\smallmatrix q_{i+1} & \pmb 0 \\ v^{(i)} & B^{i}\endsmallmatrix \)$ where we can make the substitution for $a^{(n)}_i$ as we did in the previous case with a slight modification; as before, we make $a_{(i,n)} = a^{(n)}$ equal to each other for $1 \leq i \leq k$, and we require that that 
$$
\frac{\min\brcs{p_{n+1} - (v^n)^T \pmb 1, p_{n+1} - (v^n)^T \pmb 1 + \pmb 1^T B_n \pmb1}}{k+1}  \geq a^{(n)} \geq \max\brcs{0,\pmb 1^TB_n\pmb 1 - v^n_i, -v^n_i, -(B_n)_{ij} - v^n \pmb 1^T}_{1 \leq i,j \leq k}
$$ in order that the resulting matrices be nonnegative, and of rank $k+1$. But $\|B_n\|_{\infty,\infty}  = \oh{\sqrt{p_{n+1}}}$ (which obviously persists after the telescoping), so the entries of $B_n/p_{n+1}$ go to zero; dividing the expressions by $p_{n+1}$, the $B_n$ entries contribute negligibly to the obstruction.  Now Lemma \oneone\  applies, and we have a realization of $G$ by \ecs\ matrices of size $k+2$. This yields:

\Lem Theorem \sixone. Let $\Arrow t;G. U \subseteq \Q$ be a simple dimension group with unique trace $t$, \st $t$ is rational-valued. If $\rk G = k+1$, then $G$ admits an \ecs\ realization by matrices of size $k+2$, and they are of the form $(**)$.

\SecT 5 Nearly ultrasimplicial dimension groups

Effros called a dimension group {\it ultrasimplicial\/} if it has a realization as ordered direct limit (with nonnegative matrices) $G \iso \lim \Arrow A_n; \Z^{f(n)}. \Z^{f(n+1)}$ where $\ker A_n = \brcs{0}$ for all $n$
 (so the obvious map  $\Z^{f(n)} \to G$ is one to one). Elliott [E] showed that the simple dimension group with unique trace $\Z[\slfrac 12] \oplus \Z$ (the trace is the projection onto $\Z[\slfrac12]$; this is the split case, covered by Proposition \twoone\ with $p_n = 2$ for all $n$ and $k=1$) is not ultrasimplicial, whereas any totally ordered group is ultrasimplicial, and Riedel [R1] showed that if $G$ is free of finite rank and with unique trace, then it is ultrasimplicial. It follows easily from Riedel's result that if $G$ is a simple dimension group with unique trace $\tau$ and $\rank \tau (G) > 1$, then $G$ is ultrasimplicial ([H$1\slfrac12$]). This   is practically the complementary class to the dimension groups considered here (which are characterized by $\rk \tau(G) =1$ and $\rk G < \infty$) among simple dimension groups with unique trace.

Motivated by the results here, we say a dimension group is {\it co-rank one ultrasimplicial\/} if there exists a realization as partially ordered groups $G \iso \lim \Arrow A_n ;\Z^{f(n)}. \Z^{f(n+1)} $ (as usual, $A_n$ have only nonnegative entries, and the free groups are equipped with the coordinatewise ordering) \st the kernel of any telescoping (with $m > n$) $A_{m} A_{m-1} \dots A_{n+1} A_n$ has rank at most one (alternatively, the map $\Z^{f(n)} \to G$ given by $x\mapsto [x,n]$ has kernel of rank at most one). Then among other things, combining the ultramatricial results ([H$1\slfrac12$], Corollary 4) with  Theorem \sixone, we obtain that any finite rank simple dimension group with unique trace is co-rank one ultrasimplicial. A simple direct limit argument extends this to 

\Lem Theorem \sevone. Every simple dimension group with unique trace is co-rank one ultrasimplicial. 

Riedel [R2] also showed that some free  rank three simple dimension groups with two pure traces are not ultrasimplicial. It is  possible that every simple dimension group is co-rank one ultrasimplicial, although I really doubt it.  

\SecT 6 Good and not-so-good traces

The previous results showed that if $G$ is a simple dimension group with  unique trace, the trace is rational-valued,  and $G$ is of rank $k+1$, then it admits an \ecs\ representation of size $k+1$. 
For this section, we drop the requirement that $G$ have unique trace. We show that if $(G,\tau)$ is a dimension group (having order units) and $\tau$ is a rational-valued trace, then $G$ admits an \ecs\ realization with $\tau$ obtained from the sequence of rows consisting of multiples of $\One^T$ if and only if $\tau$ is good (as defined in [BeH] and below; when the trace is unique, it is automatically good). However, even in the finite rank case, the argument does not yield {\it bounded\/} \ecs\ realizations. 

Suppose $G = \lim \Arrow M_i ; \Z^{n(i)}. \Z^{n(i+1)}$ is an \ecs\ representation of the dimension group $G$, with the $i$th  matrix having row sum  $c_i$.  Then \ecs{} merely says that $\One_{n(i+1)}^T M_i = c_i\One_{n(i)}^T$. This allows us to define a trace $\tau$ on $G$, via $\tau([w,j] ) = \One_{n(j)}^Tw/\prod_{i=1}^{j-1} c_i$.  We call this trace the {\it trace associated to the representation of $G$ via $M_i$.} Obviously $\tau$ is {\it faithful\/} (that is, $\ker \tau \cap G^+ = \brcs{0})$.

Different \ecs{} realizations of the same group $G$ can yield inequivalent traces, moreover, some of the traces so obtained can be pure, while others need not, and their value groups may differ. For example, consider the situation with
 $M_j =\(
\smallmatrix 1 & 2^j \\ 2^j & 1 \\ \endsmallmatrix\)$, a well-known
construction with two pure traces; the trace obtained from this \ecs{} representation is not pure.

For a particular \ecs{} realization, the value group of the trace is $\tau (G) = \cup \frac {1}{\prod_{i=1}^j c_i}\Z$. Thus $\tau(G) \subseteq \Q$.  In particular, if each $c_j$ is a power of the same integer $k$, then $\tau(G) = \Z[1/k]$. 

The set of order units of $G$ will be denoted $G^{++}$. Following [BeH], a trace $\Arrow \tau;G.\R$ is {\it good\/} if for all $b \in G^+$, $\tau([0,b]) = \tau(G) \cap [0,\tau(b)]$; it is {\it order unit good\/} when this property holds for all $b$ in  $G^{++}$. Notation that is not explained here will probably be found in [BeH].

\Lem Theorem \sixtwo. Suppose that $G$ is a dimension group with order unit and let $\tau$ be a trace of $G$. \hfill
\item{(a)} Suppose there is  \ecs{} realization of $G$ implementing  $\tau$. Then $\tau$ is a faithful good trace with rational values.
\item{(b)}
 Let $(G,u,t)$ be a countable dimension group with order
unit and trace such that $t(G)$ is a subgroup of the rationals. If  $t$  is
faithful and good (as a trace), then there
exists a realization of $G$ as a direct limit of simplicial groups whose
realizing matrices have the equal column sum property \st $t$ is the corresponding trace.

It is possible that the hypothesis in (b) that $t$ be good can be weakened to refinability ([BeH]) of $t$. 

Towards (a), we have already observed that $\tau$ is a faithful rational-valued trace. To show that $\tau$ is good, we have an elementary lemma. 

\Lem Lemma \sixthr. Let $n$ be a positive integer, and $\Z^n$ the simplicial group of
rank $n$. Define a trace $t$ on $\Z^n$ by $t(v) = \One_n^T v$ (so the vector $v$ is
sent to the sum of its coordinates). Then $t$ is good.

\Pf Select nonnegative vectors $a= (a_i)^T$, $ b =(b_i)^T$ in $\Z^n$ \st
$\sum a_i < \sum b_i$. We fix $b$ and systematically alter $a$. Let
$S_-(a) = \Set{i}{a_i > b_i}$,  $S_+ (a) = \Set{j}{a_j < b_j}$, and
$S_0(a) = \Set{i}{a_i = b_i}$. Obviously $S_+(a)$ is not empty.  If
$S_-(a)$ is empty, we are finished; otherwise, we proceed by induction on
$\sum_{S_-(a)} (a_i - b_i)$. Select $j \in S_+(a)$, and $k \in S_- (a)$,
and define $a'$ by subtracting $1$ from $a_j$ and adding $1$ to $a_k$, and
leaving the rest of the entries unchanged. Then $\One_n^T a' = \One_n^T a$,
$S_+ (a) \subset S_+(a') \cup S_0 (a')$, and $S_-(a') \subseteq S_-(a)$,
and moreover, $\sum_{S_-(a')} (a_i' - b_i) \leq \sum_{S_-(a)} (a_i - b_i)
- 1$. The transformation $a \mapsto a'$ is repeated until the $S_-$-set is
empty, and we are done.\qed

Now $\tau$ is the inverse limit of traces obtained as in this lemma, so is the limit of good traces, hence is good. This concludes the proof of (a). \vskip 5pt

\noindent {\it Proof\/}  of (b). Start with an arbitrary $\Z^n$ with basis $\brcs{e_j}$ and map $e_j \mapsto
g_j \neq 0$ in $G^+$.  Let $t(G) = \cup_k \prod_{i \leq k} m(i)^{-1} \Z$
(i.e., the $m(i)$s are the successive factors realizing the supernatural
number of $t(G)$); set $\Cal M_k =  \prod_{i \leq k} m(i)$. There exists $k$
such that each $t(g_j) = a_j/\Cal M_k$ for some positive integer $a_j$.

By goodness, there exists $h \in G^+$ \st $t(h) = 1/\Cal M_k$, so again by
goodness, there exists $h_{jl}$ in $G^+$ \st $g_j = \sum_{l=1}^{a_j}
h_{jl}$. This allows us to create a simplicial map $\Z^n \to \Z^{\sum
a_j}$
by sending  $e_j \mapsto \sum_{l=1}^{a_j} E_{jl}$, and we also have the
obvious map from $\Z^{\sum a_j}$ to $G$ via $E_{jl} \mapsto h_{jl}$; then
the maps to $G$ are compatible.

The upshot of this preliminary construction is that all the basis elements
of the new simplicial group are sent to the same value under $t$.
Now we apply the usual construction (as in [EHS]), that is,  adjoin the
next pre-selected generator of the positive cone, make it the image of a
map, and fix up the kernel, so we arrive at the following (all maps are positive):
$$\diagram
\Z^{n(1)}&\rTo^{} & \Z^{n'} \\
\dTo^{} &&\dTo^{} \\
G&\rTo^{=} & G\\
\enddiagram%     the vertical and diagonal strokes indicate +ve maps
 % G
%endof diagram
$$
\st the kernel of the left vertical map is contained in the kernel of
horizontal map, and $\Z^{n(1)}$ is the $\Z^{\sum a_j}$ of two paragraphs
above. We extend the horizontal map to a better simplicial group.

The standard basis elements $E_i$ of the left simplicial group map to
elements $h_i$ with the property that $t(h_i) = 1/\Cal M_k$ for some $k$. Let
$F_j$ be the standard basis elements of the right, say with images $g_j$
(we have re-indexed the bases). Then $h_i= \sum b(i,j) g_j$ for some
integers $b(i,j)$, so that $1/\Cal M_k = \sum b(i,j)t(g_j)$. Then applying the
method of the preliminary construction, we obtain a map $\Z^{n'} \to
\Z^{n''}$ (together with a map to $G$) such that the images of the new
basis elements all have value at $t$ equalling $1/\Cal M_{k'}$ for some $k'
\geq k$ (we can make sure that $k' > k$ for infinitely many iterations of
this process).

So we are in the following situation:
$$\diagram
\Z^{n(1)}&\rTo^{} & \Z^{n''} \\
\dTo^{} &&\dTo^{} \\
G&\rTo^{=} & G\\
\enddiagram%diagram
$$
with the generators of the left group mapping to $1/\Cal M_k$  under  $ t$
 and the generators of the right group mapping to $1/\Cal M_{k'}$ under $t$;
and of course, the  kernel of the vertical map from the preceding $\Z^n$
is contained in the kernel of $ \Z^n \to \Z^{n(i)}$, etc.
The image of $E_i$ in the right group is the $i$th column of the
transition matrix; if the image of $E_i$  is $\sum c_j F_j$,  applying
$t$,  we obtain $1/\Cal M_k = \sum c_j/\Cal M_{k'}$. Hence $\sum c_j$, the column
sum, is independent of  the choice of column. So the transition matrix has
equal column sums.

Now we repeat this process with the new $ \Z^{n''} $ (adjoin the next
element of the positive cone etc). Since this sequence of transition
matrices just obtained intertwines the sequence built up via the [EHS]
method, both give the same dimension group as limit. \qed

\comment
Let's see, where did we use goodness: if  t  is a good trace, then for any
a,b in G^+ with t(a) < t(b), there exists a' in G^+ such that  0 \leq a'
\leq b and  t(a') = t(a) (actually, it says more: a need not be in G^+,
merely that t(a) > 0). Moreover, goodness implies refinement, which says
that whenever  t(b) = \sum t(a_i}
 with b and all the a_i s  in G^+, then there exist  a_i' in G^+ with
t(a_i) = t(a_i') and b = \sum a_i'.

*with a little work, this can be weakened to refinable, I think.
\endcomment

An example given in [BeH] is a simple rank two dimension group with two
pure traces, such that the value groups are both $\Z[\slfrac12]$ and their
 kernels are discrete. In fact, in that example, there are no additive
functions (let alone traces, pure or impure) $\Arrow t; G.\Q$ \st the
kernel is not cyclic; in particular, none of the countably many traces with rational value groups is good (by [BeH, 1.8], the kernel of a good trace $\tau$ has dense range in $\tau^{\vdash}$).

To prove this, we note that  $G$ is strongly indecomposable and an
extension of a cyclic group by $\Z[\slfrac12]$;  now  if $\ker t$ were not
cyclic, there would be (up to isomorphism) a noncyclic subgroup of $\Q$
sitting inside $G$. Applying one of  traces to this subgroup, we see that
it must be disjoint from the kernel, so that its image in $\Z[\slfrac12]$
is an isomorphic copy. But this forces the supernatural number of the
subgroup to be $2^{\infty}$, hence the subgroup is $2$-divisible, hence the
restriction of the trace is of finite index, and therefore we have  a
splitting from a finite index subgroup of $G$, which is impossible, as $G$
is strongly indecomposable.

Hence the kernel of any trace of $G$ with rational values is either zero
or cyclic. Since $G$ is simple, this means that no trace with rational
values can be good, and thus $G$ cannot be represented by an \ecs{} limit.

\comment
For a simple dimension group $G$, a trace $\tau$ (not necessarily rational-valued) is good iff the image of $\ker \tau$ under the affine representation $G \to \Aff S(S,u)$ is norm-dense in $\tau^{\vdash} := \Set{h \in \Aff S(G,u)}{h(\tau) = 0}$ (for general dimension groups with no discrete traces, this is equivalent to order unit good). 

 \Lem Lemma. Let $G$ be a dimension group and $U \subseteq S(G,u)$. Suppose there is a subgroup $H$ of $G$ with the following properties:
\item{¥} $H \cap \ker U = 0$
\item{¥} the norms $\sup_{\tau \in S(G,u)} |\tau(h)|$ and $\sup_{\tau \in
Z} |\tau(h)|$ are equivalent on $H$ (yield the same topology)
%\item{¥} the closure of $H \oplus \ker U$ is a .
{\par}\noindent Then $\ker U$ has dense image in $Z^{\vdash} \cap \overline{H \oplus \ker U}$.

\Pf Given $j$ in $Z^{\vdash}$, there exist $g_n \in H \oplus \ker U$ \st
$\| \hat g_n - j\| \to 0$. We may write $g_n = h_n + k_n$ where $h_n \in
H$ and $k_n \in \ker U$. For every $\tau$ in $Z$, we have $|\tau(h_n) | =
|\tau(g_n) - j(\tau)| \leq \| \hat g_n - j\| $. Hence \wrt the $Z$-norm,
$\hat h_n$ converges to $0$; since $h_n \in H$, the second hypothesis
implies $\hat h_n \to 0$ \wrt the usual norm. Thus $\hat k_n \to j$, as
desired.
\qed
\endcomment

It is amusing to ask when other positive maps $\Z^n \to \Z$ are good or
(better, for our purposes, order unit good, since a limit of order unit
good traces is still order unit good, and if the limit group happens to be
simple, the limit trace is then good).
In fact, no others are good, but some others are order unit good.

\Lem Lemma \sixfou. Let $w = (c(i)) \in \Z^{1\times n}$ be a nonnegative row, for which
$\gcd \brcs{c(i)} = 1$. Then the trace $\Z^n \to \Z$ given by $v \mapsto wv$ is
good iff all the nonzero $c(i)$ are $1$.

\Pf We may discard the zero entries, and so reduce to the case wherein all
the $c(i)> 0$. If they are not all equal, by permuting the entries, we may
assume $c(1) < c(2)$. Set $b = (0,1,0,\dots,0)^T$ and $a =
(1,0,\dots,0)^T$, so that $wa = c(1) < c(2) = wb$. However, $b$ is an
atom, so the value of any nonnegative less than $b$ at the trace is zero.
\qed

There is   a characterization of order unit good traces in $\Z^n$, but it is far more complicated. 

\comment
Note however, that if $n =2$, then $w = (1,2)$ is order unit good: denote
the trace $v \mapsto wv$ by $t$. It is sufficient to show that  $(k,1)^T$
and $(1,k)^T$ are $t$-good for all $k$ (since every order unit in $\Z^2$
is a sum of vectors of this form, applying [BeH, 1.3]). Now $t((k,1)^T) = k+2$,
and then $\tau((l,\delta)^T)  = l + 2\delta$ and allowing $0 \leq l \leq
k$ and $\delta\in \brcs{0,1}$, the values of these elements run over all
the integers $\brcs{0,1,2,\dots,k+2}$. Similarly, $t((1,k)^T) = 2k+1$, and
the elements $(\delta,l)^T$ under $t$ run over all the integers in
$\brcs{0,1,2,\dots, 2k+1}$.

However, $w = (1,3)$ is not order unit good: take $a = (2,0)^T$ and $b =
(1,1)^T$. Neither is $w = (2,3)$: same choices for $a$ and $b$. In
general, if $w = (c(i))$ and $b = \One n^T$, then $b$ is $t$-good if and
only if the set of partial sums $\brcs{\sum_{i\in S} c(i)}_{S \subset
\brcs{1,2,\dots, n}}$ contains all the sums of the form $\sum a_i c(i)$
that are less than $\sum c(i)$. For example, if $w = (1,1,\dots,1, n) \in
\Z^n$, then $\One n^T$ is $t$-good, and likely $t$ is order unit good. One
the other hand, with $w = (1,1,\dots, 1,n+1) \in \Z^n$, $t$ is not order
unit good (take $a = (n,0,\dots,0)$).

Let  $G = \Z^n$ with the usual ordering. Then every trace is given by a
row $v=  (a(i))\in (\R^+)^n$, the trace being $\tau_v$ given by $\tau_v (x)
= vx$. Since positive scalar multiples don't affect goodness or order unit
goodness, we can always replace $v$ by $\lambda v$ for any positive real
$v$. We can determine necessary and sufficient conditions under which $\tau_v$ is order unit good. 

Select an order unit $u = (u(i))$ in $\Z^n$; being an order unit simply
means that $u(i) > 0$ for all $i$. For $u$ to be $\tau_v$ order unit good,
necessary and sufficient is that if there exists $w \in \Z^n$ \st $0 < vw
< vu$, then there exists $w'$ in the interval in $\Z^n$ with $0 \leq w'
\leq u$ \st $vw' = vw$. In particular, since the interval in $\Z^n$
bounded above by $u$ is finite (it has at most $\prod (u(i)+1)$ elements),
$u$ being order unit good entails that $\Set{vw}{w \in \Z^n}\cap (0, vu)$
(the latter is usual interval in $\R$) be finite.
 
We make the first observation:  if $v$ has among its ratios of entries
(that is, $a(i)/a(j) \not\in \Q$ for some $i,j$ with $a(j) \neq 0$), then $v$
cannot be  order unit good. Simply note that if an irrational ratio
occurs, then $v\cdot \Z^n$ is a dense subgroup of $\R$, hence its
intersection with any interval is infinite.
 
So we can assume that all the nontrivial ratios in the entries are
rational; we may thus assume (by multiplying by a suitable scalar) that
the entries of $v$ are nonnegative integers and with greatest common
divisor one. We may also assume there are no zero entries by removing them
and correspondingly reducing the $n$.
 
Any permutation of the entries of $v$ does not change the order unit
goodness of $v$, so we may assume the entries are increasing from left to
right. Index the {\it set\/} of distinct integers that appear as $N(1) <
N(2) < \dots < N(r)$; each has corresponding multiplicity, $n(i) > 0$ (so
$\sum N(i)n(i) = n$). Now we show that $u= \pmb 1 = (1,1,\dots,1)$ is
$\tau_v$ order unit good if and only if
 
$$\eqalign{
N(1) & =1, \qquad \text{
 and for all $i > 1$, }\cr
N(i) &\leq 1 + \sum_{k < i} n(k)N(k).\cr
}\tag*$$
 
Since $\gcd{N(i)} = 1$, there exists $w\in \Z^d$ \st $vw = 1$; since $vu >
1$, if $u$ is to be $\tau_v$-order unit good, then $N(1) = 1$. In
addition, if for some $j$, $N(j) >  1 + \sum_{k < j} n(k)N(k)$, then
$N(j)-1$  cannot be expressed as a sum $\sum p(i)N(i)$ with $p(i) \leq
n(i)$, violating order unit goodness. Hence the conditions are necessary
for $u$ to be $\tau_v$-order unit good (and also for $\tau_v$ to be order
unit good).
 
Conversely, suppose the conditions on $v$ hold; we can prove $u = \pmb 1$
is $\tau_v$-order unit good by induction on $r$ (the number of distinct
integers appearing in $v$) and $n(r)$, and this is a simple consequence
about the sums forming an interval.
 
Next we show that if $\pmb 1$ is $\tau_v$-order unit good, then every
order unit $u$ is $\tau_v$-order unit good, i.e., $\tau_v$ is order unit
good. This is by induction on the entries of $u$; it suffices to show
order unit goodness is preserved if we add $1$ to any one of the entries,
and again a simple induction using the interval result works. We thus have
the following.
 
\Lem Proposition \sixfiv. Let $v:= (a(i)) \in( \R^n)^+\setminus\brcs{0}$. Then $\tau_v = v\cdot -$
is an order unit good   trace on $\Z^n$ iff up to multiplication by a positive scalar, $(a(i))$
is a unimodular $n$-tuple of nonnegative integers satisfying (*).
 
The fastest-growing choices are of course (up to reordering) $v =
(1,2,4,8, \dots)$ but we also have
$(1,1,3,3,9,9, 27,27, \dots)$ and $(1,1,1,4,4,4, 16, \dots)$, etc.
 
In particular, $v= (1,1,3)$ will yield an order unit good trace, but
$(1,1,4) $ will not. Up to reordering and scalar multiplication, the only
order unit good traces on $\Z^2$ are given by $(0,1)$, $(1,1)$, and
$(1,2)$.
 
For a different choice of order unit, we can sometimes do better. For
example, if $u = (2,1)^T$ (rather than $(1,1)^T$), we can ask for which
$v$ is $u$ $\tau_v$-order unit good. This time,  $u$ is not permutation
invariant. If $v = (1,3)$, the values of $vu'$ with $0 \leq u' \leq u$
vary over $1 = v\cdot (1,0)^T$, $2 = v\cdot (2,0)$, $3 = v\cdot (0,1)$, $4
= v\cdot (1,1)$, and $5 = v\cdot (2,1)$; as all integer values are
attained with nonnegative elements less than or equal to $u$, $u$ is
indeed $\tau_v$ order unit good, even though $\tau_v$ is not order unit
good.
\endcomment
 
\SecT 7 Introductory section on \ers

As usual,  $\One_s$ denote the  column of size $s$ all of whose entries are $1$. When
$s$ is understood, it may be deleted.
 
Let $G$ be a dimension group (with order unit) that is not simplicial, and $H$ be a rank one
subgroup \st $G^{++} \cap H \neq 0$. Suppose we have an order isomorphism
of $G$ with a limit of maps,
$$
G \iso \lim \Arrow A_n; F_n .F_{n+1},
$$
where $F_n = \Z^{f(n)}$ is the usual simplicially ordered free abelian
group of columns of size $f(n)$, and $A_n$ are $f(n+1) \times f(n)$
matrices with nonnegative integer entries, and suppose in addition, we
have the following properties:
\item{(a)} for all $n$, there exists a (positive) integer $p_{n+1}$ \st
$A_n \One_{f(n)} = p_{n+1} \One_{f(n+1)}$;
\item{(b)} the isomorphism from $G$ to the direct limit sends the subgroup
$H$ to $\cup_n [\One_{f(n)},n]\Z$.
 
We make a couple of observations. Condition (a)
says that each $A_n$ has all of its row sums equal (to $p_{n+1}$); we say
the matrix {\it $A_n$ satisfies \ers\/} when this occurs. Condition (a)
also implies $[\One_{f(n)} ,n]\Z \subseteq [\One_{f(n+1)}]\Z$, so the union
of rank one groups is an ascending union of rank one groups (and thus is
always a group, and rank one). We also note that $\One_{f(n)}$ is an order
unit in $F_n$ and its image under $A_n$ is an order unit in $F_{n+1}$ (by
(a)). Hence $[\One_{f(n)},n]$  is an order unit in the direct limit.
Moreover, if $G_0$ denotes the direct limit, and $H_0$ denotes $\cup_n
[\One_{f(n)},n]\Z$, then $G_0/H_0$ is torsion-free (just observe that if
$kg_0 \in H_0$, then $g_0$ must be represented by an element of the form
$t[\One_{f(n)},n]$). We call the sequence (or $G_0$)  {\it an \ers\
realization of $G$ \wrt $H$\/} when (a) and (b) hold. This of course forces
$G/H$ to be torsion-free and $H \cap G^{++} \neq \brcs{0}$.  Moreover,
$p_{n+1} >1$ for infinitely many $n$, or else the limit is simplicial,
which we forbid; hence $H$ is not cyclic.
 
Sometimes, if $H$ is understood, or we are talking about whether there
exists an $H$ for which an \ers\ realization exists \wrt $H$, we say an
{\it \ers\ realization for $G$ exists.} If the matrix sizes, $\brcs{f(n)}$ are bounded,
then there is a telescoping so that they are all equal, say of size $s$,
and then the matrices have $\One_s$ as a common right eigenvector. In that
case, we say that $G$ has a {\it bounded\/} (or {\it size $s$}) \ers\
realization (\wrt $H$).
 
For example, if as an abelian group, $G \iso U \oplus \Z^k$ where $U
\subseteq \Q$, then there is only one choice for $H$, namely $U$, and an \ers\ realization
also requires that none of the traces kill $U$. If instead the underlying
group of $G$ is $\Z[1/3] \oplus \Z[1/2]$ and the only trace is given by
summing (that is, $(a,b) \mapsto a+b$, so $G$ is a simple dimension group
with unique trace, and the trace has kernel $\brcs{(m,-m)}_{m\in \Z} \iso
\Z$), then there are exactly two choices for $H$, $(\Z[1/2], 0)$ and
$(0,\Z[1/3])$. On the other hand, if $G$ has the same underlying group,
but has   as pure traces the projections on each coordinate, then
$G$ is a simple dimension group with two pure traces, but there are no
candidates for $H$ (so no \ers\ realizations exist for $G$).
 
If $G$ is simple with unique trace $\tau$, the conditions on
$H$ are equivalent to $\tau(H) \neq 0$ (equivalently, since $H$ is rank
one,  $\ker \tau \cap H = \brcs{0}$) and $G/H$ is torsion-free. The last
is a pink herring%
\plainfootnote{*}{not as misdirecting as a red herring.} because for
every rank one subgroup $H_0$ of a torsion free group $J$, there is a
unique rank one subgroup $H$ of $J$ \st $H_0 \subseteq H$ and $J/H$
 is torsion-free.
 
Our  results on \ers\ realizations show that for simple dimension groups
with unique trace, the obvious necessary conditions are sufficient, and we
obtain a bound on the size in terms of the rank. All our dimension groups
are countable.
 
\Lem Theorem \sevtwo. Let $G$ be a simple dimension group with unique trace
$\tau$, together with a noncyclic rank one subgroup $H$ \st $\tau(H) \neq
0$ and $G/H$ is torsion-free.
\item{(a)} If $\rk G = k+1$, then there exists an \ers\ realization of $G$
\wrt $H$ of size $k+2$.
\item{(b)} There exists an \ers\ realization of $G$ \wrt\ $H$.
 
Part (a) (proved in the next section as part of \eigfiv) includes an explicit bound in terms of the rank (which is sharp:
some of these dimension groups cannot be realizedÑeven without the \ers\
propertyÑat the same size as their rank). Part (b) (established in section 9) is a routine
consequence of (a), and of course permits infinite rank (which means that
the $f(n)$ have to be unbounded).

We have a huge class of \ers\ representations available: begin with an \ecs\ realization of a dimension group by square matrices, for example as obtained in \sixone, and take the sequence of transposes. The resulting dimension groups are not that closely related to the original ones from which they emanated. For example, although the dimension group defined by the transposes obtained from the previous construction will have unique trace; generically, {\it this is not rational-valued.} (This will become clear later.)\plainfootnote{$^{2}$}{\rm It is not true in general that if $G$ is a limit of square strictly
positive matrices (so is a simple dimension group) and $G$ has unique
trace, then the limit dimension group of their transposes need have 
unique trace (although it is simple). This is left as an exercise to the reader,
but with a hint: first do it for upper triangular $2 \times 2$ matrices where the number of traces---corresponding to certain
eigenvectors---can easily
be made to change by transposition, then perform a perturbation so the matrices are strictly positive.}

We  have to enter the looking-glass world of torsion-free abelian groups, and as a result, intuition goes out the window. For example, the group $G = \Z[\slfrac12] \oplus \Z[\slfrac13]$ is a simple-minded direct sum of two rank one groups; however, the addition map $\Z[\slfrac12] \oplus \Z[\slfrac13] \to \Z[\slfrac16]$ ($(a,b) \mapsto a+b$) is onto and has kernel isomorphic to $\Z$ (explicitly, $(1,-1)\Z$); hence we have a nonsplit extension of $G$, $\Z \to G \to \Z[\slfrac16]$, by rank one groups, completely different from the direct summands. More generally, if $\brcs{m(i)}_{i=1}^k$ are pairwise relatively prime integers each exceeding one with $m = \prod m(i)$, then $G = \oplus \Z[\slfrac 1{m(i)}]$ is an extension of $\Z^k $ by $\Z[\slfrac 1m]$.

\SecT 8 Transposes

Suppose $J$ is an abelian group, and is given as an extension $0 \to L \to
J \to M \to 0$, with $\Arrow \tau; J.M$ denoting the quotient map. We say
the extension is {\it nearly{} split\/}%
\plainfootnote{*}{Nearly split is almost the same as {\it quasi-split\/} used in abelian group theory (that there exist a map $\Arrow \sigma; M.J$ \st $\tau \sigma$ is $n$ times the identity for some nonnegative integer $n$), and when $J \subseteq \Q$, the definitions coincide. However, quasi-split is also used in other contexts, and I thought it would be confusing here. Different is the notion of {\it almost split\/}, used in representation theory of finite-dimensional algebras. 

In [R], nearly split is defined for extensions of nonabelian groups; it agrees with the definition here when restricted to torsion-free abelian groups. The equivalence classes of nearly split extensions are closed under Baer sums and differences, hence form a subgroup of Ext, although a very small one. We never use the additive structure of the group of extensions. 

} 
if there exists a subgroup $J_0$ of $J$ \st $L \subseteq J_0$, $J_0 = L
\oplus H_0$ for some subgroup $H_0$ of $J$ missing $\ker \tau$ and
$|J/J_0| < \infty$. Equivalently, there exists a subgroup $H_0$ of $J$ \st
$H_0 \cap L = \brcs{0}$ and $\tau(H_0)$ is of finite index in $M$.
 
 In the following, the norms on rows are the maximum of the absolute values, and the norms on matrices are the maximum absolute column sums.

\Lem Lemma \fivsix. Let $\Arrow t; G. V $ be an onto group
homomorphism from a  torsion-free group $G$ of rank $s$  to a dense subgroup $V$ of the reals. Let $H$ be a noncyclic rank one subgroup of
$G$ \st $\ker t \cap H = \brcs{0}$ and $G/H$ is torsion-free. Then there exists a realization of
$G$ as an abelian group, as the direct limit of matrices of the form
$$
\lim \Arrow M_n :=  \( \matrix  p_{n+1} & u^n \\ \pmb 0 & B_n\\
\endmatrix\); \Z^s . \Z^s
$$
 with $p_{n+1} > 1$, $B_n \in
\Z^{(s-1)\times(s-1)}$, $\det B_n \neq 0$,  and $u^n \in \Z^{1 \times({s-1})}$ \st
\item{(i)}  $H \iso
\lim \Arrow \times p_{n+1}; \Z.\Z$
\item{(ii)} $G/H$ is given as $\lim \Arrow
B_n; \Z^{s-1}. \Z^{s-1}$, each $B_n$ of nonzero determinant, and the trace is given up to rational multiple by a sequence of rows of the form $r^i = (1/p_2\dots p_i, \rho_i)$ satisfying $r^{i+1} M_i = r^i$, with $t[a,i] = r^i a$. 
\item{(iii)} The isomorphism of $G$ with the direct limit identifies $H$ with $\cup_{k \in \N} [(1,0,0,\dots,0)^T,k]\Z$.
\item{(iv)}
 $\|B_n\| \leq {p_{n+1}^{1/8s}}/({s!})^{2/s}$ and $\| u^n\| \leq p_{n+1}^{1/4}$.
{\par}
\noindent  Moreover, if $G/H$ is free, then $\ker t$ is free; if additionally, $t(G)$ is rank one, then  the image of $\ker t$   in $G/H$ is of finite index, the extension $\ker t \to G \to t(G)$ is nearly split, and we can take  $B_n = \I_{s-1}$.
 
\Rmk When we change the matrices $B_n$ to the identity, the corresponding
$u^n$ will also change.

\Pf We can write $V$ first as countably generated, say by $\brcs{l_n} \subset\R$, and $t(H) = \cup (1/q_{n+1})\Z$ where $q_n > 1$ divides $q_{n+1}$. and form the subgroups $V_n = (1/q_{n+1})\Z + \sum_{i=1}^n l_i \Z$, so that $V_n \subseteq V_{n+1}$. Next, consider $\ker t$; we can write this as an increasing union of free abelian groups, $J_n \subset J_{n+1}$, all having the same rank as $\rk \ker t = s - \rk V$ (this is true of any finite rank torsion-free abelian group). Select $h_n' \in H$ and $g_n \in G$ \st $t(h_n') = 1/q_{n+1}$ and $t(g_n) = l_n$, and form the  group $G_n$ generated by $\brcs{J_n, h_n', g_1, g_2, \dots, g_n}$; this is finitely generated, hence being a subgroup of a torsion-free group, is free; moreover, its rank must $\rk J_n + \rk V = \rk \ker t + \rk V = s$.  

Then $G_n \subseteq G_{n+1}$, and since $\ker t  = \cup J_n \subset \cup G_n$, and $\cup G_n \to V$ is onto, it follows that $G = \cup G_n$. Now define $H_n = H \cap G_n$; this is cyclic and its image under $t$ contains (possibly strictly) $h_n'\Z$. We may choose its generator, $h_n$, so that $t(h_n) >0$ (which of course uniquely determines it). Since $G/H$ is torsion-free, so is $G_n/H_n = G_n/(h_n \Z)$.  Hence for each $n$, there is an ordered $\Z$-basis whose first entry is $h_n$.

\comment

We can write $U = \cup_{i=1}^j (c_1\cdot c_2 \cdot\dots \cdot c_j)^{-1}\Z$ for
some sequence of positive integers, $c_j > 1$. 

Write $s(j) = c_1\cdot \dots \cdot c_j$, so that $U = \cup \frac 1{s(n)}\Z$. 
Set $G_n' = t^{-1}(1/s(n)\Z)$; since $t(G_n')$ is cyclic, the extension splits, and obviously we can find $g_n \in G_n'$ \st $t(g_n) = 1/s(n)$. Define $H_n = G_n' \cap H$; this is rank one, and since $t(H_n) \subseteq t(G_n')$, each $H_n$ is cyclic; it thus has a  generator $h_n$; then $t(h_n) = z_n/a(n)$ where $(z_n,a(n)) =1$, $a(n)$ divides $s(n)$, and if $h \in H_n$, then $t(h)/t(h_n)$ is an integer.

Since $\ker t$ is a rank $s-1$ torsion-free abelian group, we may find an
increasing chain of free, rank $s-1$ subgroups, $S_n$ \st $\ker t = \cup
S_n$. Define $G_n$ to be the subgroup of $G$ generated by $\brcs{F_n,h_n; g_1, \dots, g_n }$. This is finitely generated, and thus free; moreover, it is contained in $G_n'$ (since $\ker t \subset G_n'$ and $G_i' \subseteq G_{i+1}'$).  Its rank cannot
exceed that of $G$, so must be $s$. There exists no $h \in G_n$ \st $h_n = mh$ for some integer $m > 10$, since this would force $h \in H$, and thus cannot occur in $G_n'$. Hence $h_n \Z$ is a direct summand of the free finitely generated abelian group $G_n$. 

We may thus find, for each $n$, an ordered $\Z$-basis for $G_n$ whose first element is $h_n$. 
\endcomment 

The matrix implementing $G_n \subseteq G_{n+1}$ \wrt the two bases is precisely of the form displayed (but without the estimates in (iv) being satisfied), where $p_{n+1}$ is uniquely determined by $h_n = h_{n+1} p_{n+1}$. Condition (i) is straightforward to verify.  We have seen that $G = \cup G_n$, so we obtain a sequence of matrices whose limit abelian group is $G$. The matrices $B_n$ are the  maps $G_n/h_n\Z \to G_{n+1}/h_{n+1}\Z$, and the limit of these is $G/H$. From the rank conditions, $\rk B_n = s-1$ for almost all $n$, so $\det B_n \neq 0$ for almost all $n$ (and so by deleting an initial segment of the direct limit, we can ensure that $\det B_n \neq 0$ for all $n$). The second part of (ii) just follows from the definitions. Condition (iii) comes from the construction.

Now we want to adjust the sequence in order to arrange that (iv) holds.

Having the original construction of $B_n$ as the quotient maps on $G_n/h_n\Z
\iso \Z^s/h_n \Z$, let $\Arrow f;\N.\N$ be any strictly increasing
function. Define $G^n = G_n + h_{f(n)}\Z$. Then $G^n \subseteq G^{n+1}$,
and $t(h_{f(n)})/t(h_{f(n+1)}) = p_{f(n+1)+1}p_{f(n+1)}\cdot \dots \cdot
p_{f(n) +2}$. In particular, we can take the basis for $G_n$ given by
$(h_n, y_{n,1}, \dots , y_{n,k})$, and observe that $(h_{f(n)},  y_{n,1},
\dots , y_{n,k})$ is a $\Z$-basis for $G^n$. The map $\Arrow M^n; G^n.
G^{n+1}$ \wrt this basis then has its first column simply
$(p_{f(n+1)+1}p_{f(n+1)}\cdot \dots \cdot p_{f(n) +2}, 0, 0,\dots,0)^T$.
Moreover, the induced map $G^n /h_{f(n)}\Z \to G^{n+1}/h_{f(n+1)}\Z$ is
naturally the same as the induced map $G^n /h_n\Z \to G^{n+1}/h_{n+1}\Z$,
that is $B_n$. Hence the form of the transition matrices $M^n$ is
$$
\( \matrix
p_{f(n+1)+1}p_{f(n+1)}\cdot \dots \cdot p_{f(n) +2} & u^n \\
\pmb 0 & B_n \\
\endmatrix\)
$$
for some (different, but relabelled) $u^n \in \Z^{1\times s}$. Thus the new $p_{n+1}$ is the product
$p_{f(n+1)+1}p_{f(n+1)}\cdot \dots \cdot p_{f(n) +2}$, which we can make
as large as we like (by choosing $f$ to grow fast), while fixing $B_n$.
Relabel the upper left corner $p_{n+1}$. Thus we can ensure that $\| B_n
\|^{2s}  \leq \sqrt {p_{n+1}}/(s!)^2$ (or smaller if we like) and $p_n $
increasing.
 
Having this, we can now ensure that $\| u^n \| < p_{n+1}^{1/4}$. Set $U_n
= \( \smallmatrix 1 & y_n \\ \pmb 0 & \I_{s-1} \\ \endsmallmatrix\)$ where $y_n \in \Z^{1\times
(s-1)}$ is to be determined. Each $U_n$ is in $\GL{s,\Z}$ and $U_n^{-1} = 
\( \smallmatrix 1 & -z_n \\ \pmb 0 & \I_{s-1}\endsmallmatrix\)$. Then $\lim \Arrow M^n;
\Z^s. \Z^s$ is isomorphic to $\lim \Arrow U_{n+1}M^n U_n^{-1}; \Z^s.\Z^s$
(via $[a,m] \mapsto [U_n a,m]$). We calculate
$$
U_{n+1}M^n U_n^{-1} = \(\matrix
p_{n+1} &  u^n - p_{n+1}y_n + y_{n+1} B_n \\
\pmb 0 & B_n\\
\endmatrix
\).
$$
Set $y_1 = \pmb 0$. Obviously,  $0 \neq |\det B_n| \leq \| B_n\|^{s-1}
\cdot (s-1)! < p_{n+1}^{1/4}$. Now $B_n^{-1}$ exists as a matrix with
rational entries, and $\det B_n \cdot (B_n)^{-1}$ is simply the adjoint
matrix of $B$, so has integer entries. Let $d_n = |\det B_n|$. Then we
have $\Z^{1 \times (s-1)} d_n B_n^{-1} \subseteq \Z^{1 \times (s-1)}$.
Applying $B_n$, we have $d_n \Z^{1 \times (s-1)} \subseteq \Z^{1\times
(s-1)}B_n$.
 
This means that for any vector  $z \in \Z^{1\times (s-1)}$, we can find $y
\in \Z^{1 \times (s-1)}$ \st $\| z - y B_n\| < d_n$ ($\leq d_n/2$ can be
arranged, but is unnecessary here). Given $y_1, \dots, y_n$, we can thus
find $y_{n+1}$ inductively so that $\| (u^n - p_{n+1}y_n) - y_{n+1}B_n\| <
d_n$. After relabelling $U_{n+1}M^n U_n^{-1}$ to $M_n$, the resulting
upper right corner entry (again called $u^n$) thus satisfies $\| u^n\| <
d_n < p_{n+1}^{1/4}$.

Each (newly relabelled) $h_n$ appears as $[(1,0,\dots,0)^T, n]$ from the $\Z$-basis construction, and since $H = \cup h_n \Z$, the identification with $H$ follows again.

\comment
It is trivial that $G = \cup F_n$; hence we can obtain a group realization
for $G$ along the lines of, but much easier than, the method of Lemma
\oneone, dealing with subsemigroups. For each $F_n$, we may select an
ordered $\Z$-basis whose first element is  (the newly relabelled element
of $H \cap F_n$) $h_n$. Then we express each basis element of $F_{n+1}$ as
$\Z$-linear combinations of the basis elements of $F_n$. Since the ranks
of the abelian groups are all $s$, the corresponding decompositions are
unique, and moreover, from $h_{n+1} = p_{n+1}h_n$, we obtain the matrices,
$$
M_n:= \( \matrix p_{n+1} & w^n \\
\pmb 0 & B_n\\
\endmatrix
\).
$$
Form $G_0 = \lim \Arrow M_n; \Z^s.\Z^s$. There is an obvious map from the
$n$th copy of $\Z^s$ to $F_n$, sending the standard basis to the current
basis for $F_n$. This is compatible with $M_n$ (by construction; we obtain
$M_n^T$ from the coefficients), and so we obtain a group homomorphism $G_0
\to \cup F_n = G$; it is clearly onto, and since the rank of $G_0$ is at
most $s$, and the rank of $G$ is $s$, it must must also be one to one
(since an onto map between torsion-free abelian groups of equal finite
rank is necessarily one to one), hence a group isomorphism.

 The matrix $B_n$ induces the
quotient map $\Arrow B_n ;G_n/h_n\Z.  G_{n+1}/h_{n+1}\Z$, and it is easy to check that the limit of this is naturally is the kernel of $t$. Now suppose that $\ker t$ is free, necessarily of rank $s-1$. It follows that 
 after discarding a
finite number whose determinant is not $\pm 1$, all the $B_n$ must be in
$\GL{s-1,\Z}$.

Now define $R_n = \( \smallmatrix 1&  \pmb 0\\ \pmb 0 &B_n \\
\endsmallmatrix\)$, and define $M'_1 = B_1 R_1^{-1}$, and $M'_n = R_1 R_2
\cdots R_{n-1}B_n (R_1 \cdots R_{n-1}R_{n})^{-1}$. The 
the direct limit of the $M_n$ is isomorphic to the original, and now the
lower right block is the $(s-1) \times (s-1)$ identity,  the lower left block
is zero, and the upper left entry is still $p_{n+1}$. This is the form we required, and moreover, the matrices have the same common right eigenvector $(1,0,\dots,0)^T$, 
the isomorphism of
abelian groups preserves the map to $U$.
\endcomment

Now we deal with the {\it Moreover\/} statement. The map $\ker t \to G/H$ is one to one; so if $G/H$ is free, then $\ker t$, being a subgroup, is free as well. Since $G$ has finite rank, $G/H$ is free of rank $s-1$. If additionally, $t(G)$ has rank one, then   $\ker t$ has rank $s-1$, the same as that of $G/H$, and since both are free, the image $\ker t$ is of finite index in $G/H$. 

Since $G/H \iso \lim \Arrow B_n; \Z^{s-1}. \Z^{s-1}$ ($\Z^{s-1}$ is an abbreviation for $\Z^s/h_n\Z$), and $G/H$ is free of maximal rank, it must happen that $|\det B_n| = 1$ for all but finitely many $n$ (from finite generation of the direct limit). If $G/H$ is free, then $\ker t \oplus H$ is of finite index in $G$: to see this, note that $G \to G/H$ splits, so there exists a subgroup $J$ of $G$ \st $H \oplus J = G$ and $J$ maps isomorphically to $G/H$. There is no guarantee that $\ker t \subseteq J$; however, the exact sequence $H \to H \oplus \ker t \to L$ (where $L$ is the image of $\ker t$ in $G/H$) yields $G/(H \oplus \ker t)$ is finite, since it embeds in $(G/H)/L$, which is finite.

Still in the case that $G/H$ is free, we may discard an initial segment of non-elements of $\GL{s-1,\Z}$, so assume each $B_n$ is in $\GL{s-1,\Z}$. Then we can systematically pre- and post-multiply the $M_n$ by matrices of the form $\diag (1 ,E_n)$ where $E_n \in \GL{s-1,\Z}$ to arrange that the lower right blocks are all the identity.
\qed
 
For the general case, the matrices $B_n$ can be put in Hermite normal form (the normal forms arising from the action of $\GL{s-1,\Z}$ on $\Z^{(s-1) \times (s-1)}$ from the left). It is not clear whether this would be useful.

An immediate observation is that $e:= (1,0,\dots,0)^T$ is a common right
eigenvector for all the matrices $M_n$ appearing there, with eigenvalue
$p_{n+1}$, and if we identify $G$ with the direct limit, then $H =
\cup [e,k] \Z$. We can also recalculate $t$  in terms of the direct limit.
 
\comment 
Every map from the direct limit to the reals is given by a sequence of
rows, $\rho_j \in \R^{1 \times s}$ (literally rows = rhos) satisfying
$\rho_j a = \rho_j M_j a$ for all $a\in \Z^s$, and the corresponding map
is given by $t(a,j) = \rho_j a$. Here the image is rational, so all
entries of the $\rho_j$ must be rational, but more is true. Since $t(F_n)
= (c_1\dots c_n)^{-1}\Z$, we can assume that the denominators of all the
entries of $\rho_j$  divide $c_1\dots c_n$, and that $1/c_1\dots c_n$
itself  can be realized by something (actually $e_n$). This means we can
write $\rho_n = r^n(c_1\cdot c_n)$ where $r^n$ is a unimodular element of
$\Z^{1\times s}$.
 
The equations then become $r^n c_{n+1} = r^{n+1}M_n$, and $t[a,j] = r^j
a/c_1\dots c_j$. Applying this to $1/p_1\dots p_n = t[v,n] = r^n
v/c_1\dots c_n $, we deduce $r^n_1 = c_n\dots c_1/p_n \dots p_1$ (the
first coordinate of $r^n$). The group homomorphism $t$ is thus given by a sequence of rows, $r^j \in \Z^{1 \times m}$, $t [a, j] = r^j a/(c_1\dots c_{j})$ \st the following diagram commutes. 
$$\diagram
\Z^{m}&\rTo^{M_j} & \Z^{m} \\
\dTo^{r^j} &&\dTo^{r^{j+1}} \\
\Z&\rTo^{\times c_{j+1}} & \Z.\\
\enddiagram
$$
\endcomment

\Lem Lemma \fivone. Let $\Arrow B_i ;\Z^d.\Z^d$ be a bunch of matrices, and let
$J$ be their limit as an abelian group. Suppose that for all $i$, the left
kernel of $B_i$, that is, $\Set{w \in \Z^{1\times l}}{wB_i = \pmb 0}$, is the
same, $\Z z$, for some $z \in \Z^{1 \times l}$; we may assume that $z$ is
unimodular. Set $W = z^{\perp} = \Set{v \in \Z^l}{zv = 0}$; then $B_i W
\subseteq W$ and form the direct limit, $\lim J_0:= \Arrow C_i;W.W$, where
$C_i = B_i|W$. Then the natural  map  $J_0 \to J$ given by $[v,s]_W
\mapsto [v,s]$, is an isomorphism (of abelian groups).
 
\Pf Since $z (B_i W) = 0$, not only is $B_i W \subseteq W$, but in fact
$B_i (\Z^l) \subset W$.
If $B_{n+t}\cdot B_{n+t-1}\cdots B_{n+1} v = \pmb 0$ for $v$ as an element
of $W$, then it is obviously true as an element of $\Z^l$, and it follows
that the map $J_0 \to J$ is well defined and one to one. Next, if $y \in
\Z^l$, then $B_{s}y \in W$, so that $[y,s] = [B_s y , s+1]$ which is in
the image of the map $J_0 \to J$. Hence the map is onto.
\qed

In the \ecs\ cases discussed in the previous section, the sequence  of vectors $(v^n)$ is compatible with
the addition operation on the \text{Ext\,} group, that is, with the Baer
sum ($(p_{n+1}, v^n  + (v^n)')$ represents the Baer sum of the extensions
arising from $(p_{n+1}, v^n)$ and $(p_{n+1}, (v^n)')$); however, many
different sequences can represent the same equivalence class, and it is
very difficult to decide when they do. The same applies here, although if $G/H$ is free, then as abelian groups (but not as extensions), $G \iso \Z^{s-1} \times H$. 

\comment
The obvious combination of the previous yields the desired result on \ers\ realizations. Take the matrix $A$ as given in (1) and realized as a primitive matrix via Lemmas \fouthr\ and \foufou, and set $z = (-(k+1),1,1,\dots,1) \in
\Z^{1 \times k+2}$. Then $W:= {}^\perp z$ is a common invariant subset
(regardless of the choice of $n$), and we pick a $\Z$-basis for it, $e_0 =
(1,1,\dots, 1)$; $e_i = (0,0,\dots, 1, 0, \dots, -1)$ with $1 \leq i \leq
k$ (we index the coordinates of $\Z^{1 \times (k+2)}$ as usual by $0, 1,2, \dots,
k+1$; with this indexing, $e_i$ has its leading $1$ in position $i$). It
is easy to see that this is a basis for ${}^{\perp}z$. We observe that
$e_0 A = p e_0$ and $e_i A = v_i e_0 + e_i$.
 
Transposing everything in sight, let $f_i = e_i^T$, and note that
$(z^T)^{\perp}$ is spanned by $\brcs{f_i}$, and is invariant under every
$A_n$. Moreover, the matrix of $A$ \wrt this common basis is $\(
\smallmatrix p&  v^T\\ \pmb 0 & \I_k \\ \endsmallmatrix\)$. Since the
basis is the same for all $A_n$, we have as abelian groups, $\lim \Arrow
A_n^T; \Z^{k+2}.\Z^{k+2} \iso \lim \Arrow \( \smallmatrix p_{n+1}& (
v^n)^T\\ \pmb 0 & \I_k \\ \endsmallmatrix\); \Z^{k+2}.\Z^{k+2}$.  

Hence as abelian groups, $H \iso G$; moreover, the inclusion map $J_0 \to
J$ preserves the map to $U$ (given by the the sequence of left
eigenvectors of $A_n^T$), so that the induced automorphism of $U$ is
either plus or minus a rational scalar multiple of the identity (and if
$U$ has no infinitely divisible primes, it can only be $\pm I$); then it
is easy to check that either the map $ H \to G$ or its negative is an
order isomorphism.\qed

\Lem Proposition \fivfiv. Let $(G,\tau)$ be a simple finite rank dimension group with unique trace, \st the latter is rational-valued, and $\ker \tau$ is free. Then $G$ admits an \ers\ representation if and only if $G$ contains a noncyclic rank one subgroup. When this holds, $G$ is realizable by a sequence of transposes of primitive matrices of the form  (1), and this is \ers\ of size $k+2$.
 
\endcomment

%%%%%%

Here $\rho$ denotes the spectral radius. The following is well known in a more general setting, dealing with projective convergence and weak ergodicity. But we do not need this generality in our situation.
 
\Lem Lemma \fivthr. Let $G = \lim \Arrow C_i; \Z^s.\Z^s$ be a sequence of primitive matrices
for which there exists positive real numbers $f(m,n)$ with $m > n$ \st
$\lim_{m\to \infty} f(m,k+1)/f(m,k) \to 1/\rho(C_k)$ for all $k$,  and for
all $n$
$$
\lim_{m \to \infty \ \&\ m > n}\frac{C_m C_{m-1} \dots C_n}{f(m,n)} = V_n
$$
exists and is nonzero. Then the candidate map $\Arrow V;G.\R^s$ via
$V[a,k] = V_k/\prod_{i=1}^{k-1} \rho(C_i)$ is well-defined, and every pure
trace of $G$ factors through it. In particular, $G$ has unique trace iff
$\rk V_n = 1$ for almost all $n$.
 
\Rmk The simplest situation in which the hypotheses hold occur when
$\rho(C_m C_{n-1}\dots C_m) = \prod_m^n \rho (C_i)$ for all $m > n$, that
is, when the spectral radius is multiplicative on the matrices. For
example, this occurs when the $C_n$ have a common right Perron
eigenvector, or a common left Perron eigenvector.
 
\Pf $[a,k] = [C_k a,{k+1}]$, and the latter is sent to $V_{k+1} C_k
a/\prod_{i=1}^{k}\rho(C_i)$. Now $C_m \dots C_{k+1} C_k/f(m,k) = (C_m
\dots C_{k+1}/f(m,k+1)) C_k (f(m,k+1)/ f(m,k))$. The left side converges
to $V_k$; the right side converges to $V_{k+1}C_k/\rho(C_k)$. Hence
$V_{k+1} C_k a = \rho(C_k)V_k a$, so $V$ is well-defined.
 
Next, suppose that $[a,k]$ is an order unit in $G$; then there exists $m>
k$ \st
$C_{m-1} C_{m-2} \dots C_k a $ is strictly positive; as $V_m$ has only
nonzero entries, this means $V_m C_{m-1} C_{m-2} \dots C_k a $ is
nonnegative and not all entries are zero, and thus $V[a,k] $ is
nonnegative and nonzero (as an element of $\R^s$), and thus $V$ is a
positive group homomorphism. Each row of $V_n$ is either zero, or induces
a trace on $G$ (via $\R^s \to \R$).  Discard any zero rows from $V_n$ (for
all sufficiently large $n$) obtaining a  newly-labelled $V_n$ which  is
now a map from $G$ to $\R^{s'}$ with $s \leq s'$ \st every row of $V_n$ is
not zero. Then the map $V$ sends order units if $G$ to order units of
$\R^{s'}$.
 
Conversely, if $[a,k] $ is an arbitrary element of $G$ \st $V[a,k] > 0$,
then there exists $m > k$ \st $\| (C_{m-1} C_{m-2} \dots C_k)/f(m,k) a  -
V_k a \| $ is smaller than the infimum of the entries of $V_{k}a$, and
thus   $(C_{m-1} C_{m-2} \dots C_k)/f(m,k) a$ is strictly positive, hence
$C_{m-1} C_{m-2} \dots C_k a$ is strictly positive, and thus $[a,k]$ is an
order unit of $G$. Now consider all the traces on $G$ obtained by
composing $V$ with any positive vector space map  $\R^s \to \R$.
What we just obtained is that these are enough to determine the order
units of $G$, and this implies that these traces include all the extreme
points in the trace space of $G$, hence the factorization for pure traces.
 
If $\rk V_n =1$, then the trace space is $0$-dimensional (after
normalization, a single point); conversely, if $G$ has unique trace, then
all the composed traces are equal up to normalization, and it follows
immediately that $\rk V_n = 1$ for almost all $n$.
\qed
 
Suppose $A_i$ are primitive matrices of the same size with common right
Perron eigenvector. Then the spectral radius is multiplicative on products
of the $A_i$, and moreover, $A_i /\rho(A_i)$ are uniformly bounded (by row
sum) by $1$, as are their products. Hence there exists a subsequence, $1 =
n(1) < n(2) < n(3) \dots$, \st  for
the sequence $\(C_i := A_{n(i+1)-1}   \cdots A_{n(1)+1}\cdot A_{n(1)}\)$, we
have for all $k$,
$$
\frac  {C_m \cdot C_{m-1} \cdot \dots \cdot C_k}{\prod_{i=k}^m \rho(C_i)}
$$
converges to a matrix, necessarily nonzero, as the row sums are all one.
Hence by suitably telescoping, we use \eigthr\ to derive the pure traces from   rows of
the limit matrices, and if $G$ has unique trace, the limit matrices
eventually have rank one, so we can pick any fixed row.
 
There is a more general uniqueness criterion in terms of projective
convergence that is very well known, applies to all unique trace dimension
groups and corresponding Bratteli diagrams, but is not easy to  use.

 \comment
\Lem Lemma \fivthr. Let $G$ be a nearly split extension of $\Z^k$ by $U$. For every sequence  $\brcs{(p_{n+1}}$ \st $U =
\lim \Arrow \times p_{n+1}; \Z.\Z$, there exists a telescoping $(q_{j+1} =
\prod_{i = n(j)}^{n(j+1)-1} p_{i+1})$ together with $w^n \in \Z^{1 \times
k}$ \st
$\lim \Arrow M_n:= \( \smallmatrix q_{n+1}&  w^n\\ \pmb 0 & \I_k \\
\endsmallmatrix\); \Z^{k+1}.\Z^{k+1}$ yields the extension (up to
isomorphism of groups with map to $\R$).
 
\Pf Let $H_0$ be a subgroup of $G$ \st $\Arrow t; G.U$ maps $H_0$ to a
finite index subgroup of $U$ and $H_0 \cap \Z^{k} = \brcs{\pmb 0}$ (the
existence of such is the definition of nearly split). Since $G/(H_0 \oplus
\Z^{k})$ is finite, there are only finitely many direct sums between $H_0
\oplus \Z^k$ and $G$. Pick a maximal one, i.e., we have $H \subset G$,
$|U/t(H)| < \infty$, $H \cap \Z^k = \brcs{\pmb 0}$ \st if $H \oplus \Z^k
\subset H' \oplus \Z^k \subset G$, then $H \oplus \Z^k = H' \oplus \Z^k$.
 
Now any subgroup of $U$ of finite index is of the form $mU$ (just look at
the supernatural number of $U$ for some integer $m$. By tossing away
finitely many primes (with finite multiplicities), we can assume there
exists $h \in H$ \st $H = U\cdot h$; that is, there exist $h_n$ \st $h_n =
p_{n+1}h_{n+1}$ with $h_1 = h$, and the map $t$ restricted to $H$ sends
$h$ to $m \in mU$ (where we identify $U$ with the subgroup of the
rationals given by $\Set{a/b}{a \in \Z \text{ and } b = \prod_{i=1}^n
p_{i+1}}$).
 
For each $n$, consider $F_n = h_{n} \Z \oplus \Z^k $, and set $G_n =
\Set{g \in G}{mg \in F_n}$. Then $G_n$ is free of rank $k+1$, and we claim
that, possibly after telescoping, we can assume that $h_n$ is not a proper
multiple of anything in $G_n$. If $h_n = sg$ for some integer $s > 1$,
then we can consider the subgroup of $G$ which is essentially $Ug$; this
maps to a larger finite index subgroup, contradicting maximality of $H
\oplus \Z^k$. Hence $h_n \Z$ is a direct summand of $G_n$, and so we can
find an ordered $\Z$-basis for $G_n$ of the form $(h_n; e^n_1, \dots,
e^n_k)$.
 
When we express the generators of $G_n$ in terms of the generators of
$G_{n+1}$, and so deduce a sequence of matrices whose direct limit is $G$
(as an abelian group; the map to $U$ takes care of itself, since the first
basis element is always sent to $p_{n+1}$ times something), we realize $G$
as the direct limit of the matrices, in block form,  $C_n :=  \(
\smallmatrix p_{n+1}&  w_0^n\\ \pmb 0 &E_k \\ \endsmallmatrix\)$ where
$w_0^n \in \Z^{1\times k}$ and $E_k$ (a $k \times k$ matrix induces the
quotient map $\Arrow E_n ;G_n/h_n\Z.  G_{n+1}/h_{n+1}\Z$; since the
matrices have a common left eigenvector implementing the map to $U$, the
limit of the $E_n$ must be isomorphic to $\Z^k$. Hence after discarding a
finite number whose determinant is not $\pm 1$, all the $E_n$ must be in
$\GL{k,\Z}$.
 
Now define $R_n = \( \smallmatrix 1&  \pmb 0\\ \pmb 0 &E_n \\
\endsmallmatrix\)$, and define $M_1 = C_1 R_1^{-1}$, and $M_n = R_1 R_2
\cdots R_{n-1}C_n (R_1 \cdots R_{n-1}R_{n})^{-1}$. It is easy to see that
the direct limit of the $M_n$ is isomorphic to the original, and now the
lower right block is the $k \times k$ identity, and the lower left block
is zero. This is the form we required, and moreover, since they have a
common left eigenvector implementing the map to $U$, the isomorphism of
abelian groups preserves the map to $U$.
\qed

\endcomment

Denote by $\Cal B(p,B, v,a)$ for $p$ and $a$ positive integers, $v \in \Z^{k}$ and $B \in \Z^{k \times k}$, the matrix
$$
\Cal B(p,B, v,a) = \(\matrix
* & * &* \\
v + a\pmb 1 & B + v \pmb 1^T + a\pmb 1\pmb 1^T & v  + (a-1)\pmb 1 \\
a & a \pmb 1^T & a \\
\endmatrix\)
$$
where $\pmb 1$ is the column of size $k$ consisting of ones, and the column sums are all $p$ (hence the entries marked with an asterisk are uniquely determined). Then it is easy to check that $\Cal B (p,B,v,a) \Cal B(p',B',v',a') = \Cal B(pp', BB', p'v + Bv,p'a)$. Setting $A_n = \Cal B (p_{n+1},B_n,v^n,A_n) $, then inductively
$$
A_n A_{n+1} \cdots A_{n+j} = \Cal B (\prod_{i=0}^j p_{n+i+1}, B_n \cdots B_{n+j}, v^{(n,j)}, A_n \prod_{i=1}^j p_{n+i+1}),
$$
where
$$
v^{(n,j)} = p_{n+1}\dots p_{n+j+1} \(\frac {v^n}{p_{n+1}} + \frac{B_n v^{n+1}}{p_{n+1}p_{n+2}} +   \frac{B_n B_{n+1}v^{n+2}}{p_{n+1}p_{n+2}p_{n+3}} + \dots + \frac{B_n \dots B_{n+j}}{p_{n+1}\dots p_{n+j+1}} \).
$$

We assume as we may that $\| B_n \| = \Oh{p_{n+1}^{1/2}}$, and that $\brcs{v^n/p_{n+1}}$ is bounded. Then $\lim_{j\to \infty} v^{(n,j)}/p_{n+1}\dots p_{n+ j+1}$ exists (provided $p_n \to \infty$); call it $V^{\infty,n}$; this forces $A_{n+j}^T A_{n+j-1}^T\dots A_n/p_{n+1}\dots p_{n+j+1}$ to converge to the rank one matrix
$$
\pmb 1_{k+2} \( 1 - (V^{\infty,n})^T\pmb 1 - \frac {(k+1)a^n}{p_{n+1}},V^{\infty,n} + \frac {a^n}{p_{n+1}}\pmb 1^T, \frac {a^n}{p_{n+1}}\).
$$
Call the row appearing in this factorization, $W^{\infty,n}$; if the $A_n$ are primitive, then it is strictly positive. The family $\brcs{W^{\infty,n}}$ satisfies $W^{\infty,n}A_n^T = W^{\infty,n-1}$, hence induces a trace on the dimension group $G = \lim \Arrow A_n^T; \Z^{k+2}. \Z^{k+2}$ via $\tau[x,m] = W^{\infty, m}x$. As $G$ has unique trace, this is it, up to scalar multiple.

Now we consider $G$ as an abelian group with trace; then we obtain a group isomorphism from the restriction to $(z^T)^{\perp}$ (where $z = (-(k+1),1,1,\dots, 1)^T$); using as ordered $\Z$-basis for the latter, the columns $\((1,\dots,1)^T, (0,1,0,\dots,0,-1)^T, \dots , (0,0,\dots,0,1,-1)^T\)$ (it is easy to check that this {\it is\/} a $\Z$-basis), the group isomorphism from \fivone\  is with the group given as $J:= \lim M_n = \Arrow\( \smallmatrix p_{n+1}& (v^n)^T \\
0 & B_n^T \\ \endsmallmatrix\); \Z^{k+1}. \Z^{k+1}$. Moreover, the effect of $W^{\infty,1}$ on the basis yields the group homomorphism obtained from the rows $R_n:= (1, V^{\infty,1})$; that is, the corresponding homomorphism from $J$ to $\R$ is given by $t[x,k] = R_k x$. Since each $M_n$ is one to one,  a group homomorphism from $J$ is uniquely determined by its affect on the first level, that is, on elements of the form $[x,1]$.

%%%%%%%%%%%

 \comment
Suppose $J$ is given as an abelian group direct limit $\lim \Arrow M_n ;
\Z^s.\Z^s$,
where $M_n = \(\smallmatrix p_{n+1} & w^n\\ 0 & B_n\endsmallmatrix\)$. Set
$e = (1,0,0,\dots,0)^T$ and $u = [e,1]$ and let $H(u)$ be the largest rank
one subgroup of $J$ containing $u$, that is, $\cup [e,j]\Z$ (noting that
$[e,1] = p_2 [e,2]$, etc). Then $G/H(u)$ is torsion-free, and in fact, is
obtainable as the direct limit
$\Arrow {B_n}; \Z^{s-1} . \Z^{s-1}$: since $e\Z$ is invariant for all the
$M_n$, each of the latter induces a homomorphism $\overline{M_n}; \Z^s/e\Z
\to \Z^s/e\Z$ given by $\overline {M_n} (v + e\Z) = M_n v + e\Z$; if we
write $v = v_1 e + (0,v')$, then $M_n v = (p_{n+1} v_1  + w^n v')e + B_n
v'$. Hence $\overline{M_n} $ on the quotient group is just $B_n$ on the
obvious copy of $\Z^{s-1}$.

The group homomorphism $t$ is given by a sequence of rows, $r^j \in \Q^{1
\times s}$ satisfying $r^j = r^{j+1}M_j$ (since $t$ is rational-valued,
the entries can be chosen in $\Q$; for general group homomorphisms to the
reals of course, the entries will be real);  $t[v,k] = r^k v$. Consider
the subgroup $(r^j)^{\perp} = \Set {x \in \Z^s}{r^j x = 0}$ of $\Z^s$. If
we write $r^j = (r^j_1, \rho^j)$ (so $\rho^j$ consists of the remaining
entries of $r^j$), the compatibility condition $r^j = r^{j+1}M_j$
means that  $r^{j}_1 = p_{j+1}r^{j+1}_1$ and $\rho^j = r^{j}_1 w^j
/p_{j+1}+ \rho^{j+1}B_j$.
In particular, $M_j((r^j)^{\perp}) \subseteq (r^{j+1})^{\perp}$. Since all
the entries of all the $r^j$ are rational, it also means $\rk
(r^j)^{\perp} = s-1$, and since $M_n$ is one to one, it follows that the
restriction of $M_j$ to this subgroup is one to one.
 \endcomment
\comment
For $v \in \Z^s$, $M^j v = (p_{n+1}v_1+ w^j v') e + B_jv'$; if $r^j v =
0$, then $r^j_1 v_1 + \rho^j v' = 0$, so $v_1 = -\rho^j v' /r^j_1$ (the
latter is not zero if $t(u) \neq 0$, which we assume). So $M_j v =
(-p_{j+1}\rho^j v'/r^j_1 + w^j v')e + B_j v'$; then the coefficient of $e$
simplifies to $(-p_{j+1}/r^j_1)(r^j_1/p_{n+1}) w^j v' + w^j v' = 0$. So
$M_j|(r^j)^{\perp} = B_j \pi$, where $\pi$ is the projection from $\Z^s $
deleting the top entry. NOOOOOO$\ker t = \lim \Arrow B_j ;
\pi(\Z^s) . \pi( \Z^{s})$. (It is tempting to conjecture that this is
isomorphic to $\lim \Arrow B_j; \Z^{s-1}. \Z^{s-1}$, but I'm pretty sure
that is not trueÑbut each maps to a finite index subgroup of the
other(?).)
 
Get embedding from limit of $B$s restricted to action of $B$s on quotient
group. If limit Bs is free or if it is extension of $\Z^{s-1}$ by
$\Z[r^{-1}]$ for some prime, get that $\ker t + H$ has finite index in
$G$, hence nearly split. And in the former case, the other extension $G
\to G/H$ splits, so $G$ is a direct sum, but not a trivial extension!
 
Explain with 3101 matrix.
\endcomment
 
%%%%%%%%%
 
In the following, the norms on rows are the maximum of the absolute values, and the norms on matrices are the maximum absolute column sums. 
 
\Lem Lemma \eigfou. Let $p_{n+1} \uparrow \infty$, let $B_i $ be $k\times k$ integer matrices \st $\det B_i \neq 0$ and $\| B_i\| = \oh{p_i^{1/2}}$, and let $z_i \in \Z^{1\times k}$ with $\| z_i \| < p_{i+1}$; let $r^1 = (1, \rho^1)$ where
$\rho^1 \in \R^{1 \times k}$. Then there exist a sequence $\brcs{w^i}$,
with $w^i \in \Z^k$ and $\| w^i\| < \|z_i\| + (p_{i+1} +  \|B_i\|)/2$ for all $i > 1$
together with group isomorphisms $\Arrow F_i; \Z^{k+1}. \Z^{k+1}$ \st the
following diagram 
$$\diagram
{}&&\Z^{k+1}&\rTo^{\(\smallmatrix p_2 & z_1 \\ 0 & B_1\\
\endsmallmatrix\)} & \Z^{k+1} & \rTo^{\(\smallmatrix p_3 & z_2 \\ 0 & B_2\\
\endsmallmatrix\)} & \Z^{k+1} &{\cdots}&\Z^{k+1 } &\rTo{M_n:= \(\smallmatrix p_{n+1} & z_n \\ 0 & B_n\\
\endsmallmatrix\)}& \Z^{k+1}& \rTo^{\dots}\\
{} &\swarrow {r^1} &&& && && &&&\\
\R \hskip-10pt& &\dTo^{F_1}&&\dTo^{F_2}&&\dTo{F_3}&&\dTo{F_n}&&\dTo{F_{n+1}}\\
{} & \hskip0pt\nwarrow r^1\\
&&\Z^{k+1} &\rTo^{\(\smallmatrix p_2 &w_1 \\ 0 & B_1\\
\endsmallmatrix\)}& \Z^{k+1}&\rTo^{\(\smallmatrix p_3 &w_2 \\ 0 & B_2\\
\endsmallmatrix\)}& \Z^{k+1}& \dots &\Z^{k+1} &\rTo^{\(\smallmatrix p_{n+1} &w_n \\ 0 & B_n \\
\endsmallmatrix\)}& \Z^{k+1}&\rTo^{\dots}\\\
\\
\enddiagram
$$
commutes,  and \st
$$
\rho^1  = \frac{w^1}{p_2} +  \frac{w^2B_1}{p_3p_2} +
\frac{w^3B_2B_1}{p_4p_3p_2} + \frac{w^4 B_3 B_2B_1}{p_5 p_4p_3p_2} + \dots. \tag2
$$
 
\Pf We will define $F_i = \(\smallmatrix 1 & y_i \\ 0 & \I_k\\
\endsmallmatrix\)$ (where $y_i  \in \Z^{1\times k}$), and then define $w_i \in \Z^{1 \times k}$ so that all the
properties hold. First set $y_1 = \pmb 0$. Now define 
$$\eqalign{
z^{\infty} &= \frac {z_1}{p_2} + \frac {z_2 B_1}{p_3p_2} + \frac{z_3B_2 B_1}{p_4p_3p_2} + \dots\cr
& = \frac{z_1}{p_2} + \sum_{i=2}^{\infty} \frac{z_i B_{i-1}B_{i-2}\dots B_1}{p_{i+1}p_i\cdots p_2}.
}$$
That the sum exists is a consequence of $\| z_i\|/p_{i+1}$ being bounded, $\|B_i\| = \oh{\sqrt{p_{i+1}}}$, and summability of $1/\sqrt {p_ip_{i-1}\dots p_2}$. 

Now define 
$$
y_{n+1} = \left[\vphantom{\int_o^{o}}p_{n+1}p_n\cdot \dots \cdot p_2 (\rho^1 - z^{\infty})(B_n\cdot B_{n-1}\cdot \dots \cdot B_1)^{-1}\right];
$$
of course, the inverses of $B_i$ exist as matrices with rational entries. Here the integer function $\left[\vphantom{X^i}\,\cdot\,\right]$ means to take the nearest integer  in each entry. Let $Y_{n+1}$ denote the thing on the right before we take the integer function; it is an element of $\R^{1\times k}$. Then obviously we have $y_{n+1} \in \Z^{1\times k}$ and $\| y_{n+1} - Y_{n+1}\| \leq \slfrac12$.

Finally, set $w^n = z_n + y_{n+1}B_n - p_{n+1} y_n$. It is easy to check that the squares in the diagram all commute. We show that $\rho^1$ is  the  infinite sum in (2). 

Let $S_n $ be the sum of the first $n$ terms on the right side of (2). When we substitute $w^{i} = z_i + y_{i+1}B_n - p_{n+1} y_i$, we find that the series partially telescopes:
$$
S_n = \( \frac{z_1}{p_2} + \sum_{i=2}^{n} \frac{z_i B_{i-1}B_{i-2}\dots B_1}{p_{i+1}p_i\cdots p_2} \) + \frac{y_{n+1}B_n B_{n-1}\cdots B_1}{p_{n+1}\cdots p_1},
$$
as follows immediately by induction. Now $Y_{n+1}B_n\cdot \dots \cdot B_1(p_{n+1}\dots p_2)^{-1} = \rho^1 - z^{\infty}$, hence $\| y_{n+1}B_n\cdot\dots \cdot B_1(p_{n+1}\dots p_2)^{-1} -(\rho^1 - z^{\infty}) \| < 1/\sqrt{p_{n+1}\cdot\dots \cdot p_1}$. Thus $\lim S_n$ exists and 
$$
\lim S_n = z^{\infty} + (\rho^1 - z^{\infty}) = \rho^1.
$$ 

Next, we estimate $\|w^i/p_{i+1}\|$. We have 
$$\eqalign{
\left\| \frac{w^n - z_n}{p_{n+1}} \right\| &=  \left\| \frac{y_{n+1} B_n - p_{n+1}y_n}{p_{n+1}}\right\|\cr
& \leq \left\| \frac{Y_{n+1}B_n - p_{n+1}Y_n}{p_{n+1}}\right\| + \left\| \frac{(Y_{n+1}- y_{n+1})B_n + p_{n+1}(Y_n - y_n)}{p_{n+1}}\right\| \cr
& \leq  0 + \left\| \frac{ B_n}{2p_{n+1}}\right\| + \frac12 \cr
}$$
Thus $\| w^n/p_{n+1} \| \leq \| z_n/p_{n+1} \|  + (\slfrac12) \( \|B_n\|/p_{n+1} + 1\)$.
\qed

\comment 
From $w^{\infty} = \lambda{-1}  (-y_1 + \rho_1)C^{-1}$, hitting it with
$\pmb 1$, we have $w^{\infty}\pmb 1 = \lambda^{-1} (-y_1 +
\rho_1)C^{-1})\pmb 1= \lambda^{-1} t$. From $\lambda^2 - \lambda -t = 0$,
we have $\lambda^{-1}t = \lambda - 1$. Hence $\lambda= 1 - w^{\infty}\pmb
1$. Thus $ (-y_1 + \rho_1)C^{-1} = w^{\infty}/(1-w^{\infty}\pmb 1)$.
\endcomment

We are permitted to telescope the bottom row, and then apply the same transformation to the resulting upper right corner entries as we did in the \ecs\ case (for the lower left corners), conjugating by a block upper triangular element of $\GL{k+2,\Z}$, to ensure we could choose the $a^n$ so that the resulting matrices are positive.
 This yields the following.

\Lem Theorem \eigfiv. Let $G$ be a simple dimension group of rank $k+1$ with unique trace $\tau$, and let $H$ be a noncyclic rank one subgroup of $G$ \st $G/H$ is torsion-free and $\tau(H) \neq 0$. Then there exists an \ers\ realization of $G$ of size $k+2$ \st the image of $H$ in the direct limit is $\cup_{j\in \Z} [(1,1,\dots,1)^T,j]\Z$. If $G/H$ is free, then the extension $0 \to \ker \tau \to G \to \tau(G) \to 0$ is nearly split.

This does not require the trace to be rational-valued; since there is no restriction on $\tau(G) $ except $\tau(H) \neq 0$, the value group, $\tau(G) $ can be an arbitrary subgroup of $\R$ containing $\tau(H)$ and of rank at most $k+1$ (when equality occurs, $G$ is totally ordered). 
 
%%%%%%%
\comment
For  integer $p > 1$, $v \in \Z^k$, and $a\in \Z$, and $\pmb 1$ the column
of size $k$ with all entries one,  let $\Cal A (p,v,a)$ be the square
matrix of size $k+1$ given by
$$
\Cal A(p,v,a) = \(\matrix
* & * &* \\
a\pmb 1 + v & a\pmb 1 \pmb 1^T + v\pmb 1^T + \I_k & (a-1)\pmb 1 + v \\
a & a\pmb 1^T, a\\
\endmatrix\),
$$
where the asterisks are uniquely determined by the sums of the columns
equalling $p$; thus the upper left entry is $p - (k+1)a - v^T \pmb 1$, the
next $k$ entries are all $p- (k+1)a - v^T \pmb 1 -1$, and the final entry
of the first row is $p- (k+1)a - v^T \pmb 1 +k$.
 
It is easy to check that $\Cal A (p,v,a)\Cal A (p',v', a') = \Cal A(pp',
p'v + v, p'a)$, and more generally if we inductively define $v^{n,0} =
v^{n+j}$ and $v^{n,i+1} = p_{n+j+1}\dots p_{n+j-i+1}v^{(n+j-i)} +
v^{n,i}$, and then define $v^{(n,j)} = v^{n,j}$ (now reflecting the
dependence on $j$), then
$$
\Cal A(p_{n+1}, v^n, a^n) \cdot \Cal A(p_{n+2}, v^{n+1}, a^{n+1})\cdot
\dots \cdot \Cal A(p_{n+j+1}, v^{n+j}, a^{n+j})  = \Cal A\(\prod_{i=1}^j
p_{n+1+i}, v^{(n,j)}, a^n \prod_{i=1}^j p_{n+1+i}\).
$$
(We have reversed the expected order, in order to facilitate transposition.)
The column sum is $\prod_{i=0}^j p_{n+1+i}$. Now assume that
$\brcs{v^i/p_{i+1}}$ is bounded. Then the normalized expression
$v^{(n,j)}/\prod_{i=0}^{j} p_{n+i+1}$ converges as $j \to \infty$. Denote
by $v^{\infty,n}$ the limit $\lim_{j\to\infty} v^{n,j}$; it is easy to see
that this is strictly positive. Then the infinite product converges (as $j
\to \infty$) to $\Cal A(1,v^{\infty,n},a^n/p_{n+1})$, and this is rank
one. Obviously, $\Cal A(p_n,v^{n-1}, a^{n-1})\Cal
A(1,v^{\infty,n},a^n/p_{n+1})  = \Cal A(1,v^{\infty,n-1},a^{n-1}/p_{n})$.
Hence if we take the leftmost (and strictly positive) column from $\Cal
A(1,v^{\infty,n},a^n/p_{n+1})$, and call it $V^{\infty,n}$, we have $\Cal
A(p_n,v^{n-1}, a^{n-1})V^{\infty,n} = V^{\infty,n-1}$.
 
Now assume the matrices are primitive, and relabel them $A_{n}^T$the
transposes; this yields an expression for $A_{n+j}^T A_{n+j-1}^T \dots
A_n^T$ as  $ \Cal A\(\prod_{i=1}^j p_{n+1+i}, v^{(n,j)}, a^n \prod_{i=1}^j
p_{n+1+i}\)^T$.  Set $W^{\infty,n} = (V^{\infty,n})^T $, so that
$W^{\infty,n}A_{n-1}^T = W^{\infty,n-1}$. Thus the family of strictly
positive rows, $\brcs{W^{\infty,n}}$ induces a trace on $G = \lim
A_n^T;\Z^{k+2}.\Z^{k+2}$ via
$\tau [x,k] = W^{\infty,k} x$. Since $G$ has unique trace, this must be it
(up to scalar multiple).
 
Now we convert this via the $z^{\perp}$ property to the trace computed
\wrt the alternative realization, ...
\endcomment 
\comment
\Lem Lemma \fivfou. Let $\Arrow t; G. U \subseteq \Q$ be a torsion-free abelian group together with a rational-valued group homomorphism. If $G$ is isomorphic to a limit of abelian groups, $\lim \Arrow M_n;\Z^{m}.\Z^m$ \st all $M_n$ have a common right eigenvector, $v$, \st $t([v,1]) \neq 0$, then the group extension $\ker t \to G \to U$ is nearly split.

\Pf Suppose that $M_n v = p_{n+1}v$ for integers $p_{n+1}$. We may assume that $v$ is unimodular, and thus $\Z^m$ admits an ordered $\Z$-basis whose first element is $v$. We may simultaneously rewrite all the $M_n$ in terms of this basis (this amounts to conjugating with the same element of $\GL{m,\Z}$, and so induces an isomorphism of abelian groups with map). Thus we may assume $M_n = \(\smallmatrix p_{n+1} & w_n \\ \pmb 0 & B_n \\ \endsmallmatrix\)$ where $w^n$ is a row of integers of size $m-1$, $\pmb 0$ is the column of size $m-1$, and $B_n $ is an $(m-1)\times (m-1)$ with integer entries, and moreover, the common right eigenvector is $v = (1,0,\dots, 0)^T$.

The group homomorphism $t$ is given by a sequence of rows, $r^j \in \Z^{1 \times m}$, $t [a, j] = r^j a/(c_1\dots c_{j})$ \st the following diagram commutes. 
$$\diagram
\Z^{m}&\rTo^{M_j} & \Z^{m} \\
\dTo^{r^j} &&\dTo^{r^_{j+1}} \\
\Z&\rTo^{\times c_{j+1}} & \Z\\
\enddiagram
$$
Here the sequence of  integers $(c_j)$  is \st $U \iso \lim \Arrow \times c_j ; \Z.\Z$; we obtain the equation $c_{j+1} r^j = r^{j+1}M_j$. 

Since $t([v,1]) \neq 0$ and thus $t([v,j]) = t([M_{j-1}M_2 M_1v,j]) \neq 0$, this forces the first coordinate of $r^j$, $r^j_1$, to be nonzero for all $j$. Since $M_j v = p_{j+1}v$, we have $ r^j v = r^{j+1}M_j v/c_{j+1} = p_{j+1}r^j v/c_{j+1}$. Since the first coordinates of the $r^j$ are not zero, this forces $p_{j+1} = c_{j+1}$. 

Set $H_0 $ to be the subgroup of $G$ generated by $\brcs{[v,j]}$. By construction, this is isomorphic to $\lim \Arrow \times p_{n+1};\Z.\Z$, which is isomorphic to $U$. Thus  $t$ maps $H_0$ into $U$, to which it is also isomorphic. If $V \subseteq U$ are isomorphic subgroups of the rationals, then $U/V$ is finite (in fact, there exists an integer $s$ \st $V = sU$). Hence the image of $H_0$ is of finite index, and thus $\ker t \oplus H_0$ is a finite index subgroup \st the restriction of $t$ splits.
\qed

If we drop the condition that the matrix sizes be bounded (that is, $m$ is replaced by $f(n) \to \infty$), it is unclear whether the the extension is still nearly split. However, if $U = \Z[\slfrac1p]$ (that is, $U$ is of the form $p^{\infty}$) for some prime $p$, then the extension is nearly split, for the trivial reason that any noncyclic subgroup of $\Z[\slfrac1p]$ is of finite index.

In any event, Lemma \fivfou\ tells us that if we want a bounded size \ers\ realization for the dimension group $G$ (with unique trace $\tau$, and that is rational-valued), then the corresponding extension, $\ker \tau \to G \to U$ must be nearly split. 
\endcomment

%%%%%%%%%%%

\comment

Now consider the matrices $A_n^T$. They have a common left eigenvector for
zero, $z = (-(k+1),1,\dots,1)$  coming from the relation $(k+1)a_0 = \sum_{i=1}^{k+1}$, which itself forces the same relation on the columns of $A_n$, thus on the rows of $A_n^T$ (of course, $z^T$ is the common right
eigenvector of zero for $A_n$), which is already unimodular. Since the rank of
$A_n^T$ is $k+1$ (we have listed the eigenvalues earlier), the left kernel of every single $A_n$ is just $z\Z$. Then
we see that an ordered  $\Z$-basis for $z^{\perp} = W$ is $((1,1,\dots,1)^T, (0,1,0,\dots,0,-1)^T, \dots, (0,0,\dots,-1,1)^T)$.
It is a routine calculation to see that the 
matrix of $A_n^T|W$ \wrt this basis for $W$ is  $C_n:=\(
\smallmatrix  p_{n+1} & (v^n)^T \\
0& \I_k \\ \endsmallmatrix\)$. The limit $J_0 \iso \lim \Arrow C_n;\Z^{k+1}.\Z^{k+1}$
is thus isomorphic to $G$ (as can be seen from the generating relations).
\endcomment

%%%%%%%%%%%

 %%%%%%%%%%%%
 
%%%%%%%%realizations
 \comment 
Suppose that $G$ is represented by the data $(p_{n+1}, v_n)$ (that is,
$g_n = p_{n+1} + v_n$), and in addition, $\sum 1/p_n = \infty$ (this is a
strong hypothesis in our context, and is  obviously highly unstable \wrt
telescoping). Then we may find a representation of $G$ by the data
$(p_{n+1}', v_n')$ in which $\| v_n'\|/p_{n+1}' \to 0$ (in general, $\sum
1/p_n' $ is finite).
\encomment
 
%%%%%%%%
If $\Z^k \to G \to U$ is nearly{} split, and the quotient is given by
$\tau$, then $G$ (as an abelian) is isomorphic to a direct limit group of
the form $\(\smallmatrix p_{n+1}  & w_n \\ 0 & \I_k \\ \endsmallmatrix\)$
where $w_n \in \Z^{1\times k}$ (and we can assume $\| w_n\|_{\infty} <
p_{n+1}$), $0$ is the column of size $k$; the direct limit of such a
sequence is nearly{} split.

definition is close enough
 
%%%%%%%%
 
{\it Sufficiently sparse\/} in the following conjecture is not defined,
but for example, if $\brcs{p_n}$ is an increasing family of primes with
$p_n/p_{n+1} \to 0$, then it should qualify as sufficiently sparse.
 
\Lem Conjecture. Suppose $k \geq 3$ and $\brcs{p_n}$ is sufficiently
sparse, and $G$ is the extension of $U = \lim \Arrow \times p_n; \Z.\Z $
with data given by $\brcs{v_n}$ (that is $g_n = p_{n+1}g_n + v_n$) where
$v_n/p_{n+1} \to (\slfrac 25, \slfrac 25, \dots, \slfrac 25)^T$, then the
dimension group $G$ with the strict ordering induced by the map $G \to U
\subset \R$ cannot be realized as a direct limit of simplicially ordered
groups (with $A_n$ nonnegative integer matrices), $\lim \Arrow A_n; \Z^d
.\Z^d$ for any integer $d$.
 
Some modification might be necessary in order to avoid easy manipulations,
e.g., instead require that the set of limit points of $\brcs{v_n/p_{n+1}}$
cover a large box in $\R^k$.
 
If this conjecture were true, it would yield counter-examples to the
following well known but nowhere explicitly-stated conjecture.
 
\itemitem{} {\it Conjecture.} If $G$ is a (simple) dimension group of
finite rank, then $G$ admits a Bratteli diagram of bounded width.
 
({\it Width\/} is sometimes (inaptly) called {\it depth.})
 
%%%%%%%%%%%%
\comment
Elementary realizations
 
Suppose in the matrix above, we can arrange that with $- \mu = \inf v_i
\wedge 0$ (so $\mu = 0$ if $v_i \geq 0$), we have
$$p > \cases
k\mu +k \| v\|_{\infty} + k +1 & \text{if $\mu > 0$} \\
(k+1) \| v\|_{\infty} + 2k+3 & \text{if $\mu = 0$}.\\
\endcases
$$
 
In the first case, we can set $a_i = \mu +1 > 0$.
If $\mu = - v_j$ for some $j$,   then $p - v^T\pmb 1 - (k+1)a_i = p -
\sum_{i \neq j} v_i - k\mu - k-1 \geq p - k \| v\|_{\infty} - k\mu -k -1> 
0$, and thus $p - v^T \pmb 1 - (k+1) a_i -1\geq 0$ for all $i$; this is
enough to ensure the matrix is nonnegative (which is good enough; a
telescoping will automatically be strictly  positive).
 
On the other hand, if $\mu = 0$, that is, $v_i \geq 0$, then set $a_i =1$;
then $p -v^T \pmb 1 - (k+1)a_i -1 \geq p - (k+1)\| v\|_{\infty} - 2k-3
\geq 0$.
 
The extra $k+1$, respectively $2k+3$, are not too serious limitations,
since by telescoping, we can make them very small \wrt $p_n$. However, the
condition $\| v^n\|_{\infty} < p_{n+1}/(k+1)$ (or its analogous condition
if $\mu^n > 0$) is severe.
 
We can also get sufficient realization conditions based on the spread of
the vector $v^n$; the {\it spread\/} of a vector $v$ is $s(v):= \max v_i -
\min v_i$. Suppose we have arranged (as we always can) that $ -p \geq \min
v_i < 0$,  but in addition that $(k-1)s(v)  +  k+2 \leq  p$ (this
obviously prevents $\min v_i = - p$). Then we set $a_j =1 -\min v_i$ for
all $j$, and we have
$$eqalign{
p - v^T\pmb 1 - (k+1)a_i -1&  = p - \sum v_i - k-2 + (k+1)\min v_i \cr
& = p - \sum_{j=1}^n (v_i - \min v_i) - k-2 \cr
& \geq p - (k-1) s(v) - k-2 \cr
& \geq 0.\cr
}$$
Hence we obtain an analogous (and slightly less restrictive)  condition
for realization, in terms of the spread of $v$. (Telescoping might remove
the extra $k+2$.)
 
For example, suppose $k =4$ and $v \sim (p/2, p/3, p/4, -p/2)$ (e.g., if
$v_1 $ is within $1$ of $p/2$, etc); then we can subtract $p$ from the
first three coordinates, creating $(-p/2, -2p/3,-3p/4,-p/2)$ which has
spread only $p/4$, so provided $p > 4k$ (since $k$ is fixed while $p$
varies and can be increased by telescoping---suitable telescoping
preserves or reduces the relative spread $s(v)/p$---this is relatively
innocuous), the vector can be realized.
\endcomment

\SecT 9 Infinite rank \ers

The following is a routine argument involving direct limits, but it allows
us to prove \sevtwo(b) via \sevtwo(a), as well as results on simultaneous \ers\ and \ecs\ realizations (\ecrs).
 
\Lem Lemma \ninone.  Let $G_n$ be a family of dimension groups and  let $\Arrow \phi^n; G_n.G_{n+1}$ be ordered group homomorphisms that send order unit to order units to order units. Let $G$ be the ordered group $\lim\Arrow \phi^n;
G_n. G_{n+1}$.
\item{(a)}Suppose  $H_n$ are  noncyclic rank one subgroups
of $G_n$ \st $H_n \cap G^{++} \neq 0$ for all $n$,  $\phi^n(H_n) \subseteq  H_{n+1}$, and each $G_n$ admits an \ers\
realization \wrt $H_n$.   Define $H = \lim 
A_n|H_n$. Then $G$ admits an \ers\ realization \wrt $H$ obtained from telescoping. 
\item{(b)}  Suppose that $(G,u) $ is given as the direct limit of $\Arrow \psi^j; G^j. G^{j+1}$ where each $G^j$ is a simple dimension group with unique trace, and each admits an \ecs\ realization. Then $G$ admits an \ecs\ realization \wrt its unique trace. 
\item{(c)} If each of the $G_n$ admit an \ers\ realization \wrt $H_n$ that is simultaneously \ecs, then $G$ admits an \ers\ realization \wrt $H$ that is also \ecs.
 
 \Pf (a) For each $n$, let $G_n \iso \lim_i \Arrow \phi_i^n; F_i^n. F_{i+n}^n$ (where $F_i = \Z^{f(i,n)}$ with the simplicial ordering) be an \ers\ realization of $G_n$ \wrt $H_n$. Let $\brcs{e_{ji}^n}$ be the standard basis of $F_i^n$. We may of course replace $\iso$ by equality. Since $\phi^n$ is positive, given $i$, there exists $ m \equiv m(i,n)$ \st for all $j \leq f(i,n)$, $\phi^n [e_{ji}^n] = [v^{j,i,m},m]$ where $v^{j,i,m}$ has all of its entries nonnegative. 

Since $\phi^n(H_n) \subseteq H_{n+1}$, we have $\phi^n [\pmb 1_{f(i,n)}, i] = p [\pmb 1_{f(l,n)}, l]$ for some integers $l > n$ and $p \geq 1$. Thus $\sum_{j=1}^{f(i,n)} [e_{ji}^n, n] = p[\pmb 1_{f(l,n+1)}]$. Hence there exists $m'\equiv m'(i,n)$ \st for all $j \leq f(i,n)$, 
$\phi^n [e_{ji}^n] = [w^{j,n+1}, m']$ where $w^{j,n+1} \geq 0$ and $\sum_{j=1}^{f(m',n+1)} w^{j,n+1}$ is a multiple of $\pmb 1_{f(m',n+1)}$. This means we can define a positive matrix $\Arrow C_i^n; F_i^n . F_{m'(i,n)}^{n+1}$ which sends $\pmb 1_{f(i,n)}$ to a multiple of $\pmb 1_{f(m',n+1)}$; in particular, $C_i^n$ has equal row sums.

Beginning with $i = 1$, we obtain a telescoping of the sequence for $G_2$ by composing and then relabelling $F^2_{m'(1,1)}$ as $F^2_1$ (and telescoping and relabelling the mappings), $F^2_{\max\brcs{m'(1,1), m'(2,1)}}$ as $F^2_2$, etc, so that now the matrices $C^1_i$ go straight down, that is map $F^1_i \to F^2_i$ (in the new notation). We may iterate this construction inductively. It is now straightforward that $G \iso \lim \Arrow C^{n,n}\circ \phi_n^n; F_n^n. F^{n+1}_{n+1}$, the order-isomorphism sending $H$ to the obvious limit of $H_n$, that is we have an \ers\ realization of $G$ \wrt $H$.

\noindent (b)  That $G$ has unique trace is trivial.  As in the preceding argument, we may telescope the various rows, and assume that each $\phi^j$ is implemented by {\it nonnegative\/}  matrices $\Arrow A_i^j; F_{i}^j. F_{i+1}^j$, and we may assume that no $A_i^j$ has any zero rows (in fact, since each $G^j$ is simple, it is easy to arrange that the matrices be strictly positive). 
 
 The element $u$ comes from some $G^j$, so we may normalize the unique trace $\tau_j$ of $G^j$ at its pre-image. Then $\tau_{j+1}\circ \phi^j = \tau^j $ (from uniqueness); hence, 
   $\tau^j[f,n] = \tau^{j+1}[A_n^j,n]$ (where $f \in F_n^j$). Since we have assumed each realization is \ecs, this says $\tau^j[e_i,n] = \tau^j[e_{i'},n]$ for all standard $\Z$-basis elements $e_i, e_{i'}$ of $F_n^j$, for all $j$, and in particular, $\tau^j [f,n] =\lambda_{j,n}\sum f_i$, where $f = \sum f_i e_i$ for some positive rational number $\lambda_{j,n}$
 
 Thus for each basis element $e_i$, we have (where $\pmb 1^T$ represents the row of the appropriate size consisting of ones)
 $$
 \lambda_{j+1,n} \pmb 1^T A_n^j e_i = \tau^{j+1}[A^j e_i, n] = \tau_j[e_i,n] = \lambda_{j,n}.
 $$
 Since the last term is independent of the choice of $i$, we have that all the $\pmb 1^T A_n^j e_i$ are the same, as $i$ varies. This means exactly that the column sums of $A_n^j$ are all equal. Now the   diagonal argument (as in (a)) can be applied.

(c) In the simultaneous case, we first ensure that the process in (a) is carried out, then apply the method of (b). 
 \qed
 \comment
\Pf For each $n$, let $G_n \iso \lim_i \Arrow \phi^n_i ; F^n_i.F^n_{i+1}$ (where
$F_i = \Z^{f(i,n)}$ with the simplicial ordering) be an \ers\ realization
of $G_n$ \wrt $H_n$. Let $\brcs{e^n_i}$ be the standard basis of $F^n_i$.
We may of course replace $\iso$ by equality. Since $\phi^n$ is positive, given
$i$, there exists $m\equiv m(i,n)$ \st $\phi^n [e^n_j, n] = [v^{j,n+1},
m]$ where $v^{j,n+1}$ has all of its entries nonnegative.
 
Since $\phi^n(H_n) \subseteq H_{n+1}$, we have $\phi^n [\One_{f(i,n)},i] =
p[\One_{f(l,n+1)}, l]$ for some $l  \in \Z$ and  integer $p \geq 1$. As
$[\sum_{j=1}^{f(n)} e^n_j,m] = p[\One_{f(l,n+1)},l]$. Hence there exists
$m'\equiv m'(i,n)$ \st (i) for all $j = 1, 2, \dots, f(i,n)$,
$\phi^n[e^n_j,n] = [w^{j, n+1},m']$ where $w^{j, n+1} \geq 0$ and
$\sum_{j= 1}^{f(m',n+1)} w^{j,n+1}$ is a multiple of $\One_{f(m',n+1)}$.
This means we can define a positive map $\Arrow C^n_i; F^n_i .
F^{n+1}_{m'(i,n)}$ which sends $\One_{f(i,n)}$ to a multiple of
$\One_{f(m',n+1)}$; in particular, $C^n_i$ has equal row sums.
 
Starting with $i =1$, we obtain a telescoping of the sequence for $G_2$,
by composing and then relabelling $F^2_{m'(1,1)}$ as $F^2_1$ (as well as
the telescoped mappings), $F^2_{\max\brcs{m'(1,1),m'(2,1)}}$ as $F^2_2$,
etc, so that now the $C^1_i$ go straight down, that is, map $F^1_i \to
F^2_i$ (the latter in the new notation). We may obviously iterate this
construction inductively. Now it is a triviality that $G \iso \lim \Arrow C^{n,n}
\phi^n_n; F^n_n . F^{n+1}_{n+1}$ (where $\phi^n_n$ is the newly relabelled
map; in the new relabelling, we are simply going down the diagonal) is an order isomorphism whose restriction to $H$ is the obvious one,
that is, we have an \ers\ realization of $G$ \wrt $H$.
\qed
\endcomment 
\Lem Lemma \nintwo. Let $G$ be a simple dimension group with unique trace $\tau$,
and let $H$ be a noncyclic rank one subgroup \st $\tau(H) \neq 0$ and
$G/H$ is torsion-free. Then we can write $G = \cup G_n$ where $G_n \subset
G_{n+1}$ are simple dimension groups with unique trace (in the relative
ordering), each of finite rank, and each containing $H$.
 
\Rmk Curiously, it is {\it not\/} true that every simple dimension group
can be written as an increasing union of finite rank simple subgroups [H5],
and moreover, this can occur in rather drastic ways, e.g., with just two
pure traces. It isn't even true that every simple dimension group can be
written as a direct limit of simple dimension groups each of whose pure
trace spaces is finite.
 
\Pf Consider the subgroup of the reals, $\tau(G)$; since this is
countable, and $\tau(H)$ is contained in it, we can find a countable set
of elements $\brcs{r_n}$ \st with $J_n:= \tau(H) + \sum_{i=1}^n r_i \Z
\subset \R$, we have $\tau(H) \subseteq J_n \subseteq J_{n+1}$ and $\cup
J_n = \tau(G)$. Select $g_i \in G$ \st $\tau(g_i) = r_i$.
 
Now $\ker \tau$ is a countable torsion-free abelian group (and nothing
else: every countable torsion-free abelian group can appear as a $\ker
\tau$); we may thus write it as an increasing union of free abelian groups
of finite rank (this is completely elementary: list the elements, then
take increasing finite subsets), say $\ker \tau = \cup T_n$, each $T_n$ of
finite rank.
 
Finally, set $G_n = T_n + H + \sum_{1 \leq i \leq n} g_i\Z$. Then $H
\subseteq G_n \subset G_{n+1} \subset \dots $. Since $\tau(H) \subseteq \tau(G_n)$, the
range of $\tau|G_n$ is dense, and it is immediate that with the relative
ordering inherited from $G$, $G_n$ is a simple dimension group with unique
trace $\tau$. Next let $G_0 = \cup G_n \subseteq  G$; we note that $\ker
\tau = \cup A_n \subset G_0$, and $\tau(G_0) = \tau(G) $ by construction.
Hence $G_0 = G$. Since each of $T_n$, $H$, and $ \sum_{1 \leq i \leq n}
g_i\Z$ is of finite rank, so is $G_n$.
\qed
 
\Lem Corollary \ninthr. [Theorem \sevtwo(b)] Let $G$ be a (countable) simple dimension
group with unique trace $\tau$, and let $H$ be a noncyclic rank one
subgroup of $G$ \st $\tau(H)\neq 0$ and $G/H$ is torsion-free. Then there
is an \ers\ realization of $G$ \wrt $H$.
 
\Pf By the preceding, we can write $G = \cup G_n$ with $H \subset G_n$,
where each $G$ is a simple dimension group with unique trace given by the
restriction of $\tau$. Since $G/H$ is torsion-free and $G_n/H \subseteq
G/H$, we have $G_n/H$ is torsion-free. Hence each $G_n$ admits an \ers\
realization \wrt $H$. The inclusion maps $G_n \to G_{n+1}$ send $H$ onto
$H$, and implement a realization of $G$ as a direct limit of the $G_n$s,
hence \nintwo\ applies.
\qed
%%%%%%%%%%

There are still questions about realizations that both \ers\ and \ecs\
(simultaneously; that is, the matrices have their row sums equal, and
their column sums equal). These will be addressed in the next two sections.
\comment
 If such a realization exists and is bounded
(that is, $\brcs{f(n)}$ is bounded), then the group extension $0 \to \ker
\tau \to G \to \tau(G)\to 0$ is nearly split (this amounts to $\tau(H)$
being of finite index in $\tau(G)$ here, although this is not the
definition), which is a very strong condition. In contrast, if $G$ is a
simple dimension group of rank $k+1$ with unique trace and that is
rational-valued, and there exists $H$ of rank one, etc, then $G$ admits
two realizations, one \ecs\ and the other \ers, both of size $k+2$, but
unless the extension is nearly split, will not admit a bounded realization
that is simultaneously \ers\ and \ecs. We do not know whether the nearly
split property is exactly what is needed to have a bounded simultaneously
\ers-  and \ecs-realization, although this is plausible. (It is true in
the stationary case, by an old theorem of Brian Marcus.)
\endcomment

\SecT 10 \ecrs\ and  nearly split extensions

A realization is \ecrs\ if it is simultaneously \ecs\ and \ers. The trace
induced by normalized multiples of the rows $\One_{f(n)}$ is automatically rational-valued, and
will be denoted $\tau$. If $G$ admits  an \ecrs\ realization wherein,
viewed as an  \ers\ realization, it  is \wrt $H$, then we shall write, {\it an
\ecrs\ realization \wrt $H$.} It is routine to see that if $G$ admits an
\ecrs\ realization \wrt $H$ that is of size $s$ (so all the matrices have
both $\One_{s}$ and $\One_{s}^T$ as their right and left Perron
eigenvectors respectively), then $|\tau(G)/\tau(H)|$ divides $s$: the
image of the trace on $G$ is $\cup (1/p_{j+1}\dots p_2) \Z$, while on the
image of the subgroup $\cup [\One_s, k]\Z$, it is $\cup (s/p_{j+1}\dots p_2)\Z$,
and the one by the other is a quotient of $\Z/s\Z$, hence has order
dividing $s$.
 
This puts a fairly stringent condition on the matrix sizes required for
bounded \ecrs\ realizations.  Of course, unbounded \ecrs\
realizations can be obtained as direct limits (obtained from unions) of
bounded ones (exactly as in the case of \ers\ realizations).

\Lem Lemma \fivtwo. Suppose $\Arrow t; G.U \subseteq \Q$ is obtained as the direct limit $ G \iso \lim \Arrow C_n ; \Z^s.\Z^s$ where the map $t$ is obtained from a common left eigenvector $w$ of all the $C_n$ (with corresponding eigenvalue $c_{n+1}$), via $t[a,n] = wa/c_1\dots c_n$. Suppose in addition, the $C_n$ have a common right eigenvector $v$ and  $wv \neq 0$. Then the extension $0 \to \ker t \to G \to U \to 0$ is nearly split. 

\Pf We note that the eigenvalue of $v$ for $C_n$ must be the same as that of $w$, $c_{n+1}$, since $vw \neq 0$. Set $H = \cup [v,n]\Z \subset G$. Then $\cup (c_2 \dots c_n)^{-1}vw\Z \subseteq H$, and this is obviously of finite index in $U = t(G) = \cup (c_1 \dots c_n)^{-1}\Z$. Thus $\ker t \oplus H$ is of finite index in $G$, so the extension is  nearly split. 
\qed

\comment
Suppose the abelian group and homomorphism $\Arrow t; J.\R$ are given as
$J = \lim \Arrow M_n; \Z^s.\Z^s$ where $M_n = \(\smallmatrix p_{n+1} & w^n
\\ 0 & B_n \\ \endsmallmatrix\)$ and $t[x,m] = \rho_m x$ (where $r_m \in
\R^{s-1}$ and $\rho_m = \rho_{m+1}M_m$). Let $H$ be the group $\cup
[e,m]\Z$ where $e = (1,0,\dots,0)^T$. Then $G/H$ is torsion-free;
moreover, $M_i$ restricted to $\ker \rho_i$ is just the restriction and
compression of $B_i$. The corresponding maps will be denoted $b_i$. We
also have that $e$ is a simultaneous right eigenvector of all the $M_n$,
hence we may consider the action of $M_n$ on the quotient group
$\Z^s/e\Z$; this amounts to a map which over the rationals would be
conjugate to $B_n$, but here is simply $\overline {B_n}$. We have the
following relations and maps
$$
\ker \rho_i   --b_i -- \ker \rho_{i+1} --   ...        = \ker t
=                                  =                                      
  da
\Z^2       --M_i--    \Z^s --                                    = J
darrow                   darrow                                    da
\Z^s/e\Z     --\overline{B_i} -- \Z^s/e\Z             = J/H
$$
This yields a one to one group homomorphism from $\ker t$ to $J/H$, and
moreover, the ranks are equal if $t$ is rational valued. Although both
$\ker t$ and $J/H$ are derived from the sequence $\brcs{B_i}$, in general,
they are not isomorphic, even in the stationary case, $B_i = B$. However,
if $|\det B_i| = 1$ for almost all $i$, then both $b_i$ and
$\overline{B_i}$ are invertible for almost all $i$, and thus both $\ker t$
and $J/H$ are free. Conversely, if either $\ker t$ or $J/H$ is free, then
$|\det B_i| = 1$ for almost all $i$, and so both groups are free.
\endcomment
 
When freeness occurs, then the index of the image of $\ker t$ in $J/H$ is
finite. In that case, $\ker t \oplus H$ has finite index in $J$, so that
the extension $\ker t \to J \to U$ is nearly split.  Thus we have the
following.
 
\Lem Lemma \tentwo. For any simple dimension group of finite rank with unique
trace, which is rational-valued, $\Arrow t; G.U \subseteq \Q$ \st $\ker t$
is free, and admits a bounded \ers\ realization, the extension $\ker t \to
G \to U$ is nearly split.
 
Although freeness of $J/H$ implies $J \to J/H$ splitsÑyielding a
group isomorphism $J \iso \Z^{s-1} \times H$Ñthis does not imply that $J
\to U$ splits. In the stationary and \ecrs\ example, $G = \lim \Arrow \(\smallmatrix 1 & 2 \\ 2 & 1 \\ \endsmallmatrix\);\Z^2.\Z^2$, we have  $\tau(G) := U= \Z[1/3]$ and $\tau(H) = 2\Z[1/3]$ (where $H = \cup_n [\One _2,n]\Z$), so the extension $0 \to \Z \to G \to U \to 0$ does not split, although  $G \iso \Z \oplus \Z[1/3]$ (as abelian groups) and
$t|H \neq 0$. As we will see in the next section, this is fairly typical.
 
Dropping the strong assumption that $\ker t$ or $J/H$ be free, a
sufficient condition for the extension $J \to U$ to be nearly split, that
is, $\ker t \oplus H$ be of finite index in $J$, is that the image of
$\ker t$ be of finite index in $J/H$.

 \comment
One peculiar condition under which we can prove this is if almost all the
$B_i$ (possibly after telescoping) satisfy $|\det B_i |$ is a power of the
same prime $r$ and there is just one nontrivial elementary divisor (that
is, the cokernel of $B_i$ is cyclic of order a power of $r$). This is
equivalent to  $J/H$ being an extension of $\Z^{s-2}$ by $\Z[r^{-1}]$. If
$A$ is a subgroup of this and of full rank and not free, then the index of
$A$ is finite (consider $A \cap Z^{s-2} \to A \to \pi(A) \subseteq
\Z[r^{-1}]$; as $A$ has full rank, its intersection with the free group
must be of rank $s-2$, hence the intersection has finite index in
$\Z^{s-2}$; if $\pi(A)$ were cyclic, then $A$ would be free, as the
extension would split, contradiction; hence $\pi(A)$ is not cyclic. Inside
$\Z[r^{-1}]$, all subgroups that are not cyclic are of finite index,
because $r$ is prime, so that $\Z[r^{-1}]/\pi(A)$ is finite, and it easily
follows that $A$ is of finite index.
So we have a very minor extension.
 
\Lem Lemma. Suppose that $\ker t$ or $J/H$ is an extension of a free group
by $\Z[r^{-1}]$ for some prime $r$. Then $\ker t \to J \to U$ is nearly
split.
 \endcomment
%%%%%%%%%%%

\comment
This yields more situations in which \ers\ forces the extension to be
nearly split. No. Thing about restriction being $b_i$ is wrong.
 
Suppose $M_n = \diag(p_{n+1,q_n})$ where $q_n = \oh{p_{n+1}}$ and for all
$i$ and $j$, $\gcd (p_i,q_j) = 1$. With $\rho_0 = (1,1)$, we have $\rho_i
= (1/p^{(i+1)})(p^{(i+1)}, q^{(i)}$, and $\rho_i^{\perp} = (q^{(i)},
-p^{(i+1)})\Z$ as is easy to see. Then the induced map to
$\rho_{i+1}^{\perp}$ given by $M_n$ is the identity! So the embedding
$\ker t \to J/H$ is of infinite index (here it is $\Z \subset U_1= \cup
\slfrac 1{q^{(i)}}\Z$. This yields a non-nearly split extension which has
a (group theoretic) realization with a common right eigenvector. Whether
there is a corresponding \ers\ realization for the simple dimension group
$G = U \oplus U_1$ with the strict order coming from the map $\Arrow t;
G.\R$ given by the sequence $(\rho_i)$ is still open.
 
But still have if the determinants of $B_i$ almost all $\pm 1$, then the
index is finite.
 
thing about single prime, etc, is wrong too.

For example, if the $B_i$ are scalar matrices (integer multiples of the
identity). Then the two abelian groups given by the respective sequences
$b_i$ and $\overline{B_i} $ are both isomorphic to $\Z^{s-1} \otimes U$
for $U$ a subgroup of the rationals (determined by the integers that
appear as the scalars), hence are isomorphic to each other, and thus the
index of the former in the latter is finite. This occurs when $\ker t$ is
rank one.
\endcomment

This places restrictions on  realizations of dimension groups by
commuting primitive matrices. For example, if $G$ is a simple dimension
group with unique trace, and it is realized by commuting nonnegative
matrices, we can telescope and kill off zero rows, and arrange that the
matrices additionally be primitive. Hence they will have common left and common
right Perron eigenvectors and the unique trace is determined from the left Perron eigenvector. Thus the corresponding extension $ 0 \to \ker \tau
\to G \to \tau(G) \to 0$ given by the image of the trace must be nearly
split. 

If for example, $\ker \tau = \Z^k$ and $\tau(G)=U $ has an
interesting supernatural number (e.g., every prime has multiplicity at
most one), then the set of isomorphism classes of nearly{} split  extensions within
the class of extensions of $\Z^k$ by $U$ is negligible. So most extensions
cannot be given by commuting families of matrices. Certainly a realization of bounded matrix size that is both \ers\ and \ecs\ qualifies for \fivtwo, as does a stationary dimension group.

In particular, Example \fivsev\ below is a simple dimension group with unique and rational-valued trace, of rank two, that cannot be realized by any sequence of (square) primitive matrices which have common right and common left eigenvectors; in particular, it cannot be realized by a bounded sequence of simultaneously \ers\ and \ecs\ primitive matrices. It {\it can\/} be realized by  a sequence of increasing size strictly positive rectangular matrices, 
$\Arrow M_n; \Z^{f(n)}. \Z^{f(n+1)}$, where $f(n) \to \infty$, and each $M_n$ is both \ers\ and \ecs, as we will see in the next section. It can be shown that for any such realization,  $f(n+1)/f(n)$ must be divisible by $3$ for infinitely many $n$. 

There is a trivial case in which the extension must be nearly split.

\Lem Lemma \fiveig. Suppose $\Arrow t; G.U$ is a finite rank  torsion-free group \st $U = \Z[\slfrac 1p]$ for some prime $p$. If $G$ contains a noncyclic rank one subgroup that is disjoint from $\ker t$, then the extension $0 \to \ker t \to G \to U \to 0$ is nearly split.

\Rmk Stationary examples show that such extensions need not be split.

\Pf Obviously, $t$ induces an embedding $H \to t(G)$. The extension is nearly split because every noncyclic subgroup of $\Z[\slfrac 1p]$ is of finite index!\qed

Let $q$ be a prime; all noncyclic subgroups of $U = \Z[\slfrac1p]$ are
therefore of finite index. Hence, if in the situation of the lemma  above,
$t(G)$ is isomorphic to $\Z[\slfrac1p]$ (and $H$ is not cyclic, which is
part of the hypotheses), then the corresponding extension is nearly split (Lemma \fiveig\ below).
 Another situation arises when the realization is by commuting matrices,
or more generally, when the implementing matrices have common left
eigenvector and common right eigenvector. In the situation arising from
positive matrices, these must be the Perron eigenvectors, hence correspond
to the same eigenvalue (for each $n$).

\Lem Example \fivsev. An example of a simple dimension group with unique, rational-valued trace, which is \ers-realizable, but for which the corresponding extension,
$0 \to \ker t \to G \to t(G) = U \to 0$, is not nearly split.

\noindent We construct a  simple example for which the subgroup $H \iso
\Z[\slfrac13]$ and $t(G) \iso \Z[\slfrac 16]$. In this case, the extension
$\Z \to G \to \Z[\slfrac 16]$ is not nearly split, but the corresponding
dimension group admits an \ers\ realization (of size three). It also
admits an \ecs\ realization, but cannot have a simultaneously \ecs\ and
\ers\ realization of bounded size (since that would imply common left and
common right eigenvectors, which entails nearly splitting, see Lemma \fivtwo\ below).
 
Construct an extension $0 \to \Z \to G_0 \to \Z[\slfrac12] \to 0$ for which
there are no $2$-divisible elements in $G_0$, equivalently, the
extension is not nearly split (as $2$ is a prime, we can get away with this). Let $\Arrow t_0;G.\Z[\slfrac12]$ denote the
map. We may regard $G_0$ as a subgroup of its divisible hull, which is of
course $\Q^2$; $t_0$ extends uniquely to a group homomorphism $\Arrow T;
\Q^2 . \Q$. Pick an element of $G$, $u \in t_0^{-1}(1)$, and form $G = G +
u\Z[\slfrac13]$ (inside $\Q^2$). Then $T$ restricts to a map,  called
$\Arrow t; G.\Q$, with values in $\Z [\slfrac12] + \Z[\slfrac 13] =
\Z[\slfrac 16]$.
 
 Now we show that $\ker t = \ker t_0 = \Z$, so that $\Z \to G \to
\Z[\slfrac16]$ is the corresponding extension. Elements of $G$ are of
the form $g= g_0 - um/3^k$ for  $g_0 \in G_0$, $m$ an integer, and $k$ a
nonnegative integer. If $t(g) = 0$, then $t(g_0) = m/3^k$; there exist
integers $l$ and nonnegative $j$ \st $t(g_0) = l/2^j$.
Hence $3^k l  = 2^j m$. This forces $3^k $ to divide $m$, so $g \in G_0$.
Hence $\ker t_0 = \ker t$.
 
Next, we show that $G$ contains no $2$-divisible elements. Select $g = g_0
+ um /3^k$ in $G$ with $g_0 \in G_0$ as in the previous paragraph. If $g$
were $2$-divisible, for all positive integers $l$, we could solve the
equations
$$
g_0 + u \frac m{3^k} = 2^l \(g_l + u \frac {m_l}{3^{k(l)}} \),
$$
where $g_l \in G_0$; we may assume that $m$ and $m_l$ are relatively prime
to $3$. Suppose for now that $k, k(l) > 0$. Fix $l$ and multiply by $3^k$.
This yields $3^k g_0 + mu = 3^k 2^l g_l + 2^l 3^{k - k(l)}m_l u$. Thus
$2^l m_l 3^{k - k(l)} u \in G_0$, so its value at $t_0$, $2^l m_l 3^{k -
k(l)} \in \Z[\slfrac12]$. Thus if $k(l) > k$, we must have $3$ dividing
$m_l$, a contradiction. Thus $k \geq k(l) $.
 
This yields $3^k(g_0 - 2^l g_l) = u(m - 2^l 3^{k-k(l)})$. Evaluating at
$t_0$, we obtain $3^k t_0( g_0 - 2^l g_l) = m - 2^l 3^{k-k(l)}$. If $k > 
k(l)$, then $3$ divides $m$ (as the values of $t_0$) lie in $\Z[\slfrac
12]$, again a contradiction. Hence $k = k(l)$, so that $m -2^l \equiv 0
\mod 3^k $. Since $m$ and $k$ are fixed, but the $\mod 3$ equivalence
classes of $2^l$ alternate between $1,2$, this is impossible.
 
Let us dispose of the remaining possibilities; first, if $k = 0$, we have
the equations
$g + mu = 2^l g_l + 2^l m_l 3^{-k(l)}u$, so $2^l m_l 3^{-k(l)}u \in G_0$;
evaluating at $t_0$, we obtain $2^l m_l 3^{-k(l)} \in \Z[\slfrac 12]$;
since $3$ does not divide $m_l$, we must have  $k(l) = 0$ for all $l$. But
then the element $g + mu = 2^l (g_l + m_l u)$
is $2$-divisible within $G_0$, a contradiction.
 
Next, if $k(l) = 0$ for one value of $l > 0$, then $um3^{-k} \in G_0$,
which forces $k = 0$ (evaluate at $t_0$ again), and we are in the
preceding case.
 
Thus $G$ contains no $2$-divisible subgroup. Since any subgroup of finite
index in $\Z[\slfrac 16] $ must be $2$-divisible, the extension cannot be
nearly split. On the other hand, the subgroup $H = u \Z[\slfrac13]$ of $G$
is $3$-divisible, so there is a group realization of the form described in
the lemma, that is, common right eigenvector, and a corresponding \ers\
realization for the dimension group. But there cannot be a dimension group
realization that is both \ers\ and \ecs\ simultaneously when the matrix
size is bounded. 
\qed

In terms of the $B_n$, a necessary condition for $G \to \tau(G)$ to split is that if $d_n = |
\det B_n|$, then $H$ is of finite index in $\cup (1/\prod_{i=1}^n d_i)\Z + t(H)$ (e.g., if $p_{n+1}$ are powers of the same prime $p$, this would force almost all the $d_n$ to be powers of $p$). But this is not sufficient. 

 \SecT 11 \ecrs 

Suppose that $(G,H)$ is a simple dimension group with noncyclic rank one subgroup \st $G/H$ is torsion free, and in addition that $G$ admits a unique trace $\tau$. Moreover, assume that $\tau(G):=U$ is a subgroup of the rationals, and $\tau(H) \neq 0$. These conditions (except the uniqueness of the trace) are necessary for an \ecrs\ realization of $G$ \wrt $H$. 

The converse is not quite true. We will show that if $U$ is $p$-divisible for some prime $p$ (that is, at least one prime has infinite multiplicity in the supernatural number of $U$), then the converse is true. However, in case $U$ is not $p$-divisible for any prime $p$, then an \ecrs\ realization exists \wrt $H$ exists if and only if $\rk G \leq |\tau(G)/\tau(H)|$. In this formulation, we allow $\infty$ as a value, and this corresponds to unbounded realizations. In the cases that $\rk G < \infty$, we have some control on the size of the realization. 

In particular, if $\tau(G)$ has no primes with infinite multiplicity, and $\rk G > 1$ (the case of $\rk G= 1$ is trivial), then the split case, $G =U \oplus \ker \tau$ with the strict ordering from the projection onto $U$, does not admit an \ecrs\ realization. In particular, if $\ker \tau$ is free of finite rank, by earlier results, then $G$ admits both an \ecs\ realization and an \ers\ realization \wrt $H$, of the same size, but no \ecrs\ realizations at all. 
 
We begin with the case that $|\tau(G)/\tau(H)| < \infty$. This of course implies that $G \to U$ is nearly split. For now, we also assume $G/H$ is free and finite rank.
Then we can write $G = H \oplus \Z^k$ (with $\tau(H) =
n \tau(G)$ for some integer $n$), but we must recall that $\ker \tau$ is not the copy of $\Z^k$ that appears as a direct summand. 

For a row or column $v$ consisting of integers, the {\it
content\/} of $v$, denoted $c(v)$ is the greatest common divisor of the
nonzero entries of $v$.
 
\Lem Lemma \eleone. Let $\lambda, p_{n+1 } > 1$ be positive integers \st for all  $n$,
$p_{n+1} \equiv 1 \mod \lambda$, and let $\rho \in \Z^k $ be a vector \st
$(c(\rho), \lambda) = 1$; set $q_n = p_2\cdot \dots \cdot p_n$. For each $n$, define $M_n= \( \smallmatrix
p_{n+1}& 0 \\ 0 & \I_k \endsmallmatrix\)$, and $r^n =
(\lambda/p_2\cdot\dots \cdot p_n, \rho)$ for $n > 1$, and $r_1 =
(\lambda,\rho) \in \Z^{1 \times (s+1)}$; with $ G = \lim M_n$ (as abelian groups) define $\Arrow t; G.\Q$ by $t[w,n] =r^n w$. Then there exist $v_n = \rho
(p_{n+1}-1)/\lambda, y_n = \rho(q_{n}-1)/\lambda \in \Z^{1\times s}$ \st
for all $n$, the following diagram commutes,
$$\diagram
{}&&\Z^{k+1}&\rTo^{M_1} & \Z^{k+1} & \rTo^{M_{2}} & \Z^{k+1} &{\cdots}&\Z^{k+1 } &\rTo{M_n}& \Z^{k+1}& \rTo^{\dots}\\
{} &\swarrow {r^1} &&& && && &&&\\
\Q \hskip-10pt& &\dTo^{\(\smallmatrix 1 & y_1 \\ 0 & \I_k \endsmallmatrix\)}&&\dTo^{\(\smallmatrix 1 & y_2 \\ 0 & \I_k \endsmallmatrix\)}&&\dTo{\(\smallmatrix 1 & y_3 \\ 0 & \I_k \endsmallmatrix\)}&&\dTo{\(\smallmatrix 1 & y_n \\ 0 & \I_k \endsmallmatrix\)}&&\dTo{\(\smallmatrix 1 & y_{n+1} \\ 0 & \I_k \endsmallmatrix\)}\\
{} & \hskip0pt\nwarrow r^1\\
&&\Z^{k+1} &\rTo^{\(\smallmatrix p_2 &v_1 \\ 0 & \I_k\\
\endsmallmatrix\)}& \Z^{k+1}&\rTo^{\(\smallmatrix p_3 &v_2 \\ 0 & \I_k\\
\endsmallmatrix\)}& \Z^{k+1}& \dots &\Z^{k+1} &\rTo^{\(\smallmatrix p_{n+1} &v_n \\ 0 & B_n \\
\endsmallmatrix\)}& \Z^{k+1}&\rTo^{\dots},\
\\
\enddiagram
$$
and in addition, $r^i M_i = r^{i-1}$ and $r^1$ is a common left eigenvector of all the matrices $\(\smallmatrix p_{n+1} &v_n \\ 0 & B_n \\
\endsmallmatrix\) $, with corresponding eigenvalue $p_{n+1}$. {\par}
If we set $H = \cup [\pmb 1_{k+1}, n] \Z$, then $t(H) = \lambda t(G)$.
 
\Pf  Set $y_1 = \pmb 0$. To $r^1$ as common left eigenvector, we must have $(\lambda/q_n)y_n + \rho = \rho/q_n$, that is
$y_n = \rho(q_n -1)/\lambda$; as $p_i \equiv 1 \mod \lambda$, $\lambda$
divides $q_n -1$, hence $y_n$ has only integer entries.
 
For the square to commute (now that we have define all the $y$s), it is
equivalent to $y_{n+1} = p_{n+1}y_n + v_n$, that is, we set $v_n = y_{n+1}
- p_{n+1}y_n = \((q_{n+1} - 1)/\lambda - p_{n+1}(q_{n}- 1)/\lambda\)\rho$,
and this simplifies to $v_n =\rho(p_{n+1}-1)/\lambda$.

 At the $n$th level, the trace is given by $(\lambda/q_n, \rho)$, so its
image is $q^{-1}_n \( \lambda \Z + q_n c(\rho)\Z\)$. Since $\gcd{(q_n,
\lambda)} = \gcd(c(\rho), \lambda) = 1$, we have $\gcd(\lambda, q_n) =
1$. Hence the range of the trace on the $n$th level is $q^{-1}_{n}\Z$, so
that $t(G) = \cup q^{-1}_{n}\Z$. On the other hand, $t[(1,0,\dots,0)^T,
n] = \lambda/q_n$. Hence the range of $t$ on $ H = \cup [(1,0,\dots,0)^T,
n] \Z$ is $\cup \lambda q^{-1}_n = \lambda t(G)$.
\qed
 
Under the assumptions of the lemma, set $v_0 = (1,1,0,\dots,0) \in \Z^k$.
There exists $E_0 \in \GL{k,\Z}$ \st $\rho E_0^{-1} =
c(\rho)(1,0,\dots,0)$. Then
$$
v_0(1-p_{n+1}) + \frac{p_{n+1}-1}{\lambda}\rho  E_0^{-1} =
\frac{p_{n+1}-1}{\lambda}\(c(\rho) + \lambda, \lambda,0,\dots,0 \).
$$
 
Now $\gcd\brcs{\lambda, c(\rho) + \lambda} = \gcd\brcs{\lambda, c(\rho)} =
1$. There thus exists $E_1 \in \GL{k,\Z}$ \st $\(c(\rho) + \lambda,
\lambda,0,\dots,0 \)E_1^{-1} = (1,1,\dots,1)$. Setting $v = v_0
E_0^{-1}E_1$ and $E = E_1 E_0$ then for all $n$, $v(1-p_{n+1}) + \rho
\(c(\rho) + \lambda, \lambda,0,\dots,0 \) E =
\frac{p_{n+1}-1}{\lambda}(1,1,\dots,1)$.
 
Let $u = (1,1,\dots,1) \in \Z^k$.
Now for any choice of integer $p$ (\st $\lambda $ divides $p-1$), we have
$$\eqalign{
\(\matrix 1 & v \\ 0 & E \\ \endmatrix \)\(\matrix p &
\frac{p-1}{\lambda}\rho \\ 0 & I \\ \endmatrix \)\(\matrix 1 & v \\ 0 & E
\\ \endmatrix \)^{-1} & =
\(\matrix p & -pvE^{-1} + (\frac{p-1}{\lambda}\rho + v) E^{-1} \\ 0 & E \\
\endmatrix \)  \cr
& = \(\matrix p & \frac{p-1}{\lambda}u \\ 0 & I \\ \endmatrix \).\cr
}$$
Hence, after conjugating every   $\(\matrix p_{n+1} & v_n \\ 0 & I \\
\endmatrix\)$ by the same matrix, we reduce to the case that the
transition matrices are $\(\smallmatrix  p_{n+1} &
\frac{p_{n+1}-1}{\lambda}u \\ 0 & \I\endsmallmatrix \)$, having $(\lambda
, \rho E^{-1})$ as common eigenvector, and since it is an eigenvector of
the matrices, it follows that $\rho E^{-1} = (1,1,\dots,1) \in \Z^k$, and
the trace on the group with homomorphism is given
by the suitably normalized eigenvector, $(\lambda,1,1,\dots,1)/q_n$ at the
$n$th level.
 
At this stage, we note that if $\lambda = k+1$, there is a simple
finishing argument. Add the first row of each matrix, that is, $(p_{n+1},
\frac{p_{n+1}-1}{\lambda}u)$ to all the other rows, and then subtract all
the columns from the first. This amounts to conjugating every one of the
matrices with same elementary matrix. The entries are suddenly strictly
positive, and since the inner product of the left and right unimodular
Perron eigenvectors is $\lambda = k+1$, and they consist  strictly
positive of strictly positive integers, they must all be exactly one. (We
will review this argument.)

We record the following elementary criterion.
 
\Lem Lemma \eletwo. Let $A$ be a primitive integer matrix of size $s$, whose
Perron eigenvalue is an integer, and let $V$ and $W$ be the corresponding
left and right Perron eigenvectors consisting of integers, \st $c(V) =
c(W) = 1$. If $VW = s$, then all row and column sums are equal.
 
\Pf The Perron eigenvectors consist of strictly positive real numbers, and
since they are all integers, each is at least one; as they are of size
$s$, the only way $VW $ is as small as $s$ is if every entry of each is
one. Hence the column and row sums are all equal.
\qed
 
In the case that $\lambda > k$, our strategy is to embroider a block of
$\lambda - k-1$ zero rows and corresponding nonzero columns (or zero columns
and nonzero rows) around each of our current matrices in such a way that
the resulting matrices still have common left and common right
eigenvectors corresponding to $p_{n+1}$, and such that their
unimodularized inner product (the $VW$ of Lemma \eletwo) is still $\lambda$. Then we conjugate all the
matrices (with the same matrix), so that as in the $\lambda = k+1$ case outlined above, the
resulting matrices are primitive.
 
The embroidered pieces actually vary in $n$ (in order to guarantee that
the eigenvectors do {\it not\/} vary in $n$), and must be carefully chosen.
 
If $\lambda < k+1$, we run into a technical difficulty when we try this, and indeed, an easy result shows that it is impossible to proceed. 
 
Suppose  $\lambda \geq  k+1$; this bifurcates
into $\lambda - (k-1) \leq k$ and $\lambda- (k-1) \geq k$ (for which the
treatments are similar).
 
We first justify the process of embroidering; this is elementary, and
completely derivative of symbolic dynamical techniques.
 
\Lem Lemma \elethr. Let $M_n:= (\smallmatrix A_n & B_n \\ C_n & D_n\\
\endsmallmatrix)$ be block partitions of $s \times s $ integer matrices
(with $A_n$ square of size $a$ and $D_n$ square of size $s-a$) of full
rank, $s$. Form the $S \times S$ matrices
$$
M_n' = \( \matrix A_n & B_n & 0\\
C_n & D_n & 0 \\
0 & X_n & 0 \\ \endmatrix\)
\qquad
M_n'' = \( \matrix A_n & B_n & Z_n\\
C_n & D_n & Y_n \\
0 & 0 & 0 \\ \endmatrix\)
$$
where  $X_n$ are  $(S-s)\times s$,  $Y_n$ are $s \times (S-s)$, and $Z_n$
are $a \times (S-s)$ integer
matrices. Then there are natural isomorphisms $G' = \lim M_n' \to G = \lim
M_n $ and $G \to G'' = \lim M_n''$, induced by $\Z^S \to \Z^s$ (projection
onto first $s$ coordinates) and the natural inclusion of $\Z^s$ in $\Z^S$.
 {\par}
\noindent Moreover, if $v= (\alpha,\beta)$ is a left eigenvector for $M_n$ (with
corresponding block decomposition), then $v' = (\alpha,\beta,\pmb 0)$ is a
left eigenvector for $M_n'$.
 
\Pf Let $V$ be the subgroup of $\Z^S$ with zeros in the top $s$ entries,
and let $\Arrow \phi; \Z^S.  \Z^s$ be the projection onto the top $s$
coordinates, so that $ VU = \ker \phi$. Then $\phi M_{n}' = M_n \phi$, so
$\phi$ induces a group homomorphism between the limit groups, which is
clearly onto. Since $\rk M_n = s$ and this is full, it easily follows that
$\rk M_n' = s$, hence $\rk G' \leq s$. As $G' \to G$ is onto, and the rank
of the latter ($s$) is at least as large as that of the former, the map
must be one to one.
 
Define $\Arrow \psi;\Z^s.\Z^{S}$ to be the inclusion (viewing $\Z^s$ as
the subgroup whose bottom $S-s$ entries are zero. Then it is trivial that
$M_n'' \psi = \psi M_n$, so $\psi$ induces a map $G'' \to G$, which is
obviously one to one. Since $M_n'' (\Z^{S}) \subset \phi(\Z^s)$, the map
is onto (in the direct limit).
 
The eigenvector property is trivial.
\qed

First consider the case $\lambda - k-1 \leq k$ (and $\lambda \geq k+1$).
Relabel our current matrices
$$M_n = \(\matrix p_{n+1} & \frac{p_{n+1}-1}{\lambda}u \\ 0 & \I \\
\endmatrix \);
$$
this has left eigenvector $(\lambda, u)$ and right eigenvector $(1,0,\dots,
0)^T$ for $p_{n+1}$ (recall $u = (1,1,\dots, 1)$). Here, $a = 1$ and $s =
k+1$. We set $X_n = \(\I_{\lambda -k-1}\  \pmb 0\)$ (the big zero is the
block of size $(\lambda-k-1)\times(k- \lambda)$, so $X_n$ is $(\lambda-k-1)
\times (k- (\lambda - k-1))$, so we have
$$
M_n' = \(\matrix
p_{n+1} & \frac{p_{n+1}-1}{\lambda}u & 0 & 0 & \dots & 0\\
0 & &0 & 0 & \dots & 0\\
0 & &0 & 0 & \dots & 0\\
\vdots& \I_k & \vdots &&&\vdots \\
0 & &0 & 0 & \dots & 0\\
\pmb 0 & \I_{\lambda-k-1} \  \pmb 0&& & \pmb 0\\
\endmatrix\)
$$
 
Now we perform the elementary column operations which simply add the first
$\lambda-k-1$ columns of the second block to their counterparts in the
third (so the columns get shifted to the right by $k$. The inverse
operation is to subtract the corresponding rows of the third block from
their counterparts in the second. The two operations together amount to
simultaneous  conjugation by the same  element of $\GL{\lambda,\Z}$, and
lead to the following matrices,
$$
 \(\matrix
p_{n+1} & \frac{p_{n+1}-1}{\lambda}(1,1,\dots,1) & 
\frac{p_{n+1}-1}{\lambda}(1,1,\dots,1)\\
\pmb 0_{1\times (\lambda-k-1)} &\pmb 0_{(\lambda-k-1)\times k} &\pmb
0_{(\lambda-k-1) \times(\lambda-k-1)}\\
\pmb 0_{k- (\lambda-k-1)} & \pmb 0_{(k- (\lambda-k-1))\times
(\lambda-k-1)} \  \I_{k- (\lambda-k-1)}&\pmb 0\\
\pmb 0 & \I_{\lambda-k-1} \  \pmb 0& \I_{\lambda-k-1} \\
\endmatrix\).
$$
Now we add the first row to each of the others, and correspondingly
subtract all the columns from the first; again, these are implemented
simultaneously in $n$ by a single product of elementary matrices, and
results in all the entries being nonnegative, and moreover, all the
matrices are primitive (since the first column and the first row are
strictly positive), and with the same zero pattern (so products will still
be primitive). Call these matrices $A_n$
 
The content one left and right eigenvectors of $M_n'$ for $p_{n+1}$ are
$V' = (\lambda, u, \pmb 0)$ and $W' = (1,0,\dots,0)^T$, hence their inner
product $V'W' = \lambda$. This is preserved by simultaneous conjugation;
as each $A_n$ is primitive of size $\lambda$, it follows from \eletwo\ above
that the left and right Perron eigenvectors of $A_n$ consist entirely of
ones, hence the column and row sums are equal. The simultaneous
conjugations obviously induce isomorphism of the groups with homomorphism
induced by the common left eigenvector, so we have a realization of $G$ as $
 \lim A_n$, which is \ecrs.
 
 In case $\lambda = k+1$, we skip the embroidery ($X_n$), and just proceed
via conjugations with elements of $\GL{k+1,\Z}$. If $\lambda-k-1 = k$, 
then there are no extra zero blocks, and the same process works. In this case, the realization is ultrasimplicial.
 
The process  for $\lambda \geq k+1$ and  $\lambda - k-1 \geq k$  (that is,
$\lambda \geq 2k+1$) is almost the same. We embroider the matrix with
$\lambda-k-1$ columns of zeros at the right (as we did before) and the
same number of rows at the bottom, and with $X_n$ being $\(\smallmatrix
\I_k \\ \pmb 0 \\ \endsmallmatrix \)$. Then we add the corresponding
columns to the third block and subtract the rows from the second
analogously with what we did before, and we can again just perform the
last operation, adding the first row to all the others and subtracting the
 columns from the first.

So far, we have the following.

\Lem Proposition \elefiv. Let $G$ be a simple dimension group of finite rank  containing a rank one noncyclic subgroup $H$ \st $G/H$ is free and $H \cap G^{++} \neq \emptyset$, and suppose $G$ has a unique trace $\tau$, and $\tau(G)$ is a rank one subgroup of $\Q$ whose supernatural number contains no primes of infinite multiplicity. Then $G$ admits an \ecrs\ realization  of size $\lambda:= |\tau(G)/\tau(H)| $ \wrt $H$ if $\lambda  \geq \rk G$. 

\Pf  If $\rk G = 1$, then there is almost nothing to do. Otherwise, $\lambda > 1$. For  subgroups $V \subset U$ of $\Q$, $|U/V| < \infty$ implies there exists $m$ \st $V = mU$, and if $U$ has no primes of infinite multiplicity, then $|U/mU| = m$. Set $U = \tau(G)$, and discard from the supernatural number all the primes (including multiples) that divide $\lambda$; the resulting subgroup $U_0$ is isomorphic to $U$, and correspondingly, $U_0/mU_0 $ is isomorphic to $U/mU$. Consider the set of primes (together with their multiplicities) dividing $U_0$; since they are all relatively prime to $\lambda$, we may telescope them to obtain sequence of positive integers representing $U_0$, $\brcs{p_{n+1}}_{n=1}^\infty$, \st $p_{n+1} \equiv 1 \mod \lambda$. 

Since $G/H$ is free, the extension  $\ker \tau \to G \to U$ is nearly split. Hence we can write $G = H \oplus \Z^k$, and the trace, given by the row $r_1$ at the first level, is of the form described in the top row of the statement of Lemma \eleone. The bottom row of the statement yields a representation of $G$ as a direct limit of abelian groups, with group homomorphism induced by the common left eigenvector, 
$[w,n] \mapsto r^1w/q_n$. 

The comment subsequent to the lemma allows us to assume that the realizing matrices are all in the form $\(\smallmatrix  p_{n+1} &
\frac{p_{n+1}-1}{\lambda}u \\ 0 & \I\endsmallmatrix \)$, having common left eigenvector $(\lambda,1,1,\dots,1)$. Now the embroidery process, together with \elethr, and subsequent simultaneous conjugation, gives an isomorphism of groups with group homomorphism to the direct limit of primitive matrices with equal row and column sums, as described above. 
\qed

Now we show that if $U$ has no primes of infinite multiplcity, then $|t(G)/t(H)| \geq \rk G$ is a necessary condition for $G$ to have a bounded \ecrs\ realization.

\Lem Lemma \elefou. Let $U$ be a noncyclic subgroup of rank one with no primes of infinite multiplicity. If $l$ is an integer exceeding $1$, then $U/lU \iso \Z/l\Z$. 

\Pf First, if $j > 1$, then $U \neq jU$, otherwise $\times j$ is a group automorphism of $U$, hence $\times 1/j$ is also an automorphism, and it easily follows that if $p$ is a prime dividing $j$, it must have infinite multiplicity in $U$. If follows that if $j$ properly divides $l$, then $jU \neq lU$. As every subgroup of finite index in $U$ is of the form $nU$ for some integer $n$, there is an obvious bijection between the intermediate subgroups $lU \subset U_0 \subset U$ and those of $l\Z \subset \Z_0 \subset \Z$, thus the map $\Z \to U \to U/lU$ has kernel $l\Z$, and is obviously onto.
\qed

Suppose $G$, with unique trace, has a realization as $\lim \Arrow A_n ;\Z^s.\Z^s$ which is \ecrs, where $H $ is identified with $\cup [\pmb 1_s, n]\Z$. Then the trace is given by the normalized constant row, and we see  immediately that $
\tau(H) = s \tau (G)$. Hence if $
\tau(G)$ has no primes with infinite multiplicity, we have $|\tau(G)/\tau(H)| = s \geq \rk G$. However, $\tau(G)/\tau(H)$ is an invariant of $(G,H)$, as $G$ has unique trace. 

\Lem Corollary \elesix. Suppose $G$ is a finite rank simple dimension group with unique trace $\tau$, \st $\tau(G)$ is a rank one noncyclic subgroup of $\R$ with no  prime divisors of infinite multiplicity. If $G$ admits an \ecrs\ representation \wrt $H$, then $\tau(G)/\tau(H)$ is finite and must be at least as large as $\rk G$.

\Pf Finiteness comes from the extension $G \to \tau(G)$ being nearly split (\tenone). The rest is from the comment just above.\qed

\Lem Theorem \elesev. Suppose that $G$ is a finite rank simple dimension group with unique trace $\tau$, having rational values, and $H$ is a rank one noncyclic subgroup \st $G/H$ is free and $\tau(G)$ is not $p$-divisible for any prime $p$. Then $G$ admits an \ecrs\ realization (\wrt $H$) if and only if $|\tau(G)/\tau(H)| \geq \rk G$. 

For example, if $G = U \oplus \Z^k$ where $U$ is an infinite multiplicity-free noncyclic subgroup of $\Q$, and we impose the strict ordering induced by the projection onto $U$, then the extension is split, and obviously $|\tau(G)/\tau(H)| = 1$; so $G$ admits an \ecrs\ realization (there is only one choice for $H$, namely $U$) if and only if $k =0$, and the latter is uninteresting. If instead, we impose as trace $\tau(u, v) = lu+v_1$ ($ru +$ the first entry of $v$), then $\tau(G) = U$, but $\tau(U) = lU$, so that $|\tau(G)/\tau(H)|=l$, then $G$ admits an \ecrs\ realization if and only if $l \geq k+1 = \rk G$.

%%%%some unbounded results
 
Now we assume $G$ simple dimension group with unique trace $\tau$,
$\tau(G)$ is rank one [and being dense, is noncyclic]
$H$ is a noncyclic rank one subgroup of $G$ \st $G/H$ is torsion-free, and
$\tau(H) \neq 0$.
We permit $\rk G$ and $\tau(G)/\tau(H)$ to be infinite. 
 
\Lem Theorem \eleeig. Suppose that $G$ is a simple dimension group with unique trace $\tau$, the value group of $\tau$ is $\tau(G) = U \subseteq \Q$, and $U$ has no primes of infinite multiplicity. Assume that $H$ is a rank one noncyclic subgroup of $G$ \st $\tau(H) \neq 0$ and $G/H$ is torsion-free. Then $G$ admits an \ecrs\ realization \wrt $H$ if either of the conditions below hold. {\par}
\item{(a)}
$|\tau(G)/\tau(H)| = \infty$, regardless of $\rk G$ (which can be infinite); in this case the realization must be unbounded. 
\item{(b)} $\infty > |\tau(G)/\tau(H)| \geq |\rk G|$, and in this case, the realization is bounded. 
 
\Pf First we note that if $U_0 \subset U$ are noncyclic rank one subgroups
of $\Q$, then there exists an infinite increasing chain of subgroups, $U_0
\subset U_1 \subset U_2 \subset \dots \subset U$ \st $ U = \cup U_i$ and
$|U/U_i| < \infty$. Applying this with $U_0 = \tau (H)$ and $U = \tau(G)$,
set $G_i^0 = \tau^{-1}(U_i)$. Moreover, $U_i/U_{0}$, being finite, is
cyclic. Hence there exists $g_i \in G_i$ \st $\tau(G_i) = \tau(H) + \tau(g_i) \Z$.
 
Since $\ker \tau$ is torsion free, we may find an increasing union of
finitely generated groups $F_1 \subseteq F_2 \subseteq F_3 \subseteq $ \st
$\ker \tau = \cup F_i$; by interposing as many equalities as we like, and
telescoping the $G_i$, we may assume $j + \rk F_j < |\tau(G_j)/\tau(H)|$.
 
Set $G_j = F_j + H + \sum_{l \leq j} g_l \Z$; then $G_j \subseteq G_{j+1}$
and $G = \cup G_j$. Moreover, $\rk G_j \leq \rk F_j + 1 + j \leq
|\tau(G_j)/\tau(H)|$. In addition, $G_j/H$ is finitely generated, and a
subgroup of $G/H$, hence is torsion-free, hence is free. Since $\tau(H)$
is dense in $\R$, $G_j$ with the relative ordering is a simple dimension
group with  unique trace,  the restriction of $\tau$. Thus $\ker \tau \cap
G_j \to G_j \to \tau(G_j)$ is nearly split, and the condition
$|\tau(G_j)/\tau(H)| \geq \rk G_j$ ensures that $G_j$ has a bounded \ecrs\
realization \wrt $H$.
 
Since $G$ is obviously the direct limit of $G_j$, by \ninone(c), $G$ has an
\ecrs\ realization \wrt $H$. In case (a), it must be unbounded (since bounded \ers\
realizations yield $|\tau(G)/\tau(H)| < \infty$). In case (b), the realization is obtained from telescoping a uniformly bounded family of realizations (using the method of \ninone(c)), so is bounded (or see the observation in the next paragraph). 
\qed
 
Now we have an elementary observation about unbounded \ecrs\ realizations, when
$\tau(G)$ has no infinite prime divisors. Suppose $G = \lim \Arrow A_n;
\Z^{f(n)}. \Z^{f(n+1)}$ is an \ecrs\ realization, with $\sup f(n) =
\infty$. The sequence of row vectors $(\pmb 1_{f(n)}^T/p_2\cdot \dots
\cdot p_n)$, where $\pmb 1_{f(n)}^TA_n = \pmb 1_{f(n+1)}^T p_{n+1}$
(defining $p_{n+1}$, the constant column sum of $A_n$), induces a trace
$\tau[y,n] = \pmb 1_{f(n)}y/p_2\cdot \dots \cdot p_n$. If we assume  that 
$G$ is simple with unique trace; necessarily, this is $\tau$. Then
$\tau(G)$ is $\cup_n 1/p_2\cdot \dots \cdot p_n$. With $H$ identified with
$\cup_n [\pmb 1_{f(n)}, n] \Z$, we see that $|\tau(G)/\tau(H)| \geq f(n)$
for all $n$ (this follows from no $p$-divisible subgroups for all primes
$p$).
 
Combining everything in sight, we have the following complete
characterization of \ecrs\ realizations when $\tau(G)$ has no
$p$-divisible elements for any prime $p$.
 
\Lem Theorem \elenin. Let $G$ be simple dimension group with unique trace, $\tau$.
Suppose that $\tau(G)$ is a subgroup of $\Q$ whose supernatural number has
no primes of infinite multiplicity. Let $H$ be a rank one noncyclic
subgroup of $G$ \st $\tau(H) \neq 0$ and $G/H$ is torsion-free. Then $G$
admits an \ecrs\ realization \wrt $H$ if and only if $\rk G \leq
|\tau(G)/\tau(H)|$; this includes the case that one or both of $\rk G$ and
 $|\tau(G)/\tau(H)|$ are infinite. Finally, every \ecrs\ realization is of size
$|\tau(G)/\tau(H)|$ (that is, unbounded if and only if $\tau(G)/\tau(H)$
is infinite).
 
When $G$ is $p$-divisible for some prime, the situation is  different; 
 no restriction on $\tau(G)/\tau(H)$ is required. 

%%%%%%%rest of ecrs situation
 
Now we asume that $t(G)$ is divisible by $p^{\infty}$ and to begin with, we also assume $G/H$ is free and$\lambda:= |t(G)/t(H)| < \infty$. If $p$ is any prime infinitely dividing $t(G)$, then it  also divides $t(H)$; hence $\gcd (\lambda,p) = 1$ for any prime $p$ dividing $H$ (which is isomorphic to $t(H)$). If $\lambda = 1$, we are in the split case, for which there is an interesting argument, obtaining a realization by commuting matrices.

Set $G = U \oplus \Z^k$ with the projection onto $U$ as the
unique trace---this is the split case---we show that   $G$ admits a bounded \ecrs\
realization (\wrt $H = U$, the only possible choice for $H$) under the assumption that  $U$ is $p$-divisible for some prime $p$.
 
Find a power, $q = p^a > k-1$. Then the matrix $M:=
\(\smallmatrix q & 0 \\ 0 & -I_k \\ \endsmallmatrix \)$  (note the
appearance of the {\it negative\/} of the identity matrix) satisfies all
the conditions of [BoH, xxx]. Hence there exists a primitive matrix $M'$
that is algebraically shift equivalent to $M$. By [M, xxx], there exists a
primitive matrix $A$ having equal row and column sums (so that $\pmb 1^T$
and $\pmb 1$ are respectively left and right Perron eigenvectors of $A$
for the eigenvalue $q$) shift equivalent to $M$. In particular, $A$ is
algebraically shift equivalent to $M$.
 
If the supernatural number has only finitely many other primes of
multiplicity at least one, then $U = \Z[1/p]$ and then $G$ admits a
stationary realization with $A_n = A$ (the argument to show this will be
included in what follows). Otherwise, we may telescope the other primes
(including their multiplicities), so that $U_0:= \lim \Arrow \times p_i ;
\Z . \Z$ ($p_i$ are products of the other primes) is relatively prime to
$p$ and $U = \Z[1/p]\otimes U_0$. Since $\gcd(p_i,p) =1$, so $\gcd(p_i,
p^a) =1$, hence by a further telescoping, we may also assume that $p_i
\equiv 1 \mod q = p^a$.
 
Now we use the following lemma to contort $A$.
 
\Lem Lemma \eleten. Let $m > 1$ be an integer, and suppose $l$ is a positive
integer with $l \equiv \pm 1 \mod m$. Then there exists $f \in \Z[x]^+$
(polynomials with nonnegative integer coefficients) \st $f(m) = l$ and
$f(-1) \in \brcs{\pm 1}$.
 
\Pf We find $f_0 \in \Z[x]^+$ \st $f_0 (m) = l$; then we modify it
inductively until $|f(-1)| = 1$. Expand $l = \sum_{i=0}^t a_i m^i$ with $0
\leq a_i < m$ as an $m$-adic expansion. Then set $f_0 = \sum a_i x^i$.
Obviously $f(m) = l$.
 
If $f_0(-1) > 1$, then $\sum a_i (-1)^i > 1$. If $a_0$ is the only
even-indexed coefficient that is strictly greater than zero, then $f(m)
\leq a_0 < m < l$, a contradiction. Hence there must exist $i = 2j$ \st
$a_i > 0$. Replace $a_i$ by $a_i-1$ and $a_{i-1}$ by $a_{i-1} + m$, to
create $f_1$. Then $f_1(m) - f_0(m) = -m^i + m^i = 0$, so $f_1(m) = l$,
and $f_1(-1) = f_0 (-1) - m-1$.
 
For any polynomial $g \in \Z[x]$, $g(m) \equiv g(-1) \mod {m+1}$. Hence $l
= f_0 (m) - f_0(-1)$ is a multiple of $m+1$; writing $l = km +1$ (if $l
\equiv 1 \mod {m}$, we have $km+1 - f_0(-1) = s(m+1)$, so $f_0(-1) = 1 + km
- s(m+1)$ and so $km \geq s(m+1)$. Also, $f_1(-1) = km -(s+1)(m+1) + 1$.
If this is negative, then $(s+1)m > km +1 \geq s(m+1)+ 1$, so $m > s$.
Also, $ s+1 > k \geq s + s/m$. This is impossible. Hence $f_{1}(-1) \geq
1$. If it equals $1$, we are done. If not, the process can be repeated,
each time reducing the value at $-1$ by $m+1$, and it must eventually hit
$1$. A similar process works if $l \equiv -1 \mod m$, except that the value
at $-1$ eventually hits $-1$.
 
If  $f(-1) < -1$, the process is similar, but easier (we do not have to
worry about large $a_0$). There must exist $i = 2j+1$ \st $a_i > 0$;
replace $a_i$ by $a_i -1$ and $a_{i-1}$ by $a_{i-1} + m$. The resulting
$f_1$ satisfies $f_1(m) = l$ and $f_1(-1) = f_0(-1) + m+1$. A similar
argument to that of the preceding allows us to conclude that $f_1(-1) < 0$
(if $l \equiv -1 \mod m$) or $f_1(-1) \leq 1$, whence either it is $\pm 1$,
or strictly less than $-1$, and the process can be iterated.
\qed
 
For each $p_i \equiv 1 \mod q$, there exists $f_i \in \Z[x]^+$ \st $f_i (q)
= p_i$ and $f_i (-1) = 1$. Set $A_n = Af_n(A)$; as $f_n$ has only
nonnegative coefficients, so does $A_n$; since each $A_n$ is a polynomial
in $A$, its large eigenvalue is $qf_n(q) = qp_n$, they commute with each
other, and have the same Perron eigenvectors, $\pmb 1^T$ and $\pmb 1$.
 
Suppose the matrix size of $A$ is $y$ (all we know is that $y \geq k+1$; otherwise, we have very little control over it).
 
Now form $M_n = M f_n (M) = \(\smallmatrix qp_n & 0\\ 0 & \I_k
\\\endsmallmatrix \)$. Suppose the algebraic shift equivalence between $M$
and $A$ is given by $X$ and $Y$; that is, $X M = A X$, $MY = Y A$, and $XY
= A^t$, $YX = M^t$ ($t$ is called the {\it lag\/}). Then for every nonzero
power of $A$, we have $XM^r = A^r X$ and $M^r Y = Y A^r$; hence for every
polynomial $g \in \Z[x]$ \st $g(0) = 0$, we have $Xg(M) = g(A) X$ and
$g(M)Y = Y g(A)$. Hence the map $\Arrow X; \Z^{k+1}. \Z^y$ induces a
group homomorphism $G = \lim M_n \to G' = \lim A_n$
by $[z,m] \mapsto [Xz,m]$, and similarly, $Y$ induces a group homomorphism
$G'   \to G$ via $[w,m] \mapsto [Yw,m]$. The products of the two group
homomorphisms are given by $\hat A$ and $\hat M$ respectively, both of
which are immediately seen to be group automorphisms of $G'$ and $G$
respectively. Hence the maps induced by $X$ and $Y$ are isomorphisms.
 
Moreover, they take the eigenspaces of nonzero eigenvalues  for $A$ to
those of $M$ (and vice versa), and in particular, they must send the
common eigenvectors for $q$, and thus send $(1,0,\dots)$ to $(1,1,\dots,1)$
and the same with the transposes. They thus induce an isomorphism of the
groups with group homomorphism. Moreover, it is easy (trivial) to see that
$G'$ has a unique trace (when given the direct limit ordering), so that
the group isomorphism is an order isomorphism from $G$ (with ordering
induced by the common left eigenvector of $M_n$ for $qp_n$) to $G'$ (with
direct limit ordering). It is trivial that $G$ is simply the split
extension. \qed

The upshot is a special case.

\Lem Lemma \elethi. Suppose $G = U \oplus \Z^k$ where $U$ is a rank one subgroup of $\Q$ that is $p$-divisible for some prime $p$, and the unique trace on $G$ is the projection onto $U$ (the split case). Then with $H = U$, there exists a bounded \ecrs\ realization of $G$ \wrt $H$ by commuting matrices. 

The matrices constructed in the other realizations need not commute. 
 A similar argument can be made to work in some non-split cases with a
prime having infinite multiplicity. We can of course extend this via the direct limit argument of \ninone(c). 

\comment
\Lem Proposition. Suppose $G = U \oplus K$ where $U$ is a rank one subgroup of $\Q$ that is $p$-divisible for some prime $p$, $K$ is an arbitrary countable torsion-free abelian group, and the unique trace on $G$ is the projection onto $U$ (the split case). Then with $H = U$, there exists a (possibly unbounded) \ecrs\ realization of $G$ \wrt $H$ by commuting matrices. 

\Pf The argument is familiar, and easier than others of its kind. We may write $K = \cup F_n$ where $F_n \subseteq F_{n+1}$ and each $F_n$ is free. Form $G_n = U \oplus F_n$; the relative ordering from $G$ is the ordering induced by the restriction of the projection, which is itself the projection onto $U$; hence $G_n$ satisfies the conditions of the previous lemma, hence admits an \ecrs\ realization \wrt $H$. Now xxx applies.
\qed
\endcomment

So we may assume that $\lambda > 1$. 
 
 If all primes infinitely divide $H$, then $H$ (and thus $t(G)$) are rational vector spaces, $\lambda = 1$, the extension splits (indeed, there is only one extension). If $H \simeq \Z[1/n]$ for some integer $n$, then the system is stationary, and the result follows from [M]. If there are only finitely many primes with finite multiplicity, we reduce to the last case immediately.
 
 Otherwise, there exist infinitely many primes each with finite and nonzero multiplicity, in addition to at least one prime $p$ with infinite multiplicity. Throwing away all the primes that divide $\lambda$ amounts to throwing away a finite set of primes with only finite multiplicity, hence doesn't change anything.
 
 There exists a power of $p$, $q = p^a$ \st $q \equiv 1 \mod \lambda$. We may also arrange, by taking a multiple of  $a$ if necessary, that $q\lambda  > k^2 + k$. We may telescope the other primes with their powers, so obtain $t(G) $ as $\Z[p^{-1}] \otimes \lim \Arrow p_i; \Z.\Z$ where $\gcd (p_i, \lambda) = \gcd (q, p_i) = 1$. Now we can implement the same isomorphism as in \eleone, with $qp_{i+1}$ replacing $p_{i+1}$ (and $(qp_{i+1}-1)/\lambda$ replacing $(p_{i+1}-1)/\lambda$), where we set $y_1 =\pmb 0$, $v_i = \rho\cdot (qp_{i+1}- 1)/\lambda $, and $y_n = \rho\cdot(qq_{n}-1)/\lambda$. This yields a group isomorphism to the limit group obtained as $\lim \( \smallmatrix  qp_{i+1} & w_i\\ \pmb 0 & \I_k\\\endsmallmatrix\)$, and the group homomorphism has been converted to left multiplication by the common eigenvector, $(\lambda, \rho)$.
 
 As in the previous case, we can simultaneously conjugate all the current matrices by $\(\smallmatrix 1 & v \\ \pmb 0 & E \\\endsmallmatrix\)$ where $v \in \Z^{1\times k}$ and $E \in \GL{k,\Z}$. This replaces the upper right entry by $\lambda^{-1}(q p_{i+1} - 1)\(\rho - \lambda v\)$. We could have previously conjugated the matrices with $(\smallmatrix 1 & \pmb 0 \\ \pmb 0 & J \\ \endsmallmatrix\ )$  where $J \in \GL{k,\Z}$, and so have assumed that $\rho = c(\rho)(1,0,\dots,0)$ (this applies to the left eigenvector as well). Now set $v = (0,1,0,\dots,0)$, so $\rho - \lambda v = (c(\rho), -\lambda)$; hence $c(\rho - \lambda v) = 1$. Thus there exists $E \in \GL{k,\Z}$ \st $(
 \rho - \lambda v)E^{-1} = (1,1,1,\dots,1) = \pmb 1_k^T$, which as before, we call $u$. 
 
 Hence we are in the situation wherein the matrices are of the form $\(\smallmatrix qp_{i+1} & (qp_{i+1} - 1)u/\lambda \\ \pmb 0 & \I_k\\ \endsmallmatrix\)$, their common eigenvector is $(\lambda , u)$ (for the eigenvalue $qp_{i+1}$), and the group homomorphism is obtained by left multiplication by suitable multiples of the eigenvector.

 Now we embroider $\lambda q - k-1$ rows and columns around the matrix; the only nonzero entries of the newly embroidered part occur in the top row, where we put $p_{i+1}\pmb 1^T$. This creates new matrices
 $$
 B_n = \(\matrix qp_{n+1} & \frac{qp_{n+1} - 1}{\lambda} \(1,1,\dots,1\)^T & p_{n+1}(1,1,\dots,1)\\
 \pmb 0 & \I_k & \pmb 0\\
 \pmb 0 & \pmb 0 & \pmb 0 \\
 \endmatrix \),
 $$
 which are of size $q\lambda$. Miraculously, $B_n$ have a common left eigenvector, $(\lambda q, q u, \lambda \pmb 1^T)$, where the third block is of size $\lambda q - k - 1$ (it is not in general true that embroidery where the right side is not zeroâ in this case, $Z_n$ in \elethr, will preserve the common left eigenvector property). 
 
 In order to perform the desired column and row operations, we need an estimate. From $q (k^2 - \lambda^2 - k + \lambda ) - k^2/p_{i+1} + k + 1 \geq k$ (easy), we see that 
 $$
 \frac{qp_{i+1} - 1}{\lambda} \geq k + \frac q{\lambda};
 $$
 hence there exists a multiple of $k$, $t = sk$, \st $q/\lambda \leq t/k \leq (qp_{i+1} - 1)/{\lambda}$. Now to each $B_n$, subtract $s$ times each of the first $k$ columns from their counterpart in the second block. The inequalities we just used are equivalent to the resulting top row consists of positive entries, and the sum of all but the first is less that $qp_{i+1}$. The inverse operation is to subtract the bottom rows from their counterparts, but this has no effect. This amounts to a conjugacy (which of course yields an isomorphism with group homomorphism), and now we simply add the top row to each of the others, and all the columns but the first from the first column. As before, the result is a primitive matrix, of size $\lambda q$; the inner product (as is easy to see) is the same as the size, so the matrix has equal row and column sums.
 \qed
 
 This yields the following rather surprising result.
 
 \Lem Proposition \eleele. Let $G$ be a simple dimension group of finite rank with unique trace $\tau$, \st $\tau(G)$ is $p$-divisible for at least one prime $p$. Suppose $H$ is a noncyclic rank one subgroup \st $\tau(H) \neq 0$ and $G/H$ is free. Then there exists an \ecrs\ realization of $G$ \wrt $H$. 
 
 And the direct limit argument, using \ninele\ and \eleeig, yields the definitive result.
 
 \Lem Theorem \eletwe. Let $G$ be a simple dimension group with unique trace $\tau$, \st $\tau(G)$ is $p$-divisible for at least one prime $p$. Suppose $H$ is a rank one noncyclic subgroup of $G$ \st $\tau(H) \neq 0$ and $G/H$ is torsion-free. Then $G$ admits an \ecrs\ realization \wrt $H$. 
 
 So we have a dichotomy: if $\tau(G)$ is not $p$-divisible for at least one prime $p$, then the condition $|\tau(G)/\tau(H)| \geq \rk G$ is necessary and sufficient (allowing $\infty$ as possible values); but if $\tau(G)$ is $p$-divisible for no primes $p$, there is no such constraint.

\SecT 12 Comments

Related to this is a  result of George Elliott [E], showing that the rank two split extension dimension group $\Z[\slfrac 12] \oplus \Z$ (with the strict ordering) cannot be realized as a limit of simplicial groups of rank two (that is, any direct limit realization requires almost all the free abelian groups to be of rank at leat three). In the latter, it was shown that this dimension group can be realized as a limit of rank three simplicial groups, and is stationary (via a size three primitive matrix algebraically shift equivalent to  $\diag (2,1)$). This is in fact what led me to think about using semigroups to obtain realization of the transfer matrices.

This paper was motivated by a question of Christian Skau: given the split extension $G = U \oplus \Z^k$, with $U \subseteq \Q$ and the projection onto $U$ yielding the ordering (so as to be a dimension group with unique trace), does it admit an \ers\ representation? (As we have seen, there is only one possible choice for the rank one noncyclic subgroup $H$ \st $G/H$ is torsion free, namely $U$ itself, so the choice of $H$ is unambiguous.) This appears as a special case, and the implementing matrices are the transposes of the matrices of the form $A$ appearing in section two (with parameters $ p = p_{n+1}$, once we ensure that $p_{n+1} > (k+1)^2$). I would like to thank Christian for his repeated insistence on solving this problem. 

Skau's question was motivated by questions concerning Tšplitz $\Z$-actions on Cantor sets (systems which admit factor maps onto odometers). A particular consequence of the results here is that among uniquely ergodic minimal actions of $\Z$ on Cantor sets, those that are strongly orbit equivalent to a Tšplitz, and those that are orbit equivalent, are characterized.

It has been known for over a decade that dimension groups which are rational vector spaces admit \ers\ realizations \wrt\ any dimension one subspace containing an order unit (this appears in [GJ]). The recipe is to begin with any realization of the dimension group, find an increasing sequence $h_n \Z \subset h_{n+1}\Z$ whose union is $H$ where $h_n$ is an order unit, telescope the realization, so that a cofinal collection of the $h_n$ appear, each at the $n$th level, say by a strictly positive vector $v_n$, apply the obvious diagonal matrix $\Delta_n \in \GL{f(n), \Q}$ so that $\Delta_n v_n$ is a  multiple of $\pmb 1$, replace the $n$th matrix $A_n$ by $\Delta_{n+1} A_n \Delta_n^{-1}$, then multiply each by a positive integer to ensure that the entries are all nonnegative integers. Since the dimension group $G$ satisfies $G \otimes \Q \iso G$, it follows immediately that the new improved direct limit yields $G$, and the elements of $H$ are implemented by constant vectors in the limit.

There is a substantial literature on realizing shift equivalence classes of integer matrices with nonnegative ones ([BoH] and the references there), corresponding to  stationary direct limits (that is, $G$ is a limit with the same matrix repeated, as an abelian group with real-valued homomorphism emanating from the largest eigenvalue and corresponding left eigenvector). Here we have a generally easier problem, since we are permitted to telescope matrices, something not allowed in the matrix realization problem. On the other hand, there are situations in  dimension group realization questions (such as $\tau(G)$ being a subgroup of the rationals with no primes of infinite multiplicity) which don't arise in the matrix realization case.

There are of course more questions to be answered. The most interesting would be to determine conditions on general dimension groups to admit an \ers\ realization. The general result characterizing \ecs\ realizations (\sixtwo(b)) suggests some sort of dual or transpose of good, this time applied to (rank one) subgroups. Thierry Giordano suggested that this might have a connection to recent work of Glasner and Host (extending earlier results of Giordano, Putnam, and Skau) on realizing dynamically the inclusion $H \subset G$ (for not necessarily rank one subgroups $H$).  

 \long\def\Rf[#1] #2, #3. #4\par%
{\vskip 2pt \itemitem{[#1]} #2, {\it #3,} #4\par\vskip2pt}

\SecT References

\Rf [BeH] S Bezuglyi \& D Handelman, Measures on Cantor sets\/{\rm:} the good, the ugly, the bad. Trans\ Amer Math Soc (to appear).

\Rf [BoH] M Boyle \& D Handelman, Algebraic shift equivalence and primitive matrices. Trans\ Amer  Math  Soc, 336 (1993),  121--149.

\Rf [EHS] {EG Effros, David Handelman, \& Chao-Liang Shen}, Dimension groups and their affine representations. Amer J Math 102 (1980) 385--407.

\Rf [E] George Elliott,  On totally ordered groups and K${}_0$. Proc Conf Ring Theory (Waterloo 1978),
Lecture Notes in Math, Springer, 734 (1979), 1--50.

\Rf [GJ] R Gjerde \& \O\ Johansen, Bratteli-Vershik models for Cantor minimal systems\/{\hskip .5pt\rm:} applications to Toeplitz flows. Ergodic theory and dynamical systems 20 (2000) 1687--1710.

\Rf [H1] D Handelman, Positive matrices and dimension groups affiliated to C*-algebras and topological Markov chains. J Operator Theory  6 (1981) 55--74. 

\Rf [H$1\slfrac12$] D Handelman, Ultrasimplicial dimension groups. Arch Math (1983) 109--115.

\Rf [H] D Handelman, Eventually positive matrices with rational eigenvectors. Ergodic theory and dynamical systems  7 (1987) 193--196. 

\Rf [H3] D Handelman, Matrices of positive polynomials. Electronic J Linear Algebra  19 (2009) 2--89.

\Rf [H4] D Handelman, Goodness of rational-valued traces on dimension groups. notes

\Rf [H5] D Handelman, Simple direct limits of dimension groups with finitely many pure traces. in preparation (despite the title, it is about simple dimension groups which are not direct limits of simple dimension groups each having finitely many pure traces). 

\Rf [M] Brian Marcus, Factors and extensions of full shifts. Monatshefte fŸr Mathematik
88 (1979) 239--247.

\Rf [R1] Norbert Riedel, Classification of dimension groups and iterating systems. Math Scand (1981). 

\Rf [R2] Norbert Riedel, A counterexample to the unimodular conjecture on finitely generated dimension groups. Proc Amer Math Soc (1981) 11--15.

\Rf [R] DJS Robinson, Groups---St Andrews 1981. London Mathematical Society Lecture Note
Series, \#71, 46--81.
\vskip 10pt

\noindent Mathematics Department, University of Ottawa, Ottawa ON  K1N 6N5, Canada; dehsg\@uottawa.ca

\end

\SecT endstuff

The following is well-known.
 
\Lem Lemma. Let $A$ be a torsion-free abelian group of finite rank. If $B$
is a subgroup of $A$ that contains an isomorphic copy of $A$, then $B$ is
of finite index in $A$.
 
\Pf Since rank is preserved by isomorphism, $\rk B = \rk A$; hence $A/B$
is torsion. Let $\Arrow \phi; A .A$ be the endomorphism that induces the
isomorphism $A \to \phi(A) \subseteq B$; again, $\rk \phi(A) = \rk A$, so
$\phi$ must be one to one. By the Cayley-Hamilton theorem (applied to the
divisible hull of $A$, a rational vector space), there exists an integer
$n$ \st $n\phi^k = \sum_{i=1}^{k-1} a_i \phi^i + n\det \phi \I $, where
each $a_i \in \Z$ as does $n \det \phi$. Thus $n\phi \cdot (\phi^{k-1} -
\sum_{i=1}^{k-1} a_i \phi^{i-1}) = n \det \phi\cdot  \I$, so that
$\phi^{k-1} - \sum_{i=1}^{k-1} a_i \phi^{i-1} = ((n\det \phi)/n)
\phi^{-1}$. Thus $n (\phi^{k-1} - \sum_{i=1}^{k-1} a_i \phi^{i-1}))  = (n
\det \phi)\phi^{-1}$. Setting $m = n\det \phi\in \Z$, we have $\psi:=m
\phi^{-1} $ is an integer polynomial in $\phi$, and thus belongs to the
endomorphism ring of $A$.
 
It suffices to assume $B = \phi(A) $. Applying $\psi$ to the inclusion
$\phi^2(A) \subset \phi(A)= B \subset A$, we have $m B \subset mA
\subseteq m\phi^{-1}(A) $. This last is contained in $A$ of course, so
applying $\phi$ to $mA \subseteq m \phi^{-1}A \subseteq A$, we have $mB
\subseteq mA \subseteq B$. Hence there is an onto map $A/mA \to A/B$, and
since $|A/mA| \leq m^{\rk A}$, $A/B$ is finite.
\qed
 
The converse (being of finite index implies it contains an isomorphic copy
of $A$, in fact an integer multiple of $A$) is true of course.

The matrices $A_n$ do not commute; however, if $p_n \neq p_m$, then $A_n
A_m - A_m A_n $ is rank one and square zero (note that while the left
Perron eigenvector for any $A_n$ is the constant row $\pmb 1$---of course,
the column sums are all $p_n$---the right eigenvector depends nontrivially
on  $p_n$; thus the matrices do not commute). However, the latter says
that they almost commute (in any event, we can always rearrange the $p_n$
and obtain the same $G$), so the order in which the $A_n$ occur is
irrelevant. Realizing $G$ with a sequence of (square) commuting matrices
can be done if $U$ is $m$-divisible for some prime $m$, and in fact, we can arrange that they are simultaneously \ecs\ and \ers, as in section 2. 
 
Now we form $H' = \lim \Arrow A_n^T; \Z^{k+2}.\Z^{k+2}$, given by the
sequence of transposes of the matrices. To analyze this, we have to
examine the right eigenspaces of each $A_n^T$ (or the left eigenspaces of
$A_n$). Again we drop the $n$. It is easy to verify that we have the
following linearly independent (over the rationals) set of eigenvectors,
with corresponding spectrum $p, 1,1,1,\dots, 1,1,1,0$
(we already knew $p$ and $0$, the latter from the rank being smaller than
the size; we suspected the remaining eigenvalues are all $1$ just from the
form of $G$).
 
Here are the  right eigenvectors of $A^T$, listed as columns corresponding
to the eigenvalues $p, 1,1,\dots, 1, 1, 0$ in that order (the eigenvalue $1$ hasalgebraic and
geometric multiplicity $k$). There is once again, a $k\times k$ identity
matrix occupying most of the space.
$$
\matrix 1 & 0 & 0 & \dots & 0 & 0&0 & -1\\
1 & 1 & 0 & \dots & \dots & 0 & 0 &p-1\\
1 & 0 & 1 & 0 & \dots & 0 & 0 & p-1\\
\vdots& 0 & 0 & \ddots & 0 &\dots & 0  & \vdots\\
\vdots & 0 & 0 & 0 & \ddots & 0  & 0 & \vdots\\
\vdots & 0 & 0 & 0 &0&\ddots & 0   & \vdots\\
1 & 0 & 0 & 0 & 0 & 0  & 1 & p-1\\
1 & -1 & -1 & -1 & \dots & -1 & -1 &-1 \\
\endmatrix
$$
We note first that the zero eigenspace depends on $p$, but the $p$- and
$1$-eigenspaces do not. Let $W$ be the (direct) sum of the $p$ and
$1$-eigenspaces; this is obviously invariant under every $A_{n}$ (corresponding to $p_n = p$).
Moreover, it is a direct summand of $\Z^{k+2}$ (there is an obvious $(k+1)
\times (k+1)$ minor in the first $k+1$ columns that has determinant one),
and it easily follows that $W$ is a maximal invariant subgroup of rank
less than $k+2$.
 
Now restore the subscripts (except for $W$, which is independent of the
choice of $p$). For each $n$, we have $A_n^T W \subset W$ (invariance) and
there exists $v \in \Z^{k+2}$ \st $W \oplus v \Z = \Z^{k+2}$ (the standard
copy) [note: we are {\it not\/} saying that we may choose $v$ to be
invariant as well---which would force it to be an eigenvector for
zero---this is impossible; in addition, since $W$ is a direct summand, we
may choose $v$ independently of $n$].
 
Index the first $k+1$ eigenvectors $z_1, \dots, z_{k+1}$. Then \wrt the
\Z-basis
$\brcs{z_i} \cup \brcs{v}$ of $\Z^{k+2}$, $A_n^T$ has the form
$$
U_n= \left[\matrix
p & 0 & \dots & \dots& 0 & * \\
0 & 1 & \dots & \dots & 0 &* \\
0 & 0 & 1 & \dots & 0 & * \\
\vdots   &\vdots& \ddots & \ddots & 0 & *\\
0 & 0 & 0 & \dots & 1 & * \\
0 & 0 & 0 & \dots & 0 & 0 \\
\endmatrix\right];
$$
that is, the upper $k+1$-square is $\diag (p, 1,1, \dots,1)$, the bottom
row is solid zeros and the rest of the final column can be anything (and
they depend on $n$). Each $U_n = B A_n B^{-1}$ for a $B$ in GL$(k+2,\Z)$,
the change of basis matrix, independently of $n$.
 
Now consider, merely as an abelian group, $J= \lim \Arrow A_n^T ;
\Z^{k+2}.\Z^{k+2}$. The conjugant $B$ implements an isomorphism between
this and $K= \lim \Arrow U_n ; \Z^{k+2}.\Z^{k+2}$. Now we show that this
is isomorphic $L = \lim \Arrow\Delta_n= \diag(p_n, 1,1,\dots,1,1);\Z^{k+1}.\Z^{k+1}$, the
latter obviously being $U \oplus \Z^k$ (but only as an abelian group).
 
Let $F_n$ denote the $n$th copy of $\Z^{k+2}$.  Let $f_j$ denote the standard basis of $F_n$ and also of $\Z^{k+1}$ (the latter with $1 \leq j \leq k+1$ only); now
define the map $\Arrow \phi_n;\Z^{k+1}.F^n$ given by $f_j \mapsto f_j$. It is trivial that the diagram
commutes
$$\diagram
\Z^{k+1}&\rTo^{\Delta_n} & \Z^{k+1} \\
\dTo^{\phi_n} &&\dTo^{\phi_{n+1}} \\
\Z^{k+2}&\rTo^{U_n} & \Z^{k+2}\\
\enddiagram
$$

The limit of the top row is obviously $\Z^k \oplus \lim \Arrow \times p_n; \Z.\Z = \Z^k \oplus U$. This yields a group homomorphism $\Arrow \Phi;\Z^k \oplus U.H' $. It is clearly one to one; moreover, since $\text{image}(\phi_{n+1}) = \text{image}(U_n)$, $\Phi$ is onto, and thus a group isomorphism. This yields a group isomorphism from between $H'$ and $G = U \oplus \Z^k$.

Now consider $H'$ as an ordered abelian group. The $1$-eigenspace  for each $A_n$ consists the columns whose top entry is $0$ and the sum of the entries is zero, and is spanned by the basis elements obtained previously. If $v$ is in the $1$-eigenspace of $A_n$, then as elements of $H'$, $[v,n] = [v,n+1]  = [v, n+2] = \dots$. Let $A$ be the subgroup of $H'$ consisting of all the elements of the form $[v,1]$ where $v$ runs over the $1$-eigenspace of $A_1$. Then it is easy to see that $A$ is a subgroup of $H$, free of rank $k$, and $A \cap (H')^+ = \brcs{0}$ (since $v$ remains unchanged after successive applications of $A_1$, $A_2$, etc, and the sum of entries of $v$ is always zero, so if it has a positive entry, it must also have a negative entry).

In addition, $H'$ has a unique trace%
\plainfootnote{$^{2}$}{\rm It is not true in general that if $H$ is a limit of square strictly
positive matrices (so is a simple dimension group) and $H$ has unique
trace, then the limit dimension group of their transposes need have 
unique trace (although it is simple). This is left as an exercise to the reader,
but with a hint: first do it for upper triangular $2 \times 2$ matrices where the number of traces---corresponding to certain
eigenvectors---can easily
be made to change by transposition, then perform a perturbation so the matrices are strictly positive.} ; this easily follows from $A_n/p_n $ converging to a rank one matrix. The basic result is, if for a sequence of nonnegative matrices, $A_1$, $A_2$, \dots, for all $N$, there exist a sequence of positive real numbers $s_{k,N}$ \st for all $N$, the sequence $A_{N+j}A_{N+j-1}\dots A_N/s_{j,N}$ converges to a rank $t$ real matrix, then the limiting dimension group has exactly $t$ pure traces.

Let $\tau$ denote the unique (up to scalar multiple) trace. Since $H'$ is simple and unperforated with unique trace, for all $w \in H'$, $\tau(w)> 0$ implies $w \in (H')^+$. Hence if $\tau(A) \neq \brcs{0}$, there would exist $a \in A$ \st $\tau(a) > 0$, hence $a \in (H')^+ \setminus \brcs{0}$, contradicting $A \cap (H')^+ = 0$. Thus $\tau(A) = 0$. Hence $\tau(H')$ is rank one.

The group isomorphism $G \to H'$ (invert the previous one) carries a copy of $U$ to another copy of $U$, which is disjoint from $A$ (obviously, since $A$ is free and $U$ is  a noncyclic subgroup of the rationals). Since $\tau$ cannot vanish on the image of $U$, and $U$ has unique unperforated ordering, up to scalar multiplication $\tau$ is plus or minus the identity on the image of $U$; the image of $U$ could possibly reverse the ordering, so we adjust the map $G \to H'$ by composing it with $h \to -h$ if necessary. The result is that the current group isomorphism $G \to H'$ is now positive (since positivity is determined completely by the $U$ component for both groups). Since  the group isomorphism restricts to an isomorphism on the kernels of the trace (from the $1$-eigenspace), the inverse is also order preserving. Hence the map $G \to H'$ is an isomorphism of ordered groups. In particular, all these dimension groups admit \ers\ realizations of size $k+1$ via the transposes of the displayed matrices.

The starting point for this investigation was a question by Christian Skau (arising from Toeplitz dynamical systems): is every dimension group of the form $G = U \oplus \Z^k$ (the split case) having as its unique trace the projection onto $U \subseteq \Q$ realizable by an \ers\ sequence? The affirmative answer is the second part of the statement of Proposition \twoone. 

One can also argue that $U$ is a fully invariant subgroup of $U \oplus \Z^k$ (since $U$ consists of the elements $g \in G$ for which the equations $g = n x_n$ can be solved for infinitely many $n$), so that $U$ in $G$ is automatically sent to the copy of $U$ in $H'$ by any group isomorphism. All the arguments here will be subsumed by the arguments in some subsequent sections.

My first attempt using this semigroup method was based on a simpler choice for the semigroup spanning $\Z^k$, that generated by $\brcs{\pm \epsilon_i}$; this yields a pleasant ECS realization of size $2k+1$ (the matrices have eigenvalues $p,1^k, 0^k$, but the Jordan form of the zero block was not clear); however, dealing with the transposes was problematic.
\end 

\comment
 \SecT 3 Slightly divisible subgroups of the rationals

If the subgroup $U$ has infinite multiplicity for at least one prime (that is, it is $p$-divisible for at least one prime $p$), then there is a realization of the partially ordered group $G = U \oplus \Z^k$ that is simultaneously \ers\ and \ecs. We cannot control the size of the matrices however, nor give them explicitly.

The supernatural number (defined up to aa equivalence) of $U$ is given by
a sequence $(n(q))$ where $q$ varies over the prime numbers and $n(q) \in
\brcs{\infty} \cup \Z^+$, and $\sum n(q) = \infty$. If $n(p) = \infty$ for some
prime
$p$, then we say that $H$ is $p$-divisible (in particular, this occurs if and only if
the PrŸfer group $\Z_{p^{\infty}}$ is embeddable in $U/\Z$ or equivalently, that $\Z[1/p] \subseteq U$).
 
A matrix is called {\it doubly stochastic\/} if all of its columns have
the same sum and all of its rows have the same sum (so if the matrix is
square, the column sums equal the row sums); we have generalized the usual
meaning of stochastic to permit the sums to be not $1$. A consequence in
the square case is that the column of the appropriate size consisting of
$1$s, $\One {}$ (or $\One n$  if we have to indicate the size), is a right
eigenvector for the row sum eigenvalue, and its transpose is a left
eigenvector for the same eigenvalue. Of course, if the matrix is primitive
(that is, all entries are nonnegative and some power is strictly
positive), then these are the right and left Perron eigenvectors
respectively.

\Lem Proposition \throne. Suppose that $G$ is a simple dimension group with unique trace, $t$, \st 
$U = t(G)$ is rational and $p$-divisible for some prime $p$, and moreover, $\Arrow t; G.U$
 splits, and $\ker t$ is free of rank $k$. Then $G$ admits a realization by primitive matrices, each of which is both \ecs\ and \ers.

Assume that  there exists a prime $p$ for which
$p^{\infty}$ divides the supernatural number of $U$.

Form the matrix $E = \diag (p,1,1,\dots, 1)$ (with $d$ ones). By [H], there exists a strictly positive matrix  $M$  (of size $m$, and we can
do this with $ m \leq d+1 + \ceil {\log_p (d+1)}$) \st $M$ is \ase\  to some
power of $E$. Hence if we construct the stationary diagram $G_M := \lim
\Arrow M;\Z^m . \Z^m$, then it follows easily that $G_M \iso Z[\slfrac1p] \oplus
\Z^d$ equipped with the strict ordering. By [M], there
exists a primitive double stochastic matrix $A$ shift equivalent  to $M$;
automatically, $G_A \iso G_M$.
 
We will find primitive matrices commuting with $A$ that will implement $G$.
For a stationary dimension group $G_A$ [H1; H3, sections 1--2, zero variables], an invariant of shift equivalence (and also of algebraic shift
equivalence for general direct limits with specified real homomorphisms)
was introduced, specifically, $\Endc (G_A)$ (in fact, it can easily be
extended to an isomorphism invariant for dimension groups, but in most
cases is trivial). There is a naturally occurring automorphism, $\hat A$
of $G_A$ ($[f,k] \mapsto [Af,k]$ with inverse $[f,k] \mapsto [f,k+1]$).
Form the collection of group endomorphisms $\Arrow \phi; G_A . G_A$ that
are $\Z[\hat A]$-module homomorphisms; this is $\Endc (G_A)$. Elements of
this ring are ordered in the natural way, and scale the (unique) trace.
 
It was shown that $\Endc (G_A)$ is obtained by forming the centralizer of
$A$ in its matrix ring over the integers, $C_{\Z} (A)$, factoring out
$\Set{B \in C_{\Z} (A)}{A^n B = \pmb 0}$ (where $n+1$ is at least as large as the matrix size), and then inverting the image of $A$ (either explicitly as a
subring of smaller size matrix ring, formally by taking the limit of
repeated multiplication by the image of $A$, $a$, or in the usual sense of
commutative algebra, since the image is not a zero divisor in the quotient
ring). The ordering is the natural one: positive elements are those of the
form $ba^k$ ($k \in -\N$) with $b$ represented by a matrix $B$ in $C_{\Z}
(A)$ for which there exists an integer $N$ \st $A^N B$ is primitive.

If we forget the ordering on $\Endc (G_A)$, and define it for any matrix $A$ as 
$$
\frac{C_{\Z}(A)}{\cup_{n=1}^{\infty}\Set{B \in \Endc (G_A)}{A^n B = \pmb 0} }[A^{-1}],
$$ (the $A^{-1}$ constructed formally from the direct limit) then it was shown in [H3, sections 1--2] to be an invariant of algebraic shift equivalence. In particular, as rings 
$\Endc (G_A) \iso \Endc G_E = \Z[\slfrac1p] \oplus \text{M}_k \Z$. Part of the algebraic shift equivalence data is a distinguished homomorphism $\Endc (G_A) \to \R$ given by the map to the reals induced by the Perron eigenvector (when $A$ is positive); as eigenvectors (for nonzero eigenvalues) are preserved by algebraic shift equivalence, this amounts to saying that the corresponding eigenvalue is distinguished in algebraic shift equivalences to non-positive matrices. 
 
Select a positive  integer $q$ relatively prime to $p$, and form $Q=
\diag (q,1,1,\dots,1)$. This is an element of the centralizer of $E$, and
moreover its value at the map to the reals given by the first coordinate
(corresponding to the trace when we perform algebraic shift equivalence of
$E$ to $A$) is $q$. Since $\Endc (G_A)$ is an invariant of algebraic shift equivalence
with a distinguished real homomorphism, there thus
exists an element $c$ of $\Endc (G_A)$ corresponding to $N$, and since $n > 0$, the element $c$ is positive. This means that there is a power of $A$ \st $(\hat A)^l c$ represents a positive matrix in the centralizer, call it $A_q$. The nonzero eigenvalues of $A_qN$ are $p^{lm}q, 1,1,\dots, 1$ arising from the multiplication by $A^l$.
 
Since $A_q$ is primitive and commutes with $A$ (we do not
claim that when $q\neq q'$,  $A_q A_{q'} = A_{q'}A_q$; this need not hold),
it easily follows that $A_q$ is doubly stochastic (with large eigenvalue
$p^{lm}r$), and moreover, the multiple eigenspace for $1$ is the same as
that for $A$.
 
 For each prime $q$ appearing in the supernatural number of $H$, take
$n(q)$ copies of $A_q$ (if $n(q) = \infty$, take countably infinitely
many), and also take countably infinitely many copies of $A$ itself. This
yields a countably infinite batch of doubly stochastic primitive
matrices, re-indexed (say) $A^{(i)}$ ($i=1,2, 3, \dots$, and form the
direct limit $G_0 = \lim_{i\in \N} \Arrow A^{(i)}; \Z^?.\Z^?$. The common
left Perron eigenvector ${\One ?}^T$ implements a map to the rational
subgroup with supernatural number given by $(q^n(q), p^{\infty})$.
 
The common eigenspace for $1$ together with the left eigenvector yields an
isomorphism as abelian groups with real map to their diagonal forms. Hence
as partially ordered abelian groups, we obtain an easy isomorphism $G_0
\iso U \oplus \Z^d$, and this preserves the unique trace, so is an order
isomorphism. This is more or less the  same algebraic shift equivalence argument used in the earlier method. This is simultaneously \ecs\ and \ers. We can make the implementing matrices commute with each other by observing that  $A^? A_q$ has its kernel equal to that of its powers (the algebraic and geometric multiplicities of zero as an eigenvalue are equal). 
 
Each of the $A^{(i)}$ is doubly stochastic, so satisfies both ERS and ECS.
\qed
\endcomment
\comment
First, we have an elementary lemma.

\Lem Lemma. Suppose $p > 2$ and $l$ is a positive integer \st $l \equiv
\pm 1\mod p-1$. Then there exists a polynomial $g \in \Z[x]^+$ \st $g(p-1) =
l$ and $g(-1) \in \brcs{\pm 1}$; moreover, we can choose $g \equiv g^{(l)}$ to have
degree at most $1 + \log_{p} l $.
 
\noindent [Obviously, this is a best possible result, since $g(p-1) \equiv g(-1) \mod p$.]
 
\Pf We will find a polynomial of suitable degree, with positive
coefficients, satisfying $g_0 (p-1) = l$ and modify it successively, at all stages
keeping the value at $p$-1, but reducing the value of $|g(-1)|$.  Let 
$$
S_l = \Set{v = (v_i)_{i=0}^{\infty}, v_i \in \Z^+,\ v_i = 0 \text{ aa}}{\sum (p-1)^i v_i  = l}
$$
To each element $v = (v_i)$ of $V$, associate $g_v \in \Z[x]^+$ via $g_v (x) =\sum v_i x^i$; then $g_v (p-1) = l$. We perform a sequence of operations on $S_l$  and their corresponding polynomials, terminating in the desired $g$. First, we note that $S_l$ is not empty: 
the $p-1$-ary expansion, $l = \sum_{i=0} \epsilon_i (k) (p-1)^i$ (here $\epsilon_i(l) \in \brcs{0,1,\dots,p-2}$) yields the sequence $w= (\epsilon_i (l)) \in S_l$. Let
$v =(v_0, v_1, \dots, v_r)$ be a sequence of nonnegative integers (truncated at $r$;  $v_{u} =0$ for all $u > r$), and
associate to $v$ the polynomial $g_v = \sum v_i x^i$.

Given $v$ in $S_l$, if  $g_v(-1) = \pm 1$, we are done. 

Suppose instead, (a) $g_v(-1) < -1$.
Then there must be $t$ \st $v_{2t+1} > 0$; perform the operation $v_{2t+1} \mapsto v'_{2t+1}= v_{2t+1} - 1$ and $v_{2t} \mapsto v'_{2t} = v_{2t} + p-1$ (leaving the other coordinates unchanged), yielding $v' \in S_l$ (equivalently, 
$1\cdot x^{2t+1}$  is replaced by $(p-1)x^{2t}$ in the polynomial); this
adds $p$ to the value at $-1$, that is, $g_{v'}(-1) = g_v(-1) + p$. Since $g_v(-1) \equiv g_v(p-1) = l \equiv \pm 1 \mod p$, it follows that $g_v(-1) \leq -p$, and thus either $g_v' (-1) = \pm1$, in which case we are done, or $g_{v'}(-1) < -1$, and we can iterate the process, which will obviously eventually terminate after at most $\ceil{|g_v(-1)|/p}$ iterations.
 
\noindent (b) $g_v(-1) > 1$. Then $\sum v_{2t} = \sum {v_{2j+1}} + np\pm 1$ for some positive integer $n$. 

If $v_{2j} > 0$ for some $j > 0$, then we perform the operation given by $v_{2j} \mapsto v'_{2j}= v_{2j} -1$ and $v_{2j-1} \mapsto v'_{2j-1} = v_{2j-1} + p-1$, creating the new $v'$, which is also obviously in $S_l$, and moreover $g_{v'} (-1) = g_v(-1) - p$, and either the latter is $\pm 1$ (terminating the algorithm), or we still have $g_{v'}(-1) > 1$, and we are in case (b) again.

If $v_{2j} = 0$ for all $j > 0$, then we have $v_0 = np \pm 1 + \sum v_{2j+1}$. We can terminate in one operation, via  $v_0 \mapsto v'_0=  v_0 - np$ (so $v'_0 \geq 0$) and $v'_1 = v_1 + n$. Then $v' \in S_l$, and moreover, $g_{v'}(-1) = g_v(-1) - lp - n  = \pm 1$. 

If we start with $w$ (arising from the $p-1$-expansion), then the degree of the corresponding polynomial is $\flo{\log_{p-1} k} +1$, and all iterations do not increase it (even the last subcase). 
\qed
 
Now suppose that $U$ is $m$-divisible, where $m$ is prime. Select a power of $m$, $p_0= m^d$ \st $p > k$. Form the diagonal matrix $(k+1) \times (k+1)$ matrix $\diag(p_0-1, -1,-1,\dots, -1)$. By [H2], all sufficiently high powers of this matrix are each algebraically shift equivalent to a primitive matrix; take an odd power (so that the nonzero non-Perron eigenvalues are still $-1$,  with multiplicity $k$. By [M], the latter is algebraically shift equivalent to a primitive matrix $Q_0$ whose right and left Perron eigenvectors are the column and row consisting entirely of $1$. In particular, the row and column sums of $Q_0$ are all the same, $(p_0-1)^d$ for some $d$. 

Now look at the remaining primes that divide $U$. If only finitely many remain and those have only finite multiplicity, $U \iso \Z[\slfrac 1m]$, and we can realize $U$ simply by repeating $Q$ infinitely often.

There are only finitely many primes that divide $p+1$. If their total multiplicity is finite, we can discard them. Otherwise, at least one has infinite multiplicity********

Otherwise, we may regroup the primes (repeated according to their multiplicity) so as to create a sequence of integers $p_1, p_2, \dots$ \st $p_i \equiv 1 \mod p$ for all $i$ and $U \iso \Z[\slfrac 1m]\otimes \lim \Arrow \times p_i;\Z.\Z$. Now choose $g_i \in \Z[x]^+$ \st $g_i(p) = p_i$ and $g_i(-1) = 1$, and then form $A_n = Q\cdot g_n(Q)$. Each $g_n(Q)$ has nonnegative entries (since $Q$ does and $g_n \in \Z[x]^+$, and so $A_n$ is), and we can even arrange that $A_n$ be strictly positive by replacing $Q$ in the product by a suitably large power of $Q$, if we wish. Of course, every $A_n$ has the same Perron eigenvectors as $Q$, so all the row and column sums of each $A_n$ are equal.

Now the same algebraic shift equivalence argument as in the earlier method shows (much more easily) that $G \iso \lim \Arrow A_n ;\Z^{?}.\Z^{?}$. This is simultaneously \ecs\ and \ers, and the transition matrices mutually commute. 
\endcomment

\comment
Suppose $(G, \tau)$ is a simple dimension group with unique
trace which can be realized as a limit $\Arrow A_i; \Z^d .\Z^d$ \st all
the $A_i$s have a common left Perron eigenvector, and a common right
eigenvector. Then $\tau(G):= U$ is a subgroup of the rationals (up to
scalar multiple), and the corresponding extension $0 \to \ker \tau \to G
\to U \to 0$ is nearly{} split.
 
\Pf Since $G$ is simple, we can telescope the sequence $\brcs{A_i}$, and
so assume that each $A_i$ is already primitive. Let $v$ and $w$ be a
choice for the  respective left and right common eigenvectors (they are
unique up to real scalar multiple), with eigenvalues $p_i$; set $q_n =
\prod_{i\leq n} p_i$. Then the map $G \to \R$ given by $[a,n] \mapsto
va/a$ (as usual, $\Z^d$ consists of columns, while $v$ is a row)  is a
trace, hence up to scalar multiple must be the only one. Hence (by
adjusting the scalar multiple for $\tau$), we may assume $\tau ([a,n]) =
va/q_n$. Hence $U = \tau(G) \subseteq \Q$.
 
Set $H_0$ to be the subgroup of the direct limit generated by
$\brcs{[w,n]}$. Since $A_{n+1} w = p_{n+1}w$,
it easily follows that $H_0 $ is an increasing union of cyclic subgroups,
$\cup_{n=1}^{\infty} [w,n]\Z$ is a rank one subgroup of $G$, and obviously
$H_0 \cap \ker \tau = \brcs{0}$. Moreover, the range of $\tau$ on $H_0$ is
$\cup (vw/q_{n+1})\Z = (vw)U$. Since both $v$ and $w$ are strictly
positive, $m:=vw$ is nonzero integer, and thus $\tau(H_0) = mU$. Thus
$\tau(H_0) $ is of finite index in $U$, so  the extension is nearly{} split.
\qed

\endcomment

\Lem Lemma. Let $p_{n+1} \uparrow \infty$, let $r^1 = (1, \rho^1)$ where
$\rho^1 \in \R^{1 \times k}$. Then there exist a sequence $\brcs{w^i}$,
with $w^i \in (\Z^k)^+$ and $\| w_i\|_{\infty} < p_{i+1}$ for all $i > 1$
together with group isomorphisms $\Arrow F_i; \Z^{k+1}. \Z^{k+1}$ \st the
following diagram 
$$\diagram
{}&&\Z^{k+1}&\rTo^{\(\smallmatrix p_2 & 0 \\ 0 & \I_k\\
\endsmallmatrix\)} & \Z^{k+1} & \rTo^{\(\smallmatrix p_3 & 0 \\ 0 & \I_k\\
\endsmallmatrix\)} & \Z^{k+1} &{\cdots}&\Z^{k+1 } &\rTo{M_n:= \(\smallmatrix p_{n+1} & 0 \\ 0 & \I_k\\
\endsmallmatrix\)}& \Z^{k+1}& \rTo^{\dots}\\
{} &\swarrow {r^1} &\dTo^{F_1}&& \dTo^{F_2}&&\dTo{F_3} &&\dTo{F_n}&&\dTo{F_{n+1}}&\\
\R \hskip-10pt& &&&&&\\
{} & \hskip20pt\nwarrow(1,w^{\infty})\\
&&\Z^{k+1} &\rTo^{\(\smallmatrix p_2 &w_1 \\ 0 & \I_k\\
\endsmallmatrix\)}& \Z^{k+1}&\rTo^{\(\smallmatrix p_3 &w_2 \\ 0 & \I_k\\
\endsmallmatrix\)}& \Z^{k+1}& \dots &\Z^{k+1} &\rTo^{\(\smallmatrix p_n &w_n \\ 0 & \I_k\\
\endsmallmatrix\)}& \Z^{k+1}&\rTo^{\dots}\\\
\\
\enddiagram
$$
commutes,  where
$$
w^{\infty}  = \frac{w^1}{p_2} + \frac1{p_2} \frac{w^2}{p_3} +
\frac{1}{p_2p_3}\frac{w^3}{p_4} + \dots.
$$
 
\Pf We will define $F_i = \(\smallmatrix 1 & y_i \\ 0 & C\\
\endsmallmatrix\)$ (where $y_i  \in \Z^{1\times k}$ and $C \in
\GL{k,\Z}$), together with $w_i \in \Z^{1 \times k}$ so that all the
properties hold. First, there exists $y_1 \in \Z^k$ \st all entries of $-
y_1 + \rho_1$ are strictly positive (rational) numbers. We can
simultaneously arrange that from the division algorithm, there exists $D
\in \GL{k,\Z}$ \st $z:= (-y_1 + \rho_1)D$ has all of its entries strictly
positive and between zero and $1$ (for example, if $\rho^1$ consists of
integers, the smallest we can realize is $(1,1,1\dots,1)$; but if $\rho^1$
has an entry which in lowest terms has denominator $d$, then we can
arrange that all the entries are at most $\slfrac 1d$); if $\rho^1$ has an irrational ratio, we can arrange that the entries are arbitrarily small. Set $C = D^{-1}
\in \GL {k,\Z}$.
 \comment
Let $t = z \pmb 1 > 0$, and let $\lambda$ be the positive root of
$\lambda(1-\lambda) = t$, that is, $\lambda = (1+ \sqrt{1 + 4t})/2$. Set
$z_1 = \lambda^{-1} z$
\endcomment
 
Now define $w_1$ to be the integer part of $p_2 z_1$ (that is, $p_2 z_1
\geq w_1$ coordinatewise, $\pmb 0 \leq w_1 \in \Z^{1\times k}$, and all
entries of $p_2z - w_1$ are strictly less than one). Now we proceed
inductively, setting $z_2 =p_2 z_1 - w_1 $, and set $w_2$ to be the
integer part of $p_3 z_2$; we construct $z_{i+1} = p_{i+1} z_i - w_i$, and
set $w_{i+1}$ to be the integer part of $p_{i+2}z_{i+1}$. This gives us a
sequence  $\brcs{w_i}$ of nonnegative elements of $\Z^{1 \times k}$. Now
we claim
$$
z_1 = \frac{w^1}{p_2} +\frac1{p_2} \frac{w^2}{p_3} + 
\frac{1}{p_2p_3}\frac{w^3}{p_4} + \dots.
$$
We observe that each $p_{i+1}z_i - w^i = z_{i+1}$ has all of its entries
in $[0,1)$, that is $\| z_i\|_{\infty} < 1$ for all $i \geq 2$ (and $\|
z_1\| \leq 1$). It is fairly obvious, but we claim that $Z_n:=z_1 -
\sum_{i=1}^n (p_{i}\dots p_2)^{-1}\frac{w_i}{p_{i+1}} $ goes to zero quickly.
 
We note that $z_1 - w_1/p_2 = z_2/p_2$, and $z_1- w_1/p_2 - w_2/p_2 p_3 =
z_2/p_2 - w_2/p_2p_3 = z_3/p_2p_3$, and inductively we see that $Z_n =
z_{n+1}/p_2p_3\dots p_{n+1}$. Hence $\| Z_n\| < 1/\prod_{i=1}^n p_{i+1}$.
Thus the series converges, and so $w^{\infty}$ exists and equals $z_1 $,
which equals $%\lambda  
(-y_1 + \rho_1)C^{-1}$.

\comment 
From $w^{\infty} = \lambda{-1}  (-y_1 + \rho_1)C^{-1}$, hitting it with
$\pmb 1$, we have $w^{\infty}\pmb 1 = \lambda^{-1} (-y_1 +
\rho_1)C^{-1})\pmb 1= \lambda^{-1} t$. From $\lambda^2 - \lambda -t = 0$,
we have $\lambda^{-1}t = \lambda - 1$. Hence $\lambda= 1 - w^{\infty}\pmb
1$. Thus $ (-y_1 + \rho_1)C^{-1} = w^{\infty}/(1-w^{\infty}\pmb 1)$.
\endcomment
 
Now define inductively $y_{i+1} = w^i C + p_{i+1} y_i$.
\qed

To deal with more general situations, e.g., rational-valued, but the kernel is not free, it would be nice if it were true that having an \ers\ realization is preserved under direct limits in which the maps were one to one; tha is,  $\Arrow \phi_j; G_j.G_{j+1}$ are one to one positive maps between simple dimension groups, and each one admits an \ers\ realization, then the direct limit, $G_0$, ought to have an \ers\ realization. But this fails if we restrict ourselves to bounded realizations, and with constant rank; take the extension we constructed before, $\Z \to G \to \Z[\slfrac16]$; this is a union of subgroups of the form $G_n:=\tau^{-1}(3^{-n}\Z + \Z[\slfrac12])$; each for each $G_n$, the restriction $\tau|G_n$ nearly splits, so each is \ers\ realizable by size three matrices, but we know that $G$ itself is not \ers-realizable at any size. (
It may, however, be \ers-realizable
with unbounded sizes.

However, if in the limit, the image of $\tau(G_n)$ is constant (that is, $\tau(G_n) = \tau(G_0)$ for all $n$; roughly speaking, only the group of infinitesimals is increasing as $n$ increases), it may still be possible to prove the union is \ers-realizable, and presumably at the same size. This would extend thm xxx to the same statement with kernel being free removed (we would require the extension to be nearly split, and it is unclear whether this is necessary, as we discussed earlier). To see this, write the kernel as a union of free groups all of the same rank, and throw in $U$ and possibly a few other elements to create a union of simple dimension groups for each of which the restriction of the trace map is nearly split.

This would avoid having to deal with the analogue of the arguments wherein $B_n$ replaces the identity matrix; there are simply not enough automorphisms of the resulting abelian group extensions. 

%%%%%useful lemmas